\documentclass[11pt]{amsart}
\usepackage{amscd,amsmath,amssymb,amsfonts,verbatim}
\usepackage[cmtip, all]{xy}
\usepackage{enumitem}
\usepackage{MnSymbol}
\usepackage[a4 paper, left=3cm, right=3cm, top = 2.5cm, bottom = 2.5cm ]{geometry}
\usepackage{comment}
\usepackage[colorlinks, 
linkcolor=black,anchorcolor=black, citecolor=black, 
linkcolor=black]{hyperref}

\newtheorem{thm}{Theorem}[section]
\newtheorem{prop}[thm]{Proposition}
\newtheorem{lem}[thm]{Lemma}
\newtheorem{cor}[thm]{Corollary}

\theoremstyle{definition}

\theoremstyle{definition}
\newtheorem{remk}[thm]{Remark}
\newtheorem{remks}[thm]{Remarks}

\newtheorem{defn}[thm]{Definition}

\newtheorem{exm}[thm]{Example}
\newtheorem{exms}[thm]{Examples}
\newtheorem{notat}[thm]{Notation}
\numberwithin{equation}{section}

{\hfill$\square$\end{defn}}
{\hfill$\square$\end{remk}}
{\hfill$\square$\end{remks}}
{\hfill$\square$\end{exm}}
{\hfill$\square$\end{exms}}
{\hfill$\square$\end{notat}}

\newcommand{\thmref}{Theorem~\ref}
\newcommand{\propref}{Proposition~\ref}
\newcommand{\corref}{Corollary~\ref}
\newcommand{\defref}{Definition~\ref}
\newcommand{\lemref}{Lemma~\ref}

\newcommand{\sA}{{\mathcal A}}
\newcommand{\sB}{{\mathcal B}}

\newcommand{\sD}{{\mathcal D}}
\newcommand{\sE}{{\mathcal E}}
\newcommand{\sF}{{\mathcal F}}
\newcommand{\sG}{{\mathcal G}}
\newcommand{\sH}{{\mathcal H}}

\newcommand{\sK}{{\mathcal K}}
\newcommand{\sL}{{\mathcal L}}
\newcommand{\sM}{{\mathcal M}}

\newcommand{\sO}{{\mathcal O}}

\newcommand{\D}{{\mathbb D}}

\newcommand{\F}{{\mathbb F}}

\renewcommand{\H}{{\mathbb H}}

\renewcommand{\L}{{\mathbb L}}

\newcommand{\N}{{\mathbb N}}

\newcommand{\Q}{{\mathbb Q}}

\newcommand{\Z}{{\mathbb Z}}

\newcommand{\fm}{{\mathfrak m}}

\newcommand{\ff}{{\mathfrak f}}
\newcommand{\fF}{{\mathfrak F}}

\newcommand{\Ker}{{\rm Ker}}
\newcommand{\gr}{{\rm gr}}

\newcommand{\surj}{\twoheadrightarrow}
\newcommand{\inj}{\hookrightarrow}
\newcommand{\red}{{\rm red}}

\newcommand{\Div}{{\rm Div}}
\newcommand{\Hom}{{\rm Hom}}

\newcommand{\Spec}{{\rm Spec \,}}

\newcommand{\Tr}{{\rm Tr}}

\newcommand{\sHom}{{\mathcal{H}{om}}}

\newcommand{\id}{{\operatorname{id}}}

\newcommand{\Sch}{{\operatorname{\mathbf{Sch}}}}

\newcommand{\<}{\langle}
\renewcommand{\>}{\rangle}

\newcommand{\Sm}{{\mathbf{Sm}}}

\newcommand{\can}{{\operatorname{\rm can}}}

\newcommand{\eff}{{\operatorname{\rm eff}}}

\newcommand{\et}{{\textnormal{\'et}}}

\newcommand{\ds}{{/\kern-3pt/}}

\renewcommand{\log}{{\operatorname{log}}}

\newcommand{\Tor}{{\operatorname{Tor}}}

\newcommand{\un}{\underline}
\newcommand{\ov}{\overline}

\renewcommand{\dim}{\text{\rm dim}}

\newcommand{\tuborg}{\left\{\begin{array}{ll}}
\newcommand{\sluttuborg}{\end{array}\right.}

\newcommand{\zar}{{\rm zar}}
\newcommand{\nis}{{\rm nis}}

\newcommand{\dlog}{{\rm dlog}}

\newcommand{\filt}{{\rm Fil}}
\newcommand{\rp}{{\rm RP}}

\newcommand{\wt}{\widetilde}

\newcommand{\cont}{{\rm cont}}

\newcommand{\coker}{{\rm Coker}}
\newcommand{\Fil}{{\rm fil}}

\newcommand{\n}{{\un{n}}}

\newcounter{elno}

\newcounter{elno-abc}   
\newenvironment{listabc}{
                         \begin{list}{\alph{elno-abc})
                                     }{\usecounter{elno-abc}}
                      }{
                         \end{list}}
\newcounter{elno-abc-prime}   

\begin{document}
\title[Duality for Hodge-Witt sheaves]{Duality for Hodge--Witt sheaves with modulus}
\author{Amalendu Krishna, Subhadip Majumder}
\address{Department of Mathematics, South Hall, Room 6607, University of California
  Santa Barbara, CA, 93106-3080, USA.}
\email{amalenduk@math.ucsb.edu}
\address{Laboratoire de Math{\'e}matiques d'Orsay, 307, rue Michel Magat,
Facult{\'e} des Sciences d'Orsay, Universit{\'e} Paris--Saclay,
F-91405 Orsay Cedex, France.} 
\email{subhamajumdar.sxc@gmail.com}




\maketitle
\vskip .4cm

\begin{quote}\emph{Abstract.}
  For an $F$-finite regular scheme of positive characteristic, we prove a
  modulus version of a duality theorem of Kato for the logarithmic Hodge-Witt sheaves. 
  As an application, we prove the modulus version of Milne's duality for the logarithmic
  Hodge--Witt cohomology over finite fields. We prove a similar duality theorem
  for semistable schemes over complete discrete valuation rings of positive
  characteristic with finite residue fields. The latter two results provide
  refinements of the duality theorems of Jannsen--Saito--Zhao and Zhao. 
  We also prove an extension of Ekedahl's duality for Hodge--Witt
  cohomology to modulus pairs over perfect fields.
\end{quote}

\setcounter{tocdepth}{1}
\tableofcontents

\vskip .4cm

\section{Introduction}\label{sec:Intro}
Duality theorems are fundamental tools to study cohomology theories in
arithmetic and algebraic geometry. The aim of this paper is to extend the duality
theorems of Kato, Ekedahl, and Milne for $p$-torsion sheaves in characteristic
$p>0$ to the setting of modulus pairs of schemes. In this section, we provide the
necessary background and state our main results. We begin with the coherent duality.

\subsection{Ekedahl's duality for modulus pairs}\label{sec:Ekedahl}
Let $k$ be a perfect field. Let $X$ be a regular and proper $k$-scheme
which is connected of dimension $d$. Let $W_m\Omega^\bullet_X$ denote the
de Rham--Witt complex of $X$ (cf. \cite{Illusie}). 
Let $i, q \ge 0$ and $m \ge 1$ be integers. Then Ekedahl's 
duality theorem \cite[Thm.~II.2.2]{Ekedahl} for the Hodge--Witt cohomology states that
there is a canonical isomorphism of finitely generated $W_m(k)$-modules
\begin{equation}\label{eqn:Ekedahl-00}
H^i_\zar(X, W_m\Omega^q_X) \xrightarrow{\cong}
\Hom_{W_m(k)}(H^{d-i}_\zar(X, W_m\Omega^{d-q}_X), W_m(k)).
\end{equation}
Ekedahl's duality theorem has several applications to the duality for the $p$-primary
{\'e}tale cohomology of smooth varieties over fields of characteristic $p > 0$.
(e.g., see \cite{Kato-duality-1}, \cite{Kato-duality-2}).

We now let $X$ be a connected, regular and proper $k$-scheme of dimension $d$,
and let $E \subset X$ be a simple normal crossing divisor with complement $U$.
We let $D \subset X$ be a Weil divisor supported on $E$ and let
$\Fil_D W_m\Omega^\bullet_U$ denote the filtered de Rham--Witt complex introduced in
\cite{KM-3}. The latter represents the de Rham--Witt complex of the modulus pair
$(X,D)$. The first result of this paper is the following extension of the
above theorem of Ekedahl to the Hodge--Witt cohomology of the modulus pair
$(X,D)$.

\begin{thm}\label{thm:Main-10}
The product operation of
$W_m\Omega^\bullet_U$ induces a canonical isomorphism of finitely generated
$W_m(k)$-modules
\[
H^i_\zar(X, \Fil_D W_m\Omega^q_U) \xrightarrow{\cong}
\Hom_{W_m(k)}(H^{d-i}_\zar(X, \Fil_{-D-E} W_m\Omega^{d-q}_U), W_m(k)).
\]
\end{thm}

\vskip.2cm

\begin{remk}\label{remk:Main-10-0}
 \thmref{thm:Main-10} recovers ~\eqref{eqn:Ekedahl-00} if $E = \emptyset$
  while it recovers the duality theorems of Nakkajima \cite[Thm.~5.3(1)]{Nakkajima}
  and Hyodo \cite[(3.3.1)]{Hyodo} when $E \neq \emptyset$ but $D = 0$ (note that
  $\Fil_0 W_m\Omega^q_U \neq W_m\Omega^q_X$ if $E \neq \emptyset$).
  \end{remk}

\subsection{Kato's duality for modulus pairs}\label{sec:Kato}
Let $p > 0$ be a prime and let $X$ be a connected and regular $F$-finite $\F_p$-scheme
having $p$-rank $d \ge 0$ (this means that $[k(X): k(X)^p] = p^d$).
Let $X_\rp$ denote the
relatively perfect site of $X$ introduced in \cite{Kato-duality-1}. This is a
Grothendieck site whose objects are schemes which are relatively perfect over $X$ and
whose covering maps are the usual {\'e}tale coverings. For $m \ge 1$, let
$\Lambda_m = {\Z}/{p^m}$. We let $D_m(X_\rp)$ denote the derived category of sheaves of
$\Lambda_m$-modules on $X_\rp$.

For any Zariski sheaf of $W_m\sO_X$-modules
$M$, we let $M^\rp$ denote the extension of $M$ to $X_\rp$, defined to be the
sheafification of the presheaf $(f \colon X' \to X) \mapsto (W_mf)^*(M)$.
Let $D^c_{m}(X_\rp)$ denote the full triangulated subcategory of $D_m(X_\rp)$
generated by
$M^{\rp}$, where $M$ is a coherent sheaf of $\sO_X$-modules on the Zariski site of
$X$. In \cite[Thm.~4.3]{Kato-duality-1}, Kato defined a dualizing functor
$\D^0_{X,m} \colon D^c_{m}(X_\rp) \to D^c_{m}(X_\rp)$ (cf. \S~\ref{sec:DCC})
and showed that this is an anti-equivalence of categories.

Suppose now that we are in the setting of \thmref{thm:Main-10}.
Assume additionally that $D \ge 0$ and let $D' = D+E$.
For $q \ge 0$, we let $q' = d-q$. We let
\[
W_m\sF^{q,\rp}_{X|D} =
\left[Z_1\Fil_{D} W_m\Omega^{q, \rp}_{U} \xrightarrow{1-C}
  \Fil_{D} W_m\Omega^{q, \rp}_{U}\right];
\]
\[
W_m\sG^{q,\rp}_{X|-D'} = \left[\Fil_{-D'} W_m\Omega^{q,\rp}_{U} \xrightarrow{1-\ov{F}} 
\frac{\Fil_{-D'} W_m\Omega^{q,\rp}_{U}}{dV^{m-1}(\Fil_{-D'}\Omega^{q-1, \rp}_{U})}\right];
\]
and $W_m\Omega^{q,\rp}_{X|D',\log} = \Ker(1-\ov{F})$. We refer to
\S~\ref{sec:DRP-site} for details about these complexes.
We prove the following modulus version of Kato's duality theorem from
\cite{Kato-duality-1}.

\begin{thm}\label{thm:Main-1}
  We have the following.
  \begin{enumerate}
  \item
    $W_m\sF^{q,\rp}_{X|D}$ and $W_m\Omega^{q,\rp}_{X|D',\log}$ are objects of
    $D^c_{m}(X_\rp)$.
    \item
There are canonical isomorphisms
  \[
  W_m\Omega^{q',\rp}_{X|D',\log} \xrightarrow{\cong} \D^0_{X,m}(W_m\sF^{q,\rp}_{X|D}); \ \
  W_m\sF^{q,\rp}_{X|D} \xrightarrow{\cong} \D^0_{X,m}(W_m\Omega^{q',\rp}_{X|D',\log})
  \]
  in $D^c_{m}(X_\rp)$.
  \end{enumerate}
\end{thm}
\vskip .3cm

Combining \thmref{thm:Main-1} with \cite[Thm.~5.1]{Kato-duality-1}, we get
the following generalization of Milne's theorem \cite[Thm.~5.2]{Milne-ENS}.

\begin{thm}\label{thm:Main-3}
  Assume in \thmref{thm:Main-1} that $X$ is a complete variety over a perfect
  field $k$ and let $\pi \colon X \to \Spec(k)$ be the structure map. Then there are
  canonical isomorphisms
  \[
  {\bf R}\pi_*(W_m\Omega^{q',\rp}_{X|D',\log}) \xrightarrow{\cong}
  {\bf R}\sHom_{\Lambda_m}({\bf R}\pi_*(W_m\sF^{q,\rp}_{X|D}), \Lambda_m);
  \]
  \[
 {\bf R}\pi_*(W_m\sF^{q,\rp}_{X|D})  \xrightarrow{\cong}
 {\bf R}\sHom_{\Lambda_m}({\bf R}\pi_*(W_m\Omega^{q',\rp}_{X|D',\log}), \Lambda_m)
 \]
 in $D_{m}(k_{\rp})$.
\end{thm}

 \vskip .2cm

\subsection{Milne's duality for modulus pairs}\label{sec:Milne}
Assume now that $k$ is a finite field and $X$ is a complete $k$-variety of dimension
$d$ over $k$ which is regular and connected. Let $W_m\Omega^\bullet_{X, \log}$ denote
the logarithmic Hodge--Witt {\'e}tale sheaves {\`a} la Bloch, Deligne and Illusie. 
The duality theorem of Milne (cf. \cite{Milne-Zeta}) for the $p$-primary
{\'e}tale cohomology says that there is a perfect pairing of finite groups
  \begin{equation}\label{eqn:Milne-D}
   H^i_\et(X, W_m\Omega^q_{X, \log}) \times H^{d+1-i}_\et(X, W_m\Omega^{d-q}_{X, \log}) \to
     \Lambda_m.
  \end{equation}

  In order to extend Milne's result to modulus pairs, assume
  in the setting of \thmref{thm:Main-10} that $D \ge 0$ such that $D_\red = E$.
For $q, m \ge 1$, we let
$W_m\Omega^q_{X|D,\log}$ be the image of the map of {\'e}tale sheaves
$\dlog \colon \Ker(\sK^M_{1,X} \surj \sK^M_{1,D}) \otimes j_* (\sK^M_{q-1,U})
\to W_m\Omega^q_X$ (cf. \cite[Defn.~1.1.1, Prop.~1.1.3]{JSZ}), where
$\sK^M_{\star, (-)}$ is the improved Milnor $K$-theory sheaf {\`a} la Gabber and Kerz
(cf. \cite{Kerz-JAG}).
We let $W_m\Omega^0_{X|D,\log} = \Ker({\sK^M_{0,X}}/{p^m} \to {\sK^M_{0,D}}/{p^m})$.

For $q \ge 0$, we let $W_m\sF^{q, \bullet}_{X|D}$ denote the two term complex
$[Z_1\Fil_D W_m\Omega^q_U \xrightarrow{1-C} \Fil_D W_m\Omega^q_U]$ of {\'e}tale sheaves
on $X$. We refer to \cite[Cor.~6.17]{KM-3} for the proof that
$W_m\sF^{q, \bullet}_{X|D}$ is well-defined.
As an application of \thmref{thm:Main-1}, we obtain the following modulus version of
Milne's duality theorem.

\begin{thm}\label{thm:Main-2}
Given integers $q, i \ge 0, m \ge 1$, there is a perfect
  pairing of finite groups
 \[
 \H^i_\et(X, W_m\sF^{q, \bullet}_{X|D}) \times
 H^{d+1-i}_\et(X, W_m\Omega^{d-q}_{X|D+E, \log}) \to \Lambda_m.
   \]
\end{thm}

\begin{remks}\label{remk:GK-duality}
  (1) It is easy to check that \thmref{thm:Main-2} coincides with ~\eqref{eqn:Milne-D}
  when $E = \emptyset$.
  
  (2) In \cite{JSZ}, Jannsen--Saito--Zhao proved a duality theorem for certain
  $p$-primary {\'e}tale cohomology on $U$. Their duality theorem  was
  extended by Gupta--Krishna in \cite{GK-Duality-1} to the case when 
$E$ is not necessarily a simple normal crossing divisor.
We show that the pairing of \thmref{thm:Main-2} is compatible with respect to
the partial order on the set of effective divisors $\Div^{\eff}_E(X)$.
In particular, we can pass to the limit over $\Div^{\eff}_E(X)$, and the resulting
perfect pairing recovers the duality theorems of \cite{GK-Duality-1} and
\cite{JSZ}.

(3) In \cite{Ren-Rulling}, Ren--R{\"u}lling introduced a different version of
Hodge--Witt cohomology with modulus and established duality theorems analogous to
Theorems~\ref{thm:Main-10}, \ref{thm:Main-1} and \ref{thm:Main-2}. However, their
cohomology groups are unrelated to those appearing in the above results, except when
$D=0$ (cf. \corref{cor:EK-Main-4}).

Another difference between \thmref{thm:Main-10} and the corresponding result in
\cite{Ren-Rulling} (see its Corollary 9.4) is that the latter requires $D$ to be
effective, whereas no such assumption is needed for the former.
\end{remks}

\subsection{Milne's duality for log schemes}\label{sec:Log-sch}
Let $(X,E)$ be as in \thmref{thm:Main-2}. Assume that $X$ is
endowed with the log structure given by the normal crossing
divisor $E$, and let $X_\log$ denote the resulting log scheme.
We let $W_m\Omega^{q}_{X_{\log}, \log}$ denote the logarithmic Hodge--Witt sheaf
on the log scheme $X_\log$ (cf. \cite[\S~2]{Lorenzon}).
This is the $q$-th Frobenius fixed point of the de Rham--Witt complex
$W_m\Omega^{\bullet}_{X_{\log}}$ of $X_\log$ 
(cf. \cite[\S~4]{Hyodo-Kato}, \cite[\S~1.2]{Geisser-Hesselholt-JAMS}).
As a consequence of \thmref{thm:Main-2}, we get the following extension of
Milne's duality to log schemes.

\begin{cor}\label{cor:Main-4}
  For $q, i \ge 0$ and $m \ge 1$,  there is a
  canonical isomorphism $W_m\Omega^{q}_{X_{\log}, \log}
  \xrightarrow{\cong}  j_* W_m\Omega^q_{U, \log}$. In particular,
  there is a perfect pairing of finite groups
  \[
  H^i_\et(X, W_m\Omega^{q}_{X_{\log}, \log}) \times 
  H^{d+1-i}_\et(X, W_m\Omega^{d-q}_{X|E, \log}) \to \Lambda_m.
   \]
    \end{cor}

\vskip .3cm

\subsection{Duality for semistable schemes}\label{sec:D-SS}
We let $R$ be an equicharacteristic complete discrete valuation ring with finite
residue field.  Let $f \colon X \to \Spec(R)$ be a connected, strictly semistable,
projective scheme of relative dimension $d$ (cf. \cite[Defn.~1.4.2]{Zhao}).
Let $X_s$ denote the closed fiber of $f$.
We let $E \subset X$ be a simple normal crossing divisor with complement $U$.
Let $D \subset X$ be an effective divisor with support $E$.
For an {\'e}tale sheaf $\sF$ on $X$, let $H^*_{\et, X_s}(X,\sF)$ denote the
{\'e}tale cohomology of $\sF$ with support in $X_s$.

In \cite{Zhao}, Zhao proved a duality theorem for certain $p$-primary {\'e}tale
cohomology of $U$. The corresponding result in the mixed
characteristic case was earlier established by Sato \cite{Sato}.
The final duality theorem of this paper is the following generalization of
Zhao's duality to a fixed modulus pair. This can also be considered as an
extension of Milne's duality to the setting of semistable schemes. 

\begin{thm}\label{thm:Main-8}
Given integers $q, i \ge 0, m \ge 1$, there is a canonical isomorphism of
  (not necessarily finite) $\Lambda_m$-modules
  \begin{equation}\label{eqn:Perfect-dvr-0}
    \H^i_\et(X, W_{m}\sF^{q, \bullet}_{X|D}) \xrightarrow{\cong}
\Hom_{\Lambda_m}(H^{d+2-i}_{\et, X_s}(X, W_{m}\Omega^{d+1-q}_{X|D+E, \log}), \Lambda_m).
  \end{equation}
  \end{thm}

\vskip.2cm

\subsection{Layout of the paper}\label{sec:Layout}
We derive the main results of this paper using several fundamental properties of the
filtered de Rham--Witt complex established in \cite{KM-3}.

To prove \thmref{thm:Main-10}, the main additional step in this paper is to equip
the filtered de Rham--Witt complex with the canonical filtration, which was
introduced by Illusie \cite{Illusie} in the classical setting. We show that,
as in the classical case, the canonical filtration on the filtered de Rham--Witt
complex admits a simple description in terms of subquotients of
$\Omega^q_X(\log E)(D)$.
This part of the proof of \thmref{thm:Main-10} is carried out in \S~\ref{sec:FDR}.

By combining Ekedahl's duality in the classical case with a key property of the
filtered de Rham--Witt complex established in \cite{KM-3}
(cf. \thmref{thm:Global-version}(8)), we reduce the proof of
\thmref{thm:Main-10} to establishing a duality theorem for these subquotients.
The proof of this latter result constitutes the technical part of \S~\ref{sec:EK-0}.

To prove \thmref{thm:Main-1}, we show that the filtered de Rham--Witt complex,
originally defined for regular F-finite schemes, admits a natural extension to all
relatively perfect schemes over such schemes, despite the fact that the latter are
generally poorly behaved. Using this extension, we define our motivic complexes
$W_{m}\sF^{q, \bullet}_{X|D}$ and $W_{m}\sG^{q, \bullet}_{X|-D'}$ on the relatively perfect
site of a regular F-finite scheme X. This is carried out in
\S~\ref{sec:RP}.

In \S~\ref{sec:Mot-RP}, we establish two distinguished triangles involving these
complexes in the derived category of sheaves of ${\Z}/{p^m}$-modules on the
relatively perfect site of $X$. These triangles allow us to reduce the proof of
\thmref{thm:Main-1} to the case $m=1$, which we prove directly using the main
results of \cite{Kato-duality-1}. This is done in \S~\ref{sec:KD*}.
Finally, we prove Theorems~\ref{thm:Main-2} and \ref{thm:Main-8} in
\S~\ref{sec:DTHW} as direct applications of \thmref{thm:Main-1}.

\subsection{Notation}\label{sec:Notn}
If $A$ is a commutative ring, an $A$-scheme in this paper
will mean a scheme over $\Spec(A)$. 
We fix a prime $p$ for all of this paper, and assume all schemes to be 
separated $\F_p$-schemes. For a scheme $X$, we let $\Sch_X$ denote the category
of schemes over $X$. We let $\Sm_X$ denote the category of smooth schemes over $X$.
We let $\Sch_{A/\zar}$ (resp. $\Sch_{A/\nis}$, resp. $\Sch_{A/ \et})$
denote the Zariski (resp. Nisnevich, resp. {\'e}tale) site of $\Sch_{A}$.
We let $\epsilon \colon
\Sch_{A/ \et} \to \Sch_{A/ \nis}$ denote the canonical morphism of sites.

For an abelian group $A$, we write $\Tor^1_{\Z}(A, {\Z}/n)$ as
$A[n]$ and $A/{nA}$ as $A/n$.
We denote by $A\{p\}$ (resp. $A\{p'\}$) the
subgroup of elements of $A$ which are annihilated by a power of (resp.
an integer prime to) $p$. 
We denote by $\N_0$ (resp. $\N$) the set of non-negative (resp.
positive) integers, and set $J_m = \{1, \ldots , m\}$.

\section{Filtered de Rham--Witt and its canonical filtration}
\label{sec:FDR}
As this paper is a sequel to \cite{KM-3}, we make extensive use of the results
established therein. We briefly recall the relevant definitions and terminology from
op. cit. and refer the reader to that paper for further details.
In this section, we recall the definition and basic properties of the filtered de
Rham--Witt complex. We then introduce a filtration on the filtered de Rham--Witt
sheaves and establish several key results concerning this filtration.

\subsection{de Rham--Witt complex of $F$-finite schemes}\label{sec:DRFF}
Let $X$ be a scheme. We shall say that $X$ is an $F$-finite scheme if the
absolute Frobenius $F_X \colon X \to X$ is a finite morphism.
We shall say that a commutative ring $A$ is $F$-finite if it is an $\F_p$-algebra
such that $\Spec(A)$ is $F$-finite. If $F$ is a connected and regular $F$-finite
scheme with function field $k(X)$, then there is $n \in \N_0$ such that
$[k(X):k(X)^p] = p^n$ and the Zariski sheaf of absolute K{\"a}hler differentials
$\Omega^1_X$ is locally free of rank $n$. We shall say in this case that $X$ is of
$p$-rank $n$. 

We refer to
\cite[Prop.~2.2]{KM-3} for several useful properties of $F$-finite rings and schemes
(and of their modules of absolute K{\"a}hler differentials) that we shall use.
A particular case of this result that will be implicitly used throughout this paper
is that if $X$ is an essentially of finite type scheme over an $\F_p$-algebra
which is a quasi-excellent henselian discrete valuation ring (hdvr) with an $F$-finite
residue field, then $X$ is $F$-finite. We shall also use the fact that
the de Rham--Witt sheaf $W_m\Omega^q_X$ is a coherent sheaf of $W_m\sO_X$-modules if
$X$ is Noetherian and $F$-finite. In this paper, we shall assume all
schemes to be $F$-finite.

Let $W_mX$ denote the scheme $(|X|, W_m\sO_X)$, where $|X|$ is the underlining
topological space of $X$. Then the transition map
$R\colon W_m \sO_X \surj W_{m-1} \sO_X$ and the Frobenius homomorphism
$F \colon W_m\sO_X \to W_{m-1}\sO_X$ define a closed immersion
$\Re \colon W_{m-1} X \inj W_m X$ and a finite morphism
$\fF \colon W_{m-1} X \to W_m X$, respectively. The standard properties of
$(R, F, V, d)$ then yield the following $W_m \sO_X$-module homomorphisms. 
\begin{equation} \label{eqn:Mor-W_mX*-1}
 W_m \Omega^r_X \xrightarrow{R} \Re_* W_{m-1}\Omega^r_X,\ \ \ \ \ \ \ \ 
 W_m \Omega^r_X \xrightarrow{F} \fF_* W_{m-1}\Omega^r_X
 \end{equation}
 \begin{equation} \label{eqn:Mor-W_mX*-2}
 \fF_* W_{m-1} \Omega^r_X \xrightarrow{V} W_m \Omega^r_X, \ \ \ \ \ \ \ \ 
 \ov{\fF}^m_*  W_m \Omega^r_X \xrightarrow{d}  \ov{\fF}^m_* W_m  \Omega^{r+1}_X,
\end{equation}
 where $\ov{\fF} \colon W_m X \to W_m X$ is the finite morphism induced by the
 absolute Frobenius $\ov{F} \colon W_m \sO_X \inj W_m \sO_X$ which has the property
 that $\ov{F} \circ R = F \colon W_{m+1} \sO_X \to W_m \sO_X$. In particular, if $m=1$,
 then $\fF$ coincides with the absolute Frobenius $F_X$. 
Since $d F= p F d$ and $p \Omega^{r+1}_X = 0$, there is a natural morphism of
$W_{m+1}\sO_X$-modules 
\begin{equation} \label{eqn:Mor-W_mX*-3}
  F^{m-1}d \colon \fF_* W_m\Omega^r_X \to \fF^{m}_* \Omega^{r+1}_X. 
\end{equation}

\subsection{The filtered de Rham--Witt complex}\label{sec:Fil-DRW}
Let  $X$ be a connected, Noetherian, regular scheme and let $E \subset X$ be
a simple normal crossing (snc)
divisor with irreducible components $E_1, \ldots , E_r$.  Let $y_i$ denote the
generic point of $E_i$ for $i \in J_r$. Let
$j \colon U \inj X$ be the inclusion of the complement of $E$ in $X$.
Let $\eta \in X$ denote the generic point of $X$.
For $m \ge 1$, we let
$W_m\Omega^\bullet_X(\log E) \subset j_*W_m\Omega^\bullet_U$ be the de Rham--Witt complex
with log poles along $E$ (cf. \cite[\S~2.1]{KM-3}). This coincides with
the de Rham--Witt complex $W_m\Omega^\bullet_{X_\log}$ of the log scheme $X_\log$,
where the latter is defined by the canonical log
structure on $X$ given by $E$ (cf. \cite[\S~4]{Hyodo-Kato},
\cite[\S~1.2]{Geisser-Hesselholt-JAMS}, \cite[Thm.~2.8]{KM-3}).
For $m \in \N$, we let $\Lambda_m = {\Z}/{p^m}$ and let
$\sD^m_\et(X)$ denote the derived category of {\'e}tale sheaves of
$\Lambda_m$-modules on $X$. For sheaves of $\Lambda_m$-modules $\sF$ and $\sG$
on any site of $X$, we let $\sHom_{\Lambda}(\sF, \sG)$ denote the internal hom sheaf.

We let $\Div(X)$ denote the set of Weil divisors on $X$ and let
$\Div_E(X)$ denote the subset of Weil divisors on $X$ with support in $E$.
We let $\Div^{\eff}_E(X)$ denote the subset
of effective Weil divisors. All these sets are endowed with the canonical partial
order of $\Div(X)$. For $D \in \Div(X)$, we let $|D|$ denote the support of $D$.
For $\n=(n_1, \ldots , n_r)\in \Z^r$ and $t \in \N$, we
define $\n/t:= (\lfloor n_1/t \rfloor, {\ldots}, \lfloor n_r/t \rfloor)$, where
$\lfloor \cdot \rfloor$ (resp. $\lceil  \cdot \rceil$) is the greatest (resp. least)
integer function. We shall say that $\n \ge n'$ (equivalently, $\n' \le \n$)
in $\Z^r$ if $n_i \ge n'_i$ for
every $i$. We shall say that $\n \ge 0$ (resp.  $\n > 0$) if $n_i \ge 0$ for every $i$
(resp. $n_i > 0$ for some $i$). We shall say that $\n < 0$ if $n_i < 0$ for each $i$.
If $D = \sum_i n_i E_i$, we let $D/p = \sum_i \lfloor {n_i}/p \rfloor E_i$.

\begin{defn}\label{defn:Log-fil-3}
  Let $D = {\sum}^r_{i=1} n_iE_i \in \Div_E(X)$ with
  $n_i \in \Z$. For a morphism $u \colon \wt{X} \to X$ of schemes, we let
  $\wt{U} = U \times_X \wt{X}$ and $\wt{D} = D\times_X \wt{X}$. We let
      \[
      \Fil_D W_m(\sO(\wt{U})) = \left\{\underline{a}=(a_{m-1},\ldots, a_0)
      \in W_m(\sO(\wt{U}))| \ 
        a_i^{p^i} \in H^0_\et(\wt{X}, \sO_{\wt{X}}(\wt{D})) \ \forall \ i\right\}.
      \]
      Let $\Fil_D W_m \sO_U$ be the presheaf on $\Sch_X$ such that
      $\Fil_D W_m \sO_U(\wt{X}) = \Fil_D W_m(\sO(\wt{U}))$. One easily checks that
      $\Fil_D\sO_U=\sO_X(D)$ and $\Fil_D W_m \sO_U$ actually restricts
      to a sheaf on $X_\et$. For $q \ge 0$, we let $\Fil_D W_m \Omega^q_U$
      be the $W_m\sO_X$-submodule of the sheaf of $W_m\sO_X$-module
      $j_*W_m\Omega^q_U$ defined by
  \[
  \Fil_D W_m \Omega^q_U =  \Fil_D W_m \sO_U \cdot W_m\Omega^q_X(\log E) +
  d(\Fil_D W_m \sO_U) \cdot W_m\Omega^{q-1}_X(\log E).
  \]
\end{defn}

If $D = {\sum}^r_{i=1} n_i E_i$ where
$\n= (n_1, \ldots , n_r)\in \Z^r$, we use $\Fil_\n W_m \Omega^q_U$
as an alternative notation for $\Fil_{D}W_m \Omega^q_U$.
For a  henselian discrete valuation field (hdvf) $K$ containing $\F_p$ with
uniformizer $\pi$, we shall write
  $\Fil_D W_m\Omega^q_U$ as $\Fil_n W_m\Omega^q_K$ if $X=\Spec(\sO_{K})$,
  $E=V((\pi))$ and $D = nE$ with $n \in \Z$.
The main properties of $\Fil_D W_m\Omega^q_U$ which we shall frequently use
in this paper are summarized in the following two results.

\begin{lem}\label{lem:Log-fil-4}
The sheaf $\Fil_D W_m \Omega^q_U$ satisfies the following properties.
\begin{enumerate}
\item
    $d(\Fil_D W_m \Omega^q_U) \subset \Fil_D W_m \Omega^{q+1}_U$.
 \item
    $F(\Fil_D W_{m+1} \Omega^q_U) \subset \Fil_D W_m \Omega^{q}_U$.
  \item
    $V(\Fil_D W_{m} \Omega^q_U) \subset \Fil_D W_{m+1} \Omega^{q}_U$.
  \item
    $R(\Fil_D W_{m+1} \Omega^q_U) = \Fil_{D/p} W_m \Omega^{q}_U$.
  \item    $\varinjlim_n \Fil_{nE} W_m \Omega^q_U \xrightarrow{\cong}
    j_*W_m\Omega^{q}_U$.
  \item
  $\Fil_D W_{m} \Omega^q_U \subset \Fil_{D'} W_{m} \Omega^q_U$ if
$D \le D'$ and $\Fil_D W_{m} \Omega^q_U \subset  W_{m} \Omega^q_X$ if
    $D < 0$.
  \item
  $\Fil_0 W_m\Omega^q_U = W_m\Omega^q_X(\log E)$.
  \item
    $\Fil_D \Omega^q_U = \Omega^q_X(\log E)(D)$. In particular,
    $\Fil_{-E} \Omega^d_U = \Omega^d_X$, where $d$ is the rank of $\Omega^1_X$.
   \end{enumerate}
\end{lem}
\begin{proof}
See \cite[Lem.~3.6, 3.10]{KM-3}.
\end{proof}

We let $Z_1 \Fil_D W_m\Omega^q_U =  \Fil_D W_m\Omega^q_U \bigcap
j_* Z_1W_m\Omega^q_U$, where 
$Z_1W_m\Omega^q_U := F(W_{m+1}\Omega^q_U) =
\Ker(F^{m-1}d \colon W_m\Omega^q_U \to \Omega^{q+1}_U)$.
Let $C \colon Z_1\Fil_{D}W_m\Omega^q_U \to W_m\Fil_{D/p}\Omega^q_U$ denote the
Cartier homomorphism (cf. \cite[Cor.~6.17]{KM-3}).
We let $B_0\Fil_{D}\Omega^q_U =0$ and $B_1\Fil_{D}\Omega^q_U =
d(\Fil_{D}\Omega^{q-1}_U)$.
We let $Z_0\Fil_{D}\Omega^q_U = \Fil_{D}\Omega^q_U$.
For $i \ge 2$, we let $Z_i\Fil_{D}\Omega^q_U$ (resp. $B_i\Fil_{D}\Omega^q_U$)
be the inverse image of $Z_{i-1}\Fil_{D}\Omega^q_U$
(resp. $B_{i-1}\Fil_{D}\Omega^q_U$) under the composite map
$Z_{i-1}\Fil_D\Omega^q_U \surj
\frac{Z_{i-1}\Fil_{D}\Omega^q_U}{B_1\Fil_{D}\Omega^q_U}
\xrightarrow{\ov{C}} \Fil_{D/{p^{i-1}}}\Omega^q_U$.

\begin{defn}\label{defn:1-c complex} 
  If $q \in \N_0$, and $D \in \Div_E(X)$ is such that $D + E \ge 0$,
  we define
\[
W_m\sF^{q,\bullet}_{X|D} =
\left(Z_1\Fil_{D} W_m\Omega^q_U \xrightarrow{1-C} \Fil_{D} W_m\Omega^q_U\right)
\]
as an object of $\sD^m_\et(X)$.
\end{defn}

It follows from \cite[Cor.~6.17]{KM-3} that this complex is well-defined.
We shall use the following key result repeatedly in this paper.

\begin{thm}\label{thm:Global-version}
Given any $D \in \Div_E(X)$, we have the following.
  \begin{enumerate}
   \item
    There is a short exact sequence of Zariski sheaves of $W_m\sO_X$-modules
    \[
    0 \to V^{m-1}(\Fil_D\Omega^q_U) + dV^{m-1}(\Fil_D\Omega^{q-1}_U) \to
    \Fil_DW_m\Omega^q_U \xrightarrow{R}  \Fil_{D/p}W_{m-1}\Omega^q_U \to 0.
    \]
    \item
      $\Ker(F^{m-1}d \colon \Fil_{D}W_m\Omega^q_U \to \Fil_{D}\Omega^{q+1}_U)=
      Z_1\Fil_{D}W_m\Omega^q_U = F(\Fil_{D}W_{m+1}\Omega^q_U) $.
  \item
    There exist isomorphisms of Zariski sheaves of $W_{m+1}\sO_X$-modules
   \[
    \ov{F} \colon \Re_*\Fil_{{D}/p}W_{m}\Omega^q_U \xrightarrow{\cong}
    \frac{\fF_*Z_1\Fil_{D}W_m\Omega^q_U}{\fF_*dV^{m-1}(\Fil_{D}\Omega^{q-1}_U)}; \
    \ov{C} \colon  \frac{\fF_*Z_1\Fil_{D}W_m\Omega^q_U}{\fF_*dV^{m-1}(\Fil_{D}\Omega^{q-1}_U)}
     \xrightarrow{\cong} \Re_*\Fil_{{D}/p}W_{m}\Omega^q_U.
    \]
  \item
    $\Ker(F^{m-1}\colon \Fil_{D}W_m\Omega^q_U \to \Fil_{D}\Omega^q_U) =
    V(\Fil_{D}W_{m-1}\Omega^q_U)$.
  \item
    $\Ker(dV^{m-1}\colon \Fil_{D}\Omega^q_U \to \Fil_{D}W_m\Omega^{q+1}_U)
    = F^m(\Fil_{D}W_{m+1}\Omega^q_U)= Z_m\Fil_D\Omega^q_U$.
  \item
    $\Ker(V \colon \Fil_{D}W_{m}\Omega^q_U \to \Fil_{D}W_{m+1}\Omega^{q}_U)
    = dV^{m-1}(\Fil_{D}\Omega^{q-1}_U)$.
  \item
$\Ker(V^m \colon \Fil_{D}\Omega^q_U \to \Fil_{D}W_{m+1}\Omega^{q}_U)
    = F^mdV(\Fil_{D}W_m\Omega^{q-1}_U)= B_m\Fil_D\Omega^q_U$.
  \item
  If $D + E \ge 0$, there exists a distinguished triangle
    \[
     W_1\sF^{q,\bullet}_{X|D} \xrightarrow{V^{m-1}} W_m\sF^{q,\bullet}_{X|D}
     \xrightarrow{R} W_{m-1}\sF^{q,\bullet}_{X|{D/p}} \xrightarrow{+}
     W_1\sF^{q,\bullet}_{X|D} [1]
    \]
    in the derived category of sheaves of $\Lambda_m$-modules on
    $X_\zar$.
  \item $\varinjlim\limits_{D\ge -E} W_m\sF^{q,\bullet}_{X|D}
    \xrightarrow{\simeq} \left(j_* Z_1W_m\Omega^q_U
    \xrightarrow{1-C} j_*W_m\Omega^q_U\right) = {\bf R}j_* (W_m\Omega^q_{U,\log}).$
     \item 
       $\Ker(dV^{m-1}\colon \Fil_{D}\Omega^q_U \to \Fil_{D}W_m\Omega^{q+1}_U/V^{m-1}\Fil_D
       \Omega^{q+1}_U) = F^{m-1}(\Fil_{D}W_{m}\Omega^q_U)$.
   \item
     $\Ker(V^{m-1} \colon \Fil_{D}\Omega^q_U \to
     \Fil_{D}W_{m}\Omega^{q}_U/dV^{m-1}\Fil_{D}\Omega^{q-1}_U )
     = F^{m}dV(\Fil_{D}W_{m}\Omega^{q-1}_U)$.
   \item
     $j^*(W_m\sF^{q,\bullet}_{X|D})$ is quasi-isomorphic to
     $W_m\Omega^q_{U, \log}$.
   \item
     The inclusion $\Fil_{-E}W_m\Omega^d_U \inj W_m\Omega^d_X$ is a bijection if
     $X$ has $p$-rank $d$.
     \end{enumerate}
\end{thm}
\begin{proof}
  See \cite[Thm.~6.19]{KM-3}.
  \end{proof}

\subsection{Canonical filtration of $\Fil_D W_m\Omega^q_U$}\label{sec:HWD-0}
Let $k$ be an $F$-finite field of characteristic $p$ and let 
$X$ be a connected, Noetherian, regular $F$-finite $k$-scheme. 
Let $E$ be an snc divisor on  $X$ with irreducible components $E_1, \ldots , E_r$.
Let $j \colon U \inj X$ be the
inclusion of the complement of $E$ in $X$. Let $K$ be the function field of $X$
and $d$ be the $p$-rank of $X$.
We fix integers $q \ge 0, \ m \ge 1$. We let $D \in \Div_E(X)$ and
define the `canonical filtration' on $\Fil_D W_m\Omega^q_U$ as follows.

\begin{defn}\label{defn:Can-filt}
For $n \in \Z$, we let
  \[
  \filt^n\Fil_D W_m\Omega^q_U =  \left\{\begin{array}{ll}
  \Ker(R^{m-n} \colon \Fil_D W_m\Omega^q_U \surj
  \Fil_{{D}/{p^{m-n}}} W_n\Omega^q_U) & \mbox{if $n < m$} \\
 0 & \mbox{if $n  \ge  m$.}
\end{array}\right.
\]
\end{defn}

The surjectivity of $R^{m-n}$ in the above definition follows by
\cite[Thm.~6.19]{KM-3} and an
elementary observation that $\lfloor{{\lfloor{n/{p^s}}\rfloor}/p}\rfloor =
    \lfloor{n/{p^{s+1}}}\rfloor$ for $n \in \Z$ and $s \ge 0$.
    We have a decreasing filtration
   \[
    \Fil_D W_m\Omega^q_U =  \filt^0\Fil_D W_m\Omega^q_U \supset
      \filt^{1}\Fil_D W_m\Omega^q_U \supset \cdots \supset    
      \filt^{m-1}\Fil_D W_m\Omega^q_U \supset \filt^m\Fil_D W_m\Omega^q_U = 0.
     \]
We let
\begin{equation}\label{eqn:Can-filt-0}
  \gr^n \Fil_D W_m\Omega^q_U =
  \frac{\filt^n\Fil_D W_m\Omega^q_U}{\filt^{n+1}\Fil_D W_m\Omega^q_U}; \ 0 \le n \le m-1.
  \end{equation}
It is easy to check that for $0 \le n \le m-1$,
there is a canonical isomorphism $W_{m}\sO_X$-modules
\begin{equation}\label{eqn:Can-filt-1}
  \begin{array}{lll}
    \gr^n \Fil_D W_m\Omega^q_U & \xrightarrow{\cong} &
    \Re^{m-n-1}_* \filt^{n} \Fil_{D/{p^{m-n-1}}} W_{n+1}\Omega^q_U \\
& = & \Ker(\Re^{m-n-1}_* \Fil_{D/{p^{m-n-1}}} W_{n+1}\Omega^q_U
    \xrightarrow{R} \Re^{m-n}_* \Fil_{D/{p^{m-n}}} W_{n}\Omega^q_U).
    \end{array}
\end{equation}

The following result is the extension of \cite[Cor.~I.3.9]{Illusie} to
the filtered de Rham--Witt complex.

\begin{lem}\label{lem:Can-filt-3}
  For $n \in \N$, we have the following.
  \begin{enumerate}
  \item
   There is a canonical exact sequence of $W_{n+1}\sO_X$-modules
  \[
  0 \to \fF^n_* \frac{\Fil_{D} \Omega^q_U}{B_n\Fil_{D} \Omega^q_U}
  \xrightarrow{V^n} \filt^{n} \Fil_{D} W_{n+1}\Omega^q_U \xrightarrow{\beta_n}
   \fF^n_* \frac{\Fil_{D} \Omega^{q-1}_U}{Z_n\Fil_{D} \Omega^{q-1}_U} \to 0.
   \]
   \item
   There is a canonical exact sequence of $W_{n+2}\sO_X$-modules
  \[
  0 \to \fF^{n+1}_* \frac{\Fil_{D} \Omega^{q-1}_U}{Z_{n+1}\Fil_{D} \Omega^{q-1}_U}
  \xrightarrow{dV^n} \fF_* \filt^{n} \Fil_{D} W_{n+1}\Omega^q_U  \xrightarrow{\alpha_n}
   \fF^{n+1}_* \frac{\Fil_{D} \Omega^{q}_U}{B_{n+1}\Fil_{D} \Omega^{q}_U} \to 0.
   \]
  \end{enumerate}
  \end{lem}
\begin{proof}
  We begin by proving (1).
  By items (1) and (7) of \thmref{thm:Global-version},
  we have an exact sequence
  \begin{equation}\label{eqn:Can-filt-3-0}
    0 \to \fF^n_* \frac{\Fil_{D} \Omega^q_U}{B_n\Fil_{D} \Omega^q_U} \xrightarrow{V^n}
    \filt^{n} \Fil_D W_{n+1}\Omega^q_U \to
    \frac{\Fil_{D} W_{n+1}\Omega^q_U}{V^n(\Fil_{D}\Omega^q_U)}
  \end{equation}
and $\coker(V^n)$ is generated by the image of $dV^n(\Fil_{D}\Omega^{q-1}_U)$.
Combining this with the items (6) and (11) of \thmref{thm:Global-version}, we get
the exact sequence (1) in which the map $\beta_n$ is defined by the
property that $\beta_n(V^n(a) + dV^n(b))$ is the image of $b$ under the quotient map
$\Fil_{D} \Omega^{q-1}_U \surj
\frac{\Fil_{D} \Omega^{q-1}_U}{Z_n\Fil_{D} \Omega^{q-1}_U}$.

To prove (2), we use item (1) and (5) of
\thmref{thm:Global-version} to get an exact sequence
\[
0 \to \fF^{n+1}_* \frac{\Fil_{D} \Omega^{q-1}_U}{Z_{n+1}\Fil_{D} \Omega^{q-1}_U}
\xrightarrow{dV^n} \fF_* \filt^{n} \Fil_{D} W_{n+1}\Omega^q_U  \to
\frac{\Fil_{D} W_{n+1}\Omega^q_U}{dV^n(\Fil_{D} \Omega^{q-1}_U)},
\]
where $dV^n$ is $W_{n+2}\sO_X$-linear (cf. \cite[(5.6)]{KM-3})
  and the image of the last map is generated by $V^n(\Fil_{D}\Omega^q_U)$.
  Combining this with items (7) and (11) of \thmref{thm:Global-version}, we get
  the exact sequence (2) in which $\alpha_n$ is defined by the property that
  $\alpha_n(V^n(a) + dV^n(b))$ is the image of $a$ under the quotient map
  $\Fil_{D} \Omega^{q}_U \surj
  \frac{\Fil_{D} \Omega^{q}_U}{B_{n+1}\Fil_{D} \Omega^{q}_U}$.
\end{proof}

Recall from \cite[Chap.~I]{Illusie} that the multiplication map
$p \colon W_m\Omega^q_K \to W_{m}\Omega^q_K$ has a factorization
$W_m\Omega^q_K \stackrel{R}{\surj}  W_{m-1}\Omega^q_K \stackrel{\ov{p}}{\inj}
W_{m}\Omega^q_K$.
By the definition of the $W_{n+1}\sO_X$-linear map $\ov{p}$,
and items (2) and (7) of \thmref{thm:Global-version}, there is a unique
$W_{n+1}\sO_X$-linear
map $\theta_n \colon \fF^n_* B_n\Fil_D \Omega^{q+1}_U \to 
\frac{\Fil_D W_{n+1}\Omega^q_U}{\ov{p}(\Fil_{D/p} W_{n}\Omega^q_U)}$ such that
$V = \theta_n \circ F^{n-1}d \colon \Fil_D W_{n}\Omega^q_U \to  
\frac{\Fil_D W_{n+1}\Omega^q_U}{\ov{p}(\Fil_{D/p} W_{n}\Omega^q_U)}$.

\begin{lem}\label{lem:Can-filt-4}
For $n \in \N_0$, there is a canonical exact sequence of $W_{n+1}\sO_X$-modules
  \[
  0 \to \fF^n_* B_n\Fil_D \Omega^{q+1}_U \xrightarrow{\theta_n} 
  \frac{\Fil_D W_{n+1}\Omega^q_U}{\ov{p}(\Fil_{D/p} W_{n}\Omega^q_U)}
  \xrightarrow{F^n} \fF^n_* Z_n \Fil_D \Omega^q_U \to 0.
\]
\end{lem}
\begin{proof}
We begin by noting that $F^n$ is surjective by item (5) of
  \thmref{thm:Global-version}.
  To show the exactness at the middle term, we let $x \in \Fil_D W_{n+1}\Omega^q_U$
  be any element whose
  image $\ov{x}$ in $\frac{\Fil_D W_{n+1}\Omega^q_U}{\ov{p}(\Fil_{D/p} W_{n}\Omega^q_U)}$
  lies in the kernel of $F^n$. By item (4) of \thmref{thm:Global-version},
  we can write $x = V(y)$, where $y \in \Fil_D W_n\Omega^q_U$.
  This implies that $\ov{x} = \theta_n(F^{n}d(y))$, and hence we are done
  by item (7) of \thmref{thm:Global-version}. It remains to prove that $\theta_n$
  is injective.

  To that end, we look at the commutative diagram
  \begin{equation}\label{eqn:Can-filt-4-0}
    \xymatrix@C1.5pc{
     \fF^n_* B_n\Fil_D \Omega^{q+1}_U \ar[r]^-{\theta_n} \ar[d] & 
  \frac{\Fil_D W_{n+1}\Omega^q_U}{\ov{p}(\Fil_{D/p} W_{n}\Omega^q_U)} \ar[d] \\
j_* \fF^n_* B_n \Omega^{q+1}_U \ar[r]^-{\theta_n}  & 
  j_* \frac{W_{n+1}\Omega^q_U}{\ov{p}(W_{n}\Omega^q_U)},}
  \end{equation}
  where the vertical arrows are the canonical maps induced by the inclusion
  $\Fil_D W_{\star}\Omega^\bullet_U \inj j_* W_{\star}\Omega^\bullet_U$.
  The bottom horizontal row is injective by \cite[Lem.~0.6]{Ekedahl} and
  N{\'e}ron--Popescu approximation. The left
  vertical arrow is injective by \cite[Lem.~6.5]{KM-3}. It follows that
  $\theta_n$ on the top row is injective. 
\end{proof}

\begin{lem}\label{lem:Can-filt-5}
  For any $n \in \N_0$, the sheaves of $\sO_X$-modules
  \[
  \ov{\fF}^n_* B_n \Fil_D \Omega^q_U, \ \ov{\fF}^n_* Z_n \Fil_D \Omega^q_U, \
  \ov{\fF}^n_* \frac{\Fil_D \Omega^q_U}{B_n \Fil_D \Omega^q_U} \ \mbox{and} \
  \ov{\fF}^n_* \frac{\Fil_D \Omega^q_U}{Z_n \Fil_D \Omega^q_U}
  \]
  are coherent and locally free.
\end{lem}
\begin{proof}
  Since $X$ is $F$-finite, we have seen before that the sheaves in question are
  all coherent $\sO_X$-modules. To prove their locally freeness, we
  can assume $n \ge 1$ by \cite[Lem.~3]{KM-3} and \lemref{lem:Log-fil-4}(8).
  We can also assume that $X = \Spec(A)$ is local.
  We now begin by noting that
  there is an exact sequence of $\sO_X$-modules
\begin{equation}\label{eqn:Can-filt-5-0}
  0 \to  \ov{\fF}^n_* \frac{Z_n \Fil_D \Omega^q_U}{B_n \Fil_D \Omega^q_U} \to
  \ov{\fF}^n_* \frac{\Fil_D \Omega^q_U}{B_n \Fil_D \Omega^q_U} \to 
 \ov{\fF}^n_* \frac{\Fil_D \Omega^q_U}{Z_n \Fil_D \Omega^q_U} \to 0.
\end{equation}

On the other hand, we have the Cartier isomorphism
$\ov{C}^n \colon \ov{\fF}^n_* \frac{Z_n \Fil_D \Omega^q_U}{B_n \Fil_D \Omega^q_U}
\xrightarrow{\cong} \Fil_{D/{p^n}} \Omega^q_U$ by \cite[Lem.~6.3, 6.4]{KM-3}, and
locally freeness of $\Fil_{D/{p^n}} \Omega^q_U$ by
\lemref{lem:Log-fil-4}(3).
It follows that $\ov{\fF}^n_* \frac{\Fil_D \Omega^q_U}{B_n \Fil_D \Omega^q_U}$ is
locally free if $\ov{\fF}^n_* \frac{\Fil_D \Omega^q_U}{Z_n \Fil_D \Omega^q_U}$ is so.

We shall now prove the locally freeness of
$\ov{\fF}^n_* \frac{\Fil_D \Omega^q_U}{Z_n \Fil_D \Omega^q_U}$ by induction on
$n \ge 1$.
For $n =1$, we use the isomorphism
$d \colon \ov{\fF}_* \frac{\Fil_D \Omega^q_U}{Z_1 \Fil_D \Omega^q_U}
\xrightarrow{\cong} \ov{\fF}_* B_1 \Fil_D \Omega^q_U$ and the exact sequence
\begin{equation}\label{eqn:Can-filt-5-1}
  0 \to \ov{\fF}_* B_1 \Fil_D \Omega^q_U \to \ov{\fF}_* Z_1 \Fil_D \Omega^q_U
  \xrightarrow{C} \Fil_{D/p} \Omega^q_U \to 0.
\end{equation}
Since $\Fil_{D/p} \Omega^q_U$ is locally free, we are left with showing that
$\ov{\fF}_* Z_1 \Fil_D \Omega^q_U$ is locally free.

To that end, we recall that if
$(x_1, \ldots , x_r, x_{r+1}, \ldots , x_N)$ is the maximal ideal of $A$ and $(\pi)$ is
  a defining ideal of $E$, then there are elements  $y_1, {\ldots},y_s \in A$ such that
  $\Omega^1_A(\log \pi)$ is a free $A$-module with basis
  $\{\dlog(x_1), \ldots, \dlog(x_r), dx_{r+1}, \ldots , dx_N,$
  $dy_1, {\ldots}, dy_s\}$ (cf. \cite[Lem.~2.3]{KM-3}). Also,
  $\{x_1,\ldots,x_N,y_1,\ldots,y_s\}$ forms a differential basis
  (cf. \cite[Prop.~2.2(7)]{KM-3}) and hence a $p$-basis
  (see \cite[\S~II, \S~III, Thm.~1]{Tyc} for the definition of $p$-basis and its
  relation with differential basis) of $A$ over $\F_p$. 

  Using the above facts about the local rings of $X$, it is easy to check that
  locally, (depending on $D$) we can find two finite dimensional $\F_p$-vector
  spaces $V_1, V_2$ and an $\F_p$-linear map $d_0 \colon V_1 \to V_2$ such that
  $d \colon \Omega^{q-1}_X(\log E)(D) \to \Omega^{q}_X(\log E)(D)$ is
  locally of the form $(V_1 \xrightarrow{d_0} V_2)\otimes_{\F_p} \sO_{X}^p$
    (cf. \cite[Lem.~1.7]{Milne-ENS}).
    In particular, $\Ker(d)$ and $\coker(d)$ are locally free $\sO^p_{X}$-modules
    of finite rank. The ring isomorphism $\ov{F} \colon \sO_X \xrightarrow{\cong}
    \sO^p_X$ now implies that they are locally free $\sO_X$-modules of finite rank.
    We can now apply \lemref{lem:Log-fil-4}(8) again to conclude that
    $\ov{\fF}_* Z_1\Fil_{D}\Omega^q_U$ and
    $\ov{\fF}_* \frac{\Fil_{D}\Omega^{q}_U}{d(\Fil_{D}\Omega^{q-1}_U)} =
    \ov{\fF}_* \frac{\Fil_{D}\Omega^{q}_U}{B_1 \Fil_{D}\Omega^{q}_U}$ are
    locally free sheaves of $\sO_X$-modules.

    For $n \ge 2$, we combine the exact sequence
\begin{equation}\label{eqn:Can-filt-5-2}
0 \to \ov{\fF}^n_* \frac{Z_1\Fil_{D}\Omega^q_U}{Z_n\Fil_{D}\Omega^q_U}
\to \ov{\fF}^n_* \frac{\Fil_D \Omega^q_U}{Z_n \Fil_D \Omega^q_U} \to
\ov{\fF}^n_* \frac{\Fil_D \Omega^q_U}{Z_1 \Fil_D \Omega^q_U} \to 0
    \end{equation}
 with the isomorphism $\ov{C} \colon 
 \ov{\fF}^n_* \frac{Z_1\Fil_{D}\Omega^q_U}{Z_n\Fil_{D}\Omega^q_U} \xrightarrow{\cong}
 \ov{\fF}^n_* \frac{\Fil_{D/p} \Omega^q_U}{Z_{n-1} \Fil_{D/p} \Omega^q_U}$
 and use induction on $n$ to complete the proof.

 Finally, we use the exact sequences
  \begin{equation}\label{eqn:Can-filt-5-3}
0 \to \ov{\fF}^n_* B_n \Fil_D \Omega^q_U \to \fF^n_* \Fil_D \Omega^q_U
\to \ov{\fF}^n_* \frac{\Fil_D \Omega^q_U}{B_n \Fil_D \Omega^q_U} \to 0; \ \ \mbox{and}
  \end{equation}
  \begin{equation}\label{eqn:Can-filt-5-4}
0 \to \ov{\fF}^n_* Z_n \Fil_D \Omega^q_U \to \ov{\fF}^n_* \Fil_D \Omega^q_U
\to \ov{\fF}^n_* \frac{\Fil_D \Omega^q_U}{Z_n \Fil_D \Omega^q_U} \to 0
  \end{equation}
  to conclude the locally freeness of $\ov{\fF}^n_* B_n \Fil_D \Omega^q_U$ and
  $\ov{\fF}^n_* Z_n \Fil_D \Omega^q_U$. This completes the proof.
\end{proof}
   
\section{Ekedahl's duality for modulus pairs}\label{sec:EK-0}
In this section,  we prove \thmref{thm:Main-10}. We continue working in the setting of
\S~\ref{sec:HWD-0}. We assume additionally that $k$ is a perfect field so that
$d = \dim(X)$. We let $S = \Spec(k)$ and
 let $f_m \colon W_mX \to W_mS$ denote the projection
 map induced by the structure map $f \colon X \to S$.  We have a commutative diagram
\begin{equation}\label{eqn:Can-filt-2}
  \xymatrix@C1.5pc{
  W_{m-n}X \ar[r]^-{\Re^{n}_X} \ar[d]_-{f_{m-n}} & W_mX \ar[d]^-{f_m} \\
  W_{m-n}S \ar[r]^-{\Re^{n}_S} & W_mS.}
\end{equation}
The horizontal arrows are closed immersions. The vertical arrows are
are proper (resp. projective) if $f$ is proper (resp. projective). 
We let $\iota = \Re^{m-1}_S$ and $\nu = \Re^{m-1}_X$.

By Ekedahl's duality (cf. \cite[Thm.~I.3.4]{Ekedahl}), there is a canonical
isomorphism $\Tr_m \colon W_m\Omega^d_X$ \\
$\cong f^{!}_m(W_m(k))[-d]$ in the bounded
derived category of the sheaves of coherent $W_m\sO_X$-modules. Furthermore,
one has a canonical isomorphism $\delta_{m,n} \colon
(\Re^{m-n}_S)^{!}(W_m(k)) \cong W_n(k)$ for
$1 \le n \le m-1$. In particular, $(\Re^{m-n}_X)^{!}(W_m\Omega^d_X) \cong
W_n\Omega^d_X$.

\begin{lem}\label{lem:Ext-0}
For any locally free $\sO_X$-module $\sM$, one has
  ${\sE}xt^i_{W_m\sO_X}(\nu_*\sM, W_m\Omega^d_X) = 0$ for all $i > 0$.
\end{lem}
\begin{proof}
  This is an elementary exercise using the isomorphisms $\Tr_m$ and $\delta_{m,1}$,
 the identities $\nu^{!} \circ f^{!}_m = (f_m \circ \nu)^{!} =
 (\iota \circ f_1)^{!} = f^{!}_1 \circ \iota^{!}$ and the adjointness of the
 pair $(\nu_*, \nu^{!})$, see the proof of \cite[Lem.~II.2.2.7]{Ekedahl}.
  \end{proof}

\begin{lem}\label{lem:Ext-1}
One has ${\sE}xt^i_{W_m\sO_X}(\Fil_D W_m\Omega^q_U, W_m\Omega^d_X) = 0$ for $i > 0$.
\end{lem}
\begin{proof}
  It is enough to show that
  ${\sE}xt^i_{W_m\sO_X}(\gr^n \Fil_D W_m\Omega^q_U, W_m\Omega^d_X) = 0$ for $i > 0$.
  But this follows by using ~\eqref{eqn:Can-filt-1},
  applying Lemma~\ref{lem:Can-filt-3}(1) to $D/{p^{m-n-1}}$ and taking
  $\Re^{m-n-1}_{X *}$
    of the resulting exact sequence for $0 \le n \le m-1$, and then
    subsequently using Lemmas~\ref{lem:Can-filt-5} and ~\ref{lem:Ext-0}
    as well as the identity $\fF^n_* = (\Re^n_X)_* \circ \ov{\fF}^n_*$.
\end{proof}

For the remainder of this section, we fix $q \in \N_0, D \in \Div_E(X)$. We let
$q' = d-q$ and $D' = -D-E$. We consider the diagram
\begin{equation}\label{eqn:Ext-4}
  \xymatrix@C1.5pc{
  \Fil_D W_m\Omega^q_U \times \Fil_{D'} W_m\Omega^{q'}_U  \ar[r]^-{\< \ , \ \>}
    \ar@<5ex>[d]^-{R_m} &
  W_m\Omega^d_X  \\
 \Fil_{D/p} W_{m-1}\Omega^q_U \times \Fil_{{D'}/p} W_{m-1}\Omega^{q'}_U
  \ar[r]^-{\< \ , \ \>}
  \ar@<7ex>[u]^-{\ov{p}} & W_{m-1}\Omega^d_X \ar[u]_-{\ov{p}},}
  \end{equation} 
where the horizontal arrows are the wedge product pairings
(cf. \cite[Lem.~3.10, Cor.~5.14]{KM-3}).

\begin{lem}\label{lem:Ext-5}
  The diagram ~\eqref{eqn:Ext-4} commutes. In particular, it induces a
  wedge product pairing
  \begin{equation}\label{eqn:Ext-5-0}
    \coker(\ov{p}) \times \Ker(R_m) \xrightarrow{\< \ , \ \>_m} W_m\Omega^d_X.
    \end{equation}
    \end{lem}
\begin{proof}
  We let $\alpha \in \Fil_{D/p} W_{m-1}\Omega^q_U$ and
  $\beta \in \Fil_{D'} W_m\Omega^{q'}_U$. Using \thmref{thm:Global-version}(1),
  choose $\alpha_1 \in
  \Fil_{D} W_{m}\Omega^q_U$ such that $\alpha = R(\alpha_1)$.
  We then get $\ov{p}(\alpha) \wedge \beta = (p\alpha_1) \wedge \beta =
  p (\alpha_1 \wedge \beta) = \ov{p}(\alpha \wedge R(\beta))$, and this
  proves the lemma.
\end{proof}

\begin{lem}\label{lem:Ext-6}
For $n \in \N_0$, ~\eqref {eqn:Ext-5-0} induces a commutative diagram of pairings 
\begin{equation}\label{eqn:Ext-6-0}
  \xymatrix@C2.5pc{
    0 \hspace*{1.5cm} 0 \ar@<-5ex>[d] & \\
    \fF^n_*B_n\Fil_D\Omega^{q+1}_U \times
    \fF^n_* \frac{\Fil_{D'}\Omega^{q'-1}_U}{Z_n\Fil_{D'}\Omega^{q'-1}_U}
    \ar[r]^-{\< \ , \ \>_{n+1}} \ar@<-5ex>[d]_-{\theta_n} \ar@<-5ex>[u] &
    W_{n+1}\Omega^d_X \ar@{=}[d] \\
    \coker(\ov{p}) \times \Ker(R_{n+1}) \ar[r]^-{\< \ , \ \>_{n+1}}
    \ar@<-5ex>[u]_-{\beta_n}
    \ar@<-5ex>[d]_-{F^n} & W_{n+1}\Omega^d_X \\
   \fF^n_* Z_n\Fil_D \Omega^q_U \times  
   \fF^n_* \frac{\Fil_{D'} \Omega^{q'}_U}{B_n\Fil_{D'} \Omega^{q'}_U}
   \ar[r]^-{\< \ , \ \>_{n+1}} \ar@<-5ex>[u]_-{V^n} \ar@<-5ex>[d] &
    W_{n+1}\Omega^d_X \ar@{=}[u] \\
   0  \hspace*{1.5cm} 0 \ar@<-5ex>[u]}
\end{equation}
\end{lem}
\begin{proof}
  Suffices to prove the lemma when $D = \emptyset$ for which we refer to
  \cite[\S~II.1, p.~201]{Ekedahl}.
\end{proof}

We next claim that the cup product pairing
$\< \ , \ \>_1 \colon \Fil_D\Omega^{q+1}_U \times
\Fil_{D'}\Omega^{q'-1}_U \to \Omega^d_X$
induces a pairing $\< \ , \ \>_1 \colon \ov{\fF}^n_*B_n\Fil_D\Omega^{q+1}_U \times
\ov{\fF}^n_* \frac{\Fil_{D'}\Omega^{q'-1}_U}{Z_n\Fil_{D'}\Omega^{q'-1}_U} \to \Omega^d_X$,
given by $(a,b) \mapsto C^n(\<a, b\>_1)$.
To prove this claim, we can reduce to the case when $D = \emptyset$ which is
done in \cite[\S~II.1, p.~201]{Ekedahl}.
One similarly checks that the cup product pairing
$\< \ , \ \>_1 \colon \Fil_D\Omega^{q}_U \times \Fil_{D'}\Omega^{q'}_U \to \Omega^d_X$
induces a pairing $\< \ , \ \>_1 \colon \ov{\fF}^n_*Z_n\Fil_D\Omega^{q}_U \times
\ov{\fF}^n_* \frac{\Fil_{D'}\Omega^{q'}_U}{B_n\Fil_{D'}\Omega^{q'}_U} \to \Omega^d_X$.
We let $\sF_n = \frac{\Fil_{D'}\Omega^{q'-1}_U}{Z_n\Fil_{D'}\Omega^{q'-1}_U}$
and $\sG_n = \frac{\Fil_{D'}\Omega^{q'}_U}{B_n\Fil_{D'}\Omega^{q'}_U}$.

\begin{lem}\label{lem:Ext-7}
  The diagram
\begin{equation}\label{eqn:Ext-7-0}
  \xymatrix@C1.6pc{
    \fF^n_*B_n\Fil_D\Omega^{q+1}_U \ar[rr]^-{\eta_1} \ar[ddd]_-{(-1)^{q+1}} & &
   \sHom_{W_{n+1}\sO_X}({\fF}^n_* \sF_n, W_{n+1}\Omega^d_X)
    \\
    & & \nu_* \sHom_{\sO_X}(\ov{\fF}^n_* \sF_n, \nu^{!} W_{n+1}\Omega^d_X)
    \ar[u]_-{\rm adjunction} \ar[u]^-{\cong} \\
    & & \nu_* \sHom_{\sO_X}(\ov{\fF}^n_* \sF_n, \Omega^d_X) \ar[u]^-{\cong} \\
    {\fF}^n_*B_n\Fil_D\Omega^{q+1}_U \ar[r]^-{\cong} &
    \nu_* \ov{\phi}^n_*B_n\Fil_D\Omega^{q+1}_U \ar[r]^-{\eta'_1} &
    \nu_* \sHom_{\sO_X}(\ov{\fF}^n_* \sF_n, \ov{\fF}^n_* \frac{Z_n\Omega^d_X}{B_n\Omega^d_X}) \ar[u]_-{C^n \circ (-)} \ar[u]^-{\cong} }
\end{equation}
commutes if we let $\eta_1$ (resp. $\eta'_1$) denote the map induced by the
uppermost pairing $\< \ , \ \>_{n+1}$ (resp. $\< \ , \ \>_1$) in \lemref{lem:Ext-6}.
\end{lem}
\begin{proof}
 By \cite[Lem.~II.2.2.4]{Ekedahl}, the composite vertical arrow on the right is
  the map induced by the composition with $V^n$. It suffices therefore to check that
  $\eta_1 = (-1)^{q+1} (V^n \circ \eta'_1)$. 
  To prove the latter identity, we let $\alpha \in B_n\Fil_D\Omega^{q+1}_U$ and
  ${\beta} \in \Fil_{D'}\Omega^{q'-1}_U$. Let $\ov{\beta}$ be the image of $\beta$
    in $\sF_n$. We write $\alpha = F^{n-1}d(\alpha')$, where $\alpha' \in
    \Fil_D W_n\Omega^q_U$. We then have $\theta_n(\alpha) = V(\alpha')$ by 
    definition of $\theta_n$ (cf. \S~\ref{sec:HWD-0}).
    On the other side, we have that $dV^n(\beta) \in \Ker(R)$ such that
    $\beta_n(dV^n(\beta)) = \ov{\beta}$.
    We thus get $\eta_1(\alpha)(\ov{\beta}) = V(\alpha') \wedge dV^n(\beta) \in
    W_{n+1}\Omega^d_X$ by definition of the pairing $\< \ , \ \>_m$.

    We now compute
    \begin{equation}\label{eqn:Ext-7-1}
      \begin{array}{lll}
        V^n(\eta'_1(\alpha)(\ov{\beta})) & = & V^n(\alpha \wedge \beta) \\
        & = & V^n(F^{n-1}d(\alpha') \wedge \beta)
         =  V(d\alpha' \wedge V^{n-1}(\beta)) \\
        & = & V\left(d(\alpha' \wedge V^{n-1}(\beta)) - (-1)^q \alpha' \wedge
        dV^{n-1}(\beta)\right) \\
        & = & pdV(\alpha' \wedge V^{n-1}(\beta)) + (-1)^{q+1} V(\alpha') \wedge
        dV^n(\beta) \\
        & = & (-1)^{q+1} V(\alpha') \wedge
        dV^n(\beta) 
         =  \eta_1(\alpha)(\ov{\beta}).
      \end{array}
    \end{equation}
    This completes the proof.
   \end{proof}

\begin{lem}\label{lem:Ext-8}
  The diagram
\begin{equation}\label{eqn:Ext-8-0}
  \xymatrix@C1.6pc{
    \fF^n_*Z_n\Fil_D\Omega^{q}_U \ar[rr]^-{\eta_2} \ar@{=}[ddd] & &
   \sHom_{W_{n+1}\sO_X}({\fF}^n_* \sG_n, W_{n+1}\Omega^d_X)
    \\
    & & \nu_* \sHom_{\sO_X}(\ov{\fF}^n_* \sG_n, \nu^{!} W_{n+1}\Omega^d_X)
    \ar[u]_-{\rm adjunction} \ar[u]^-{\cong} \\
    & & \nu_* \sHom_{\sO_X}(\ov{\fF}^n_* \sG_n, \Omega^d_X) \ar[u]^-{\cong} \\
    {\fF}^n_*Z_n\Fil_D\Omega^{q}_U \ar[r]^-{\cong} &
    \nu_* \ov{\fF}^n_*Z_n\Fil_D\Omega^{q}_U \ar[r]^-{\eta'_2} &
    \nu_* \sHom_{\sO_X}(\ov{\fF}^n_* \sG_n, \ov{\fF}^n_* \frac{Z_n\Omega^d_X}{B_n\Omega^d_X}) \ar[u]_-{C^n \circ (-)} \ar[u]^-{\cong} }
\end{equation}
commutes if we let $\eta_2$ (resp. $\eta'_2$) denote the map induced by the
lowermost pairing $\< \ , \ \>_{n+1}$ (resp. $\< \ , \ \>_1$) in \lemref{lem:Ext-6}.
\end{lem}
\begin{proof}
  As in the proof of \lemref{lem:Ext-7}, it suffices to show that
  $\eta_2(\alpha)(\ov{\beta}) = V^n(\alpha \wedge \beta)$ if 
  $\alpha \in Z_n\Fil_D\Omega^{q}_U, \
  {\beta} \in \Fil_{D'}\Omega^{q'}_U$ and $\ov{\beta}$ is the image of $\beta$
  in $\sG_n$. To this end, we write $\alpha = F^n(\ov{\alpha'})$ with
  $\alpha' \in \Fil_D W_{n+1}\Omega^q_U$ and $\ov{\alpha'}$ the image of $\alpha'$
  in $\coker(\ov{p})$. We then get
  $\eta_2(\alpha)(\ov{\beta}) = \alpha' \wedge V^n(\beta) =
  V^n(F^n(\alpha') \wedge \beta) = V^n(\alpha \wedge \beta)$, and 
  this completes the proof.
\end{proof}

For complexes of sheaves of $\sO_Z$-modules $\sA^\bullet$ and
$\sB^\bullet$ on a locally ringed space $(Z, \sO_Z)$, we let
${\bf R}\sHom_{\sO_Z}(\sA^\bullet, \sB^\bullet)$ denote the derived internal
hom complex (defined, for instance, in \cite[\S~6(A)]{Spaltenstein},
where it is denoted by ${\bf R}\sHom^\bullet_{\sO_Z}(\sA^\bullet, \sB^\bullet)$).
For any coherent sheaf $\sF$ of $W_m\sO_X$-modules, we let
$D^0_m(\sF) = \sHom^\bullet_{W_m\sO_X}(\sF, W_m\Omega^d_X)$, and for any bounded complex
$\sF^\bullet$ of coherent sheaves of $W_m\sO_X$-modules, we let
$D_m(\sF^\bullet) = {\bf R}\sHom_{W_m\sO_X}(\sF^\bullet, W_m\Omega^d_X)$.

We consider the pairings $\ov{\fF}^n_* B_n\Fil_D\Omega^{q}_U \times
\ov{\fF}^n_* \frac{\Fil_{D'}\Omega^{q'}_U}{Z_n\Fil_{D'}\Omega^{q'}_U} \to
\Omega^d_X$ and
$\ov{\fF}^n_* Z_n\Fil_D \Omega^q_U \times
\ov{\fF}^n_* \frac{\Fil_{D'} \Omega^{q'}_U}{B_n\Fil_{D'} \Omega^{q'}_U}$ \\
$\to \Omega^d_X$,
given by $(a, \ov{b}) \mapsto C^n(a \wedge b)$ and
$(a', \ov{b'}) \mapsto C^n(a' \wedge b')$ (cf. \lemref{lem:Ext-6}).
We let $\wt{\eta}_1 \colon \ov{\fF}^n_* B_n\Fil_D\Omega^{q}_U \to
D^0_1\left(\ov{\fF}^n_* \frac{\Fil_{D'}\Omega^{q'}_U}{Z_n\Fil_{D'}\Omega^{q'}_U}\right)$
and $\wt{\eta}_2 \colon \ov{\fF}^n_* Z_n\Fil_D \Omega^q_U \to
D^0_1\left(\ov{\fF}^n_* \frac{\Fil_{D'} \Omega^{q'}_U}{B_n\Fil_{D'} \Omega^{q'}_U}\right)$
denote the induced maps.

\begin{lem}\label{lem:Ext-9}
  For $n \in \N_0$ and $i =1,2$, the map $\wt{\eta}_i$ is an isomorphism
  of $\sO_X$-modules.
\end{lem}
\begin{proof}
  We first assume that $\wt{\eta}_2$ is an isomorphism and show that $\wt{\eta}_1$
  is an isomorphism. For this,
we look at the commutative diagram of sheaves of $\sO_X$-modules
  \begin{equation}\label{eqn:Ext-9-0}
    \xymatrix@C1.5pc{
      0 \ar[r] & \ov{\fF}^n_* B_n\Fil_D\Omega^{q}_U \ar[r] \ar[d]_-{\wt{\eta}_1} &
      \ov{\fF}^n_* \Fil_D\Omega^{q}_U \ar[r] \ar[d] &
      \ov{\fF}^n_* \frac{\Fil_{D} \Omega^{q}_U}{B_n\Fil_{D} \Omega^{q}_U} \ar[r]
      \ar[d]^-{\wt{\eta}^{\vee}_2} & 0 \\
      0 \ar[r] &
   D^0_1\left(\ov{\fF}^n_* \frac{\Fil_{D'}\Omega^{q'}_U}{Z_n\Fil_{D'}\Omega^{q'}_U}\right)
   \ar[r] &
   D^0_1\left(\ov{\fF}^n_* \Fil_{D'}\Omega^{q'}_U\right) \ar[r] &
   D^0_1\left(\ov{\fF}^n_* Z_n\Fil_{D'}\Omega^{q'}_U\right) \ar[r] & 0.}
    \end{equation}

  Since all the sheaves involved in the lemma are coherent
  locally free $\sO_X$-modules by \lemref{lem:Can-filt-5},
  our assumption about $\wt{\eta}_2$ (applied to $D'$ in place of $D$)
  implies that the right vertical arrow in
  ~\eqref{eqn:Ext-9-0} is an isomorphism. The middle vertical arrow is
  an isomorphism, as shown in the Step~2 of the proof of \cite[Thm.~4.1.1]{JSZ}).
  Since the two rows of ~\eqref{eqn:Ext-9-0} are exact, it follows that its left
  vertical arrow is an isomorphism.

  We now show that $\wt{\eta}_2$ is an isomorphism by induction on $n \ge 0$.
  The case $n =0$ is the middle vertical arrow of ~\eqref{eqn:Ext-9-0} and
 the case $n=1$ is shown in the Step~2 of the proof of \cite[Thm.~4.1.1]{JSZ}).
 Using \cite[Cor.~6.6]{KM-3}, the argument for the general case becomes
 identical to the one for the $D = \emptyset$ case, given in the proof of
 \cite[Thm.~1.4]{Nygaard} (especially, look at the diagram at the bottom of p.~251).
 \end{proof}

\begin{remk}\label{remk:Ext-9-1}
  The reader may have noticed that, in the above proof, we needed to interchange
  $D$ and $D'$ in order to show that the right vertical arrow in \eqref{eqn:Ext-9-0}
  is bijective. This was possible because, throughout our discussion, we have
  allowed $D$ to be arbitrary, rather than restricting it to be effective.
\end{remk}

We shall now prove the main results of this section.
The following result is an extension of \cite[Thm.~II.2.2]{Ekedahl} to
Hodge--Witt sheaves with modulus.

\begin{thm}\label{thm:Ek-Main-1}
  Let $q \ge 0, m \ge 1$ be integers and $D \in \Div_E(X)$. Let
  $q' = d-q$ and $D' = -D-E$. Then the multiplication operation of
  $W_m\Omega^\bullet_U$ induces canonical isomorphisms
  \[
  \Fil_D W_m\Omega^q_U \xrightarrow{\cong}  D_m(\Fil_{D'} W_m\Omega^{q'}_U) \ \
  \mbox{and} \ \
  \Fil_{D'} W_m\Omega^{q'}_U \xrightarrow{\cong} D_m(\Fil_D W_m\Omega^q_U).
  \]
\end{thm}
\begin{proof}
  It is enough to prove the first isomorphism as the other one
  then follows by swapping $D$ with $D'$.
  Using the adjointness of the pair of functors $(\Re_{X *}, \Re^{!}_X)$ and the
  isomorphism $\Re^{!}_S W_m(k) \cong W_{m-1}(k)$, one quickly checks that
  $\Re_{X *} D_{m-1}(\sF^\bullet) \cong D_m(\Re_{X *} \sF^\bullet)$ for any
  bounded complex $\sF^\bullet$ of coherent $W_{m-1}\sO_X$-modules.
  Using this isomorphism, we get a commutative diagram of exact triangles
  \begin{equation}\label{eqn:Ek-Main-1-0}
    \xymatrix@C1.5pc{
      0 \ar[r] & \Re_{X *} \Fil_{D/p} W_{m-1}\Omega^q_U \ar[r]^-{\ov{p}} \ar[d] &
      \Fil_{D} W_{m}\Omega^q_U \ar[r] \ar[d] & \coker(\ov{p}) \ar[r] \ar[d] & 0 \\
      0 \ar[r] & \Re_{X *}D_{m-1}(\Fil_{{D'}/p}W_{m-1}\Omega^{q'}_U) \ar[r]^-{R^\vee}
      & D_m(\Fil_{D'}W_m\Omega^{q'}_U) \ar[r] & D_m(\Ker(R_m)) \ar[r] & 0,}
  \end{equation}
  in which $R_m \colon \Fil_{D'}W_m\Omega^{q'}_U \surj
  \Re_{X *} \Fil_{{D'}/p}W_{m-1}\Omega^{q'}_U$ and
  the vertical arrows are induced by the duality pairings.

  Using induction on $m$, the proof of the theorem is reduced 
  to showing that the right vertical
  arrow in this diagram is an isomorphism. Using Lemmas~\ref{lem:Can-filt-5},
~\ref{lem:Ext-0}, ~\ref{lem:Ext-6},
  ~\ref{lem:Ext-7} and ~\ref{lem:Ext-8} (note that composition with $C^n$
  in the right columns of ~\eqref{eqn:Ext-7-0} and ~\eqref{eqn:Ext-8-0} are
  isomorphisms), the proof is further reduced to showing that
  for every $n \ge 0$, the duality maps $\wt{\eta}_1$ (with $q$ replaced by
  $q+1$) and $\wt{\eta}_2$ are isomorphisms.
   But this follows from \lemref{lem:Ext-9}.
  \end{proof}

The following result is an extension of \cite[Cor.~II.2.2.23]{Ekedahl} to
Hodge--Witt cohomology with modulus. 

\begin{cor}\label{cor:Ek-Main-2}
Assume in \thmref{thm:Ek-Main-1} that $X$ is proper over $k$ and let $i \in \N_0$.
Then the multiplication operation of
$W_m\Omega^\bullet_U$ induces a canonical isomorphism of finitely generated
$W_m(k)$-modules
\[
H^i_\zar(X, \Fil_D W_m\Omega^q_U) \xrightarrow{\cong}
\Hom_{W_m(k)}(H^{d-i}_\zar(X, \Fil_{D'} W_m\Omega^{q'}_U), W_m(k)).
\]
\end{cor}
\begin{proof}
  This is a direct consequence of \thmref{thm:Ek-Main-1}, the Grothendieck
  duality for the proper map $f_m \colon W_mX \to W_mS$
  (cf. \cite[Thm.~4.30]{Lee-Nakayama}, \cite[Tag.~0FVU]{SP}), and the injectivity
  of $W_m(k)$ as a $W_m(k)$-module. The latter fact easily follows once one
  observes that $W_m(k)$ is an Artinian local ring in which
  the annihilator of the ideal $V(W_{m-1}(k))$ is
  $V^{m-1}(k)$ which is a 1-dimensional vector space over $k$.
\end{proof}

For $D = \sum_i n_iE_i \ge 0$ and $q \in \N_0$,
we let $W_m\Omega^q_{X|D} = \Ker\left(W_m\Omega^q_{X} \surj
W_m\Omega^q_{D}\right)$ (cf. \cite[Defn.~3.5]{KM-3}) and $W_m\Omega^q_{(X, -D)} =
\Ker\left(W_m\Omega^q_{X} \to \ \stackrel{r}{\underset{i =1}\bigoplus}
W_m\Omega^q_{n_iE_i}\right)$ (cf. \cite[Defn.~8.1]{Ren-Rulling}).
It follows from \lemref{lem:Log-fil-4}(8) that
$\Fil_{-E} W_m\Omega^q_U \subset W_m\Omega^q_{X|E} \subset W_m\Omega^q_{(X, -E)}$.
As an application of \thmref{thm:Ek-Main-1}, we obtain the following
stronger relation between these sheaves.

\begin{cor}\label{cor:EK-Main-4}
  The inclusions $\Fil_{-E} W_m\Omega^q_U \subset W_m\Omega^q_{X|E} \subset
  W_m\Omega^q_{(X, -E)}$ are isomorphisms.
\end{cor}
\begin{proof}
  For $q = 0$, this can be checked directly. For $q \ge 1$,
  we look at the commutative diagram
  \begin{equation}\label{eqn:EK-Main-4-0}
    \xymatrix@C1.5pc{
      \Fil_{-E} W_m\Omega^q_U  \ar[r] \ar[d] &  D_m(\Fil_0 W_m\Omega^{d-q}_U) \ar[d] \\
      W_m\Omega^q_{(X, -E)} \ar[r] &  D_m(W_m\Omega^{d-q}_{(X,E)}),}
  \end{equation}
  where $W_m\Omega^{d-q}_{(X,E)}$ is the sheaf considered in \cite[\S~6.2]{Ren-Rulling},
  and the horizontal arrows are the duality homomorphisms.
  The right vertical arrow is an isomorphism by \lemref{lem:Log-fil-4} and
  the comments following [op. cit., Cor.~1]. The top horizontal arrow is
  an isomorphism by \thmref{thm:Ek-Main-1}, and the bottom horizontal arrow is an
  isomorphism by [op. cit., Thm.~1]. We conclude that
  the left vertical arrow is an isomorphism of sheaves.
\end{proof}

\section{Filtered de Rham--Witt on the RP site}\label{sec:RP}
In this section, we extend the filtered de Rham--Witt complex to
Kato's relatively perfect site of a regular $F$-finite scheme and prove
some results about this extension which will be needed in the proof of
\thmref{thm:Main-1}.

\subsection{de Rham--Witt complex on the RP site}\label{sec:DRRP}
This subsection is mostly the recollection of relatively perfect schemes over a
regular scheme and their de Rham--Witt complex from \cite{Kato-duality-1} and
  \cite{Kato-duality-2}.
  
Let $X$ be a connected, Noetherian, regular $F$-finite scheme.
 Recall from \cite[\S~1]{Kato-duality-1}
that a scheme $T$ over $X$ is relatively perfect over $X$ (i.e., $T$ is a relatively
perfect $X$-scheme) if the diagram of absolute Frobenius morphisms
\begin{equation}\label{eqn:Rel-perf-0}
  \xymatrix@C1.7pc{
    T \ar[r]^-{F_T} \ar[d]_-{f} & T \ar[d]^-{f} \\
    X \ar[r]^-{F_X} & X}
\end{equation}
is Cartesian, where $f \colon T \to X$ is the projection map.

It is easy to check that a relatively perfect $X$-scheme is reduced.
By \cite[Lem.~1.3, Prop.~5.2]{Kato-duality-1}, an {\'e}tale $X$-scheme is
relatively perfect, and a relatively perfect $X$-scheme is flat and formally
{\'e}tale over $X$.
By \cite{Kunz} (see also \cite[Tag~0EBU]{SP}), a relatively perfect $X$-scheme is a
regular scheme if and only if it is locally Noetherian. We let $X_{\rp}$ be the
site whose underlying category is that of relatively perfect $X$-schemes and whose
underlying Grothendieck topology is the {\'e}tale topology.
We let $\nu \colon X_{\rp} \to X_{\et}$ denote the morphism of sites induced by the
canonical inclusion of the category of {\'e}tale $X$-schemes into the category of
relatively perfect $X$-schemes.

As a relatively perfect $X$-scheme
may be badly behaved in general (it may not even be Noetherian),
the following result plays a key role
in the study of the de Rham--Witt complex of such a scheme.

\begin{prop}\label{prop:RP-0}
  Let $f \colon T \to X$ be a relatively perfect morphism and let
  $W_mf \colon W_mT \to W_mX$ be the induced morphism. We then have the following.
  \begin{enumerate}
  \item
    $W_mf \colon  W_mT \to W_mX$ is a relative perfect and flat morphism.
    In particular, it is formally {\'e}tale.
  \item
    The diagrams
    \begin{equation}\label{eqn:RP-0-0}
      \xymatrix@C1.6pc{
        W_{m-1}T \ar[r]^-{\phi} \ar[d]_-{W_{m-1}f} & W_mT \ar[d]^-{W_mf} & &
          W_{m}T \ar[r]^-{\ov{\fF}} \ar[d]_-{W_mf} & W_mT \ar[d]^-{W_mf} \\
     W_{m-1}X \ar[r]^-{\phi}  & W_mX  & &
     W_mX \ar[r]^-{\ov{\fF}} & W_mX}
    \end{equation}
    are Cartesian for $\phi \in \{\Re, \fF\}$.
  \item
    If $f \colon X' \to X$ is any morphism and $T' = T\times_X X'$, then
    $T' \to X'$ is relatively perfect and $W_mT' \cong W_mT \times_{W_mX} W_mX'$.
  \end{enumerate}
\end{prop}
\begin{proof}
  See \cite[Prop.~2.17]{Shiho} for all claims except that $W_mf$ and $T' \to X'$
  are relatively perfect morphisms. The latter claims follow from
  \cite[Lem.~1.2, 1.3]{Kato-duality-2}.
  Note also that $\fF =  \Re \circ \ov{\fF}$ so that item (2) for
  $\Re$ and $\ov{\fF}$ implies the same for $\fF$.
\end{proof}

Recall that for any scheme $T$, the $p$-typical de Rham--Witt complex
(cf. \cite{Illusie}) $W_m\Omega^\bullet_T$ is the universal Witt complex on $T$.
When $T$ is a relatively perfect $X$-scheme, its de Rham--Witt complex has
the following properties.

\begin{lem}\label{lem:RP-1}
  Let $f \colon T \to X$ be a relatively perfect morphism and let
  $T_U = T \times_X U$ for any open $U \subset X$. Then we have the following.
  \begin{enumerate}
  \item
    The base change map $W_m\Omega^q_X \otimes_{W_m\sO_X} (W_mf)_*W_m\sO_T \to
    (W_mf)_*W_m\Omega^q_T$ is an isomorphism of $(W_mf)_*W_m\sO_T$-modules.
  \item
    The canonical map $W_m\Omega^q_X \to (W_mf)_*W_m\Omega^q_T$ is injective.
  \item
    The canonical map $W_m\Omega^q_T \to (W_mj^T)_*W_m\Omega^q_{T_U}$ is injective,
    where $j^T \colon T_U \inj T$ is the inclusion. 
  \end{enumerate}
\end{lem}
\begin{proof}
The item (1) is precisely \cite[Prop.~2.18]{Shiho}. To prove item (2), we look
  at the commutative diagram
  \begin{equation}\label{eqn:RP-1-0}
    \xymatrix@C1.7pc{
 W_m\Omega^q_X  \ar[r] \ar[d] &   (W_mf)_*W_m\Omega^q_T \ar[d] \\
 W_m\Omega^q_{\eta} \ar[r] & (W_mf')_*W_m\Omega^q_{T_\eta},}
  \end{equation}
  where $f' \colon T_\eta \to \Spec(k(\eta))$ is projection map for the generic
  fiber of $f$ and the vertical arrows are the restriction maps.

  Since the left vertical arrow in ~\eqref{eqn:RP-1-0} is injective
  (cf. \cite[Prop.~2.8]{KP-Comp}), it suffices to show that the
  bottom horizontal arrow is injective. We can therefore assume that $X$ is the
  spectrum of a field.
  In particular, we can assume that $f$ is faithfully flat. This implies by
  virtue of \propref{prop:RP-0} that $W_mf \colon W_mT \to W_mX$ is also faithfully
  flat. In this case, the result follows from \cite[Thm.~7.5]{Matsumura},
  \propref{prop:RP-0} and item (1). Item (3) follows easily from
  by comparing item (1) for $X$ and $U$, and using \propref{prop:RP-0}.
\end{proof}

We shall use the following results about relatively perfect morphisms.

\vskip.2cm

\begin{lem}\label{lem:RP-20}
  Let $f \colon X' \to X$ be a relatively perfect morphism and let
  $E' = E \times_X X'$. Assume that $X'$ is Noetherian and regular.
  Then $(X', E')$ is an snc-pair and $f \colon (X', E') \to (X,E)$ is a
  morphism of snc-pairs.
\end{lem}
\begin{proof}
  By looking at the complement of all but one connected component of $X'$ and $E'$,
  we can
assume that both $X'$ and $E'$ are connected. It follows from \cite{Kunz-1} and
\cite[Lem.~7.2]{Kato-89} that the Krull dimension of $X'$ is finite and coincides
with $\dim(X)$. We let this number be $d$. It follows from \propref{prop:RP-0}
that $E' \to E$ is a flat and relatively perfect morphism of connected Noetherian
schemes of dimension $d-1$. Let $\{E'_1, \ldots , E'_s\}$ be the set of
irreducible components of $E'$.
To prove the first part of the lemma, we only need to show that
$E'$ is an snc divisor on $X'$. By \cite[Lem.~8.1]{KM-2}, it suffices to show
that $E'_{i_1} \bigcap \cdots \bigcap E'_{i_m}$ is
a disjoint union of regular schemes of dimension $d-m$ if $1 \le i_1 < \cdots <
i_m \le s$.

We let $F_i = E_i \times_X X'$. Since $F_i$ is relatively perfect over $E_i$, and the
latter is regular, we see that $F_i$ is a disjoint union of connected regular divisors
on $X'$. Furthermore, each irreducible component of $F_i$ is an irreducible
component of $E'$. Since $E'$ is connected, we conclude that each $E'_i$ is regular
and no two components of $E'$ lie over one component of $E$. We now let
$Z = E_{i_1} \bigcap \cdots \bigcap E_{i_m}$ and
$Z' = E'_{i_1} \bigcap \cdots \bigcap E'_{i_m}$.
We can assume after a permutation that $f(E'_{i_j})$ is an open subscheme of
$E_{i_j}$ (this uses flatness of $E'_{i_j}$ over $E_{i_j}$) for
$1 \le i_1 < \cdots < i_m \le s$.

If $Z$ is an empty set, then so is $Z'$ and there is nothing to prove.
We therefore assume that $Z \neq \emptyset$. Let $W$ be an irreducible
component of $Z$ and $W' = f^{-1}(W)$. Since $W$ is regular and
$W'$ is relatively perfect over $W$, it follows that $W'$ is a disjoint
union of regular closed subschemes of $X'$. It is also clear by induction
on $m$ that each irreducible component of $W'$ has dimension $\dim(W) = d-m$. We
conclude from this that $Z'$ is either empty or is a disjoint union of connected
and regular closed subschemes of $X'$ of dimension $d-m$. This completes the proof
of the first part of the lemma. The second part is clear.
\end{proof}

Recall from \cite[\S~2]{Kelly-Morrow} that an $\F_p$-algebra $A$ is called
Cartier smooth if the following hold.
\begin{enumerate}
\item
  $\Omega^1_A$ is a flat $A$-module.
\item
  $H_n(\L_{A/{\F_p}}) = 0$ for $n > 0$, where $\L_{A/{\F_p}}$ is the cotangent
  complex of $A$.
\item
  The inverse Cartier map $\ov{F} \colon \Omega^n_A \to H^n(\Omega^\bullet_A)$
  is an isomorphism for each $n \ge 0$.
\end{enumerate}
An $\F_p$-scheme is called Cartier smooth if it is locally Cartier smooth.

\begin{lem}\label{lem:RP-7*}
  Let $Z$ be a Cartier smooth $\F_p$-scheme and let $g \colon Z' \to Z$
  be a flat and relatively perfect morphism of schemes.  Then $Z'$ is a Cartier
  smooth $\F_p$-scheme.
\end{lem}
\begin{proof}
  We can assume $Z = \Spec(A)$ and $Z' = \Spec(B)$.
  Since $B$ is relatively perfect over $A$, it is formally {\'e}tale over $A$.
  In particular, the canonical map of $B$-modules $\Omega^n_{A} \otimes_{A} {B}
  \to \Omega^n_{B}$ is an isomorphism for every $n \ge 0$, and so is the
  map of complexes of $B$-modules
  $F_*(\Omega^\bullet_{A}) \otimes_{A} {B} \to F_*(\Omega^\bullet_{B})$, where
  $F$ is the absolute Frobenius. It follows from this that
  conditions (1) and (3) are satisfied for $B$.
The condition (2) follows from the exact triangle
  (cf. \cite[Tag~08QX]{SP})
  \begin{equation}\label{eqn:RP-7-0}
    \L_{A/{\F_p}} \stackrel{\L}{\otimes}_A {B} \to \L_{{B}/{\F_p}} \to
    \L_{{B}/{A}}
    \end{equation}
  in the derived category of $B$-modules together with the
  fact that $B$ is flat over $A$, and
  $\L_{{B}/{A}} = 0$  (cf. \cite[Prop.~3.17]{Rojas-Stabler}).
\end{proof}
  
\begin{cor}\label{cor:RP-7}
  Let $Z$ be a locally Noetherian regular $\F_p$-scheme and let $g \colon Z' \to Z$
  be a relatively perfect morphism of schemes. Then $Z'$ is a Cartier
  smooth $\F_p$-scheme.
\end{cor}
\begin{proof}
Under the hypotheses of the corollary, $Z'$ is flat over $Z$ by
  \cite[Prop.~5.2]{Kato-duality-1}. Furthermore, it follows using the
  N{\'e}ron-Popescu approximation that $Z$ is a Cartier smooth $\F_p$-scheme.
We can therefore apply \lemref{lem:RP-7*}.
\end{proof}

\subsection{Filtered de Rham--Witt complex on the RP site}\label{sec:FDR-RP}
Let $X$ be a connected, Noetherian, regular $F$-finite scheme and let $E \subset X$ be
an snc divisor with irreducible components $E_1, \ldots , E_r$.  Let
$j \colon U \inj X$ be the inclusion of the complement of $E$ in $X$.
Let $d$ denote the $p$-rank of $X$.

\begin{defn}\label{defn:RP-2}
  Let $f \colon X' \to X$ be a relatively perfect morphism and
  $U' = X' \times_X U$. Given $q \in \N_0$ and $m \in \N$, we let
  \[
  W_m\Omega^q_{X'}(\log E) := f^*(W_m\Omega^q_{X}(\log E)) \cong
  W_m\Omega^q_{X}(\log E) \otimes_{W_m\sO_X} W_m\sO_{X'}.
  \]
Given $D \in \Div_E(X)$, we let
\[
\Fil_D W_m \Omega^q_{U'} :=  f^*(\Fil_D W_m \Omega^q_U) \cong
\Fil_D W_m \Omega^q_U \otimes_{W_m\sO_X} W_m\sO_{X'}.
\]

Here, $M \otimes_{W_m\sO_X}W_m\sO_{X'}$ is a shorthand notation for $(W_mf)^{-1}M \otimes_{(W_mf)^{-1}W_m\sO_X}W_m\sO_{X'}$ for any quasi-coherent sheaf $M$ on $W_mX$. It is clear from Proposition~\ref{prop:RP-0} and ~\lemref{lem:RP-1} that $\Fil_D W_m \Omega^q_{U'}$ is a subsheaf of $W_m\sO_{X'}$-modules inside $(W_mj')_* W_m \Omega^q_{U'}$ and $\Fil_0 W_m \Omega^q_{U'} = W_m\Omega^q_{X'}(\log E)$, where $j' \colon U' \inj X'$ is the inclusion. We shall write $\Fil_D W_1 \Omega^q_{U'}$ (resp. $W_1\Omega^q_{X'}(\log E)$) as $\Fil_D \Omega^q_{U'}$ (resp. $\Omega^q_{X'}(\log E)$).
We let $Z_1\Fil_D W_{m} \Omega^q_{U'} := \Fil_D W_{m} \Omega^q_{U'} \bigcap
j'_* Z_1W_m\Omega^q_{U'}$, where $ Z_1W_m\Omega^q_{U'}:=\Ker(F^{m-1}d:W_m\Omega^q_{U'}\to \Omega^{q+1}_{U'})$. Note that $\Fil_D W_{m} \Omega^q_{U'} = W_{m} \Omega^q_{X'}$
when $E = \emptyset$.
In the remainder of this paper, we shall write $u_* W_m\Omega^q_{V}$ simply as $W_m\Omega^q_{V}$
for any open immersion $u \colon V \inj T$ of schemes.
\end{defn}

\begin{lem}\label{lem:RP-3}
Let $f \colon X' \to X$ be a relatively perfect morphism and
let $j' \colon U' = X' \times_X U \inj X'$ be the inclusion. Let $D \in \Div_E(X)$.
As subsheaves of the de Rham--Witt complex
$(j'_* W_m \Omega^\bullet_{U'}, F,V, R, d)$ of $U'$,
$\Fil_D W_m \Omega^\bullet_{U'}$ satisfy the following.
\begin{enumerate}
\item
  $\Fil_D W_m \Omega^q_{U'}$ is a sheaf of finitely generated $W_m\sO_{X'}$-modules.
  \item
    $\Fil_D W_m \Omega^q_{U'} = \Fil_D W_m \Omega^q_{U}\cdot W_m\sO_{X'}$
    as a subsheaf of $W_m \Omega^q_{U'}$.
  \item
  $\Fil_D W_m \Omega^q_{U'} \cdot \Fil_{D'} W_m \Omega^{q'}_{U'} \subset
  \Fil_{D+D'} W_m \Omega^{q+q'}_{U'}$.
\item
  $\Fil_D \Omega^q_{U'} \cong \Omega^q_{X'}(\log E)(D')$, where
  $D' = f^*(D)$.
  \item
    ${\varinjlim}_{n} \Fil_{nE} W_m \Omega^\bullet_{U'} \xrightarrow{\cong}
  W_m \Omega^\bullet_{U'}$.
  \item
    $\Fil_D W_m \Omega^q_{U'} = \Fil_D W_m \sO_U \cdot W_m\Omega^q_{X'}(\log E) +
    d(\Fil_D W_m \sO_U) \cdot W_m\Omega^{q-1}_{X'}(\log E)$.
  \item
    The map $\theta \colon {\phi}_*\Fil_D W_m \Omega^q_{U} \otimes_{W_{m+1}\sO_X}
    W_{m+1}\sO_{X'} \to  {\phi}_* \Fil_D W_m \Omega^q_{U'}$, defined by
    $\theta(w \otimes b) = w\phi^*(b)$,
    is an isomorphism of sheaves of $W_{m+1}\sO_{X'}$-modules, where $\phi \in \{\Re, \fF\}$.
  \item
    As $W_m\sO_{X'}$-submodule of $W_m\Omega^d_{U'}$, one has
    $\Fil_{-E} W_m \Omega^d_{U'} = W_m \Omega^d_{X'}$.
\end{enumerate}
\end{lem}
\begin{proof}
Each of the items from (1) to (6) follows directly from the definition of
$\Fil_D W_m \Omega^q_{U'}$ and the corresponding property of
$\Fil_D W_m \Omega^q_{U}$ using Proposition~\ref{prop:RP-0} and ~\lemref{lem:RP-1}
(cf. \cite[Lem.~3.6, 3.10]{KM-3}).
Item (7) follows from \propref{prop:RP-0} using the flat base change theorem
\cite[Tag 02KH]{SP}. Item (8) follows directly from
 \thmref{thm:Global-version}(13) and \lemref{lem:RP-1}.
 \end{proof}

Let $f \colon X' \to X$ be a relatively perfect morphism and
let $j' \colon U' = X' \times_X U \inj X'$ be the inclusion.
Then $X'$ is endowed with the log structure 
given by the canonical inclusion of the {\'e}tale sheaves of monoids
$\alpha_{X'} \colon \sM_{X'} = \sO_{X'} \bigcap j'_*\sO^\times_{U'} \inj \sO_{X'}$.
We let $X'_\log = (X', \sM_{X'})$ denote the resulting log scheme. We then have a 
strict morphism of log schemes $f \colon X'_\log \to X_\log$.
We let $W_m\Omega^\bullet_{X'_\log}$ denote the de Rham--Witt complex of the
log scheme $X'_\log$ (cf. \cite{Hyodo-Kato}).
We refer to \cite[\S~1.1]{Geisser-Hesselholt-JAMS} for the definition of a
Witt complex over a log scheme.

\begin{prop}\label{prop:RP-4}
  $W_m\Omega^\bullet_{X'}(\log E)$ is a Witt complex over $X'_\log$ and there is a
  canonical isomorphism of Witt complexes $f^* \colon W_m\Omega^\bullet_{X'}(\log E)
  \xrightarrow{\cong} W_m\Omega^\bullet_{X'_\log}$ over $X'_\log$.
\end{prop}
\begin{proof}
  By \lemref{lem:RP-1}, we have $W_m\Omega^\bullet_{X'}(\log E) \inj
  W_m\Omega^\bullet_{U'}$. We first show that this is an inclusion of
  Witt complexes over $X'$. For this, we recall that
  $W_m\Omega^q_{X'}(\log E) = W_m\Omega^q_{X}(\log E)\cdot W_m\sO_{X'}
  \subset W_m\Omega^\bullet_{U'}$. Since $d(W_m\sO_{X'}) \subset
  W_m\Omega^1_{X'} \subset W_m\Omega^1_{X'}(\log E)$, it follows from
  \lemref{lem:RP-3}(3) that $d(W_m\Omega^q_{X'}(\log E)) \subset
  W_m\Omega^{q+1}_{X'}(\log E)$. A similar argument, using the multiplicativity of
  $F$, shows that $F(W_{m+1}\Omega^q_{X'}(\log E)) \subset 
  W_{m}\Omega^q_{X'}(\log E)$. To show that $W_m\Omega^\bullet_{X'}(\log E)$ is
  preserved by $V$, it suffices to show, using \lemref{lem:RP-3}(7) (for $\phi=\fF$) that
  $V(wF(b)) \in W_{m+1}\Omega^q_{X'}(\log E)$ if $w \in W_{m}\Omega^q_{X}(\log E)$
  and $b \in W_{m+1}\sO_{X'}$. But this is clear since
$V(wF(b)) = V(w)b$. 

Next, to endow  $W_\star\Omega^\bullet_{X'}(\log E)$ with the structure of a Witt complex
over $X'_\log$, we let
$\dlog \circ [\cdot]_m \colon \sM_{X'} \to W_m\Omega^1_{X'}(\log E)$ be the map which
is locally defined by $a \mapsto \dlog[a]_m$.
This is well-defined because locally we can write $a = a'\pi^{n}$ where
$a' \in \sO^\times_{X'}$ and $n \in \Z$, and $W_m\Omega^1_{X'} =
W_m\Omega^1_{X} \cdot W_m\sO_{X'} \subset W_m\Omega^1_{X'}(\log E)$.
  We take $\lambda \colon (W_\star \sO_{X'}, \sM_{X'}) \to
  (W_\star\Omega^0_{X'}(\log E), \sM_{X'})$ to
  be the identity map. Since $W_\star\Omega^\bullet_{X'}(\log E)$ is preserved by the
  Frobenius and Verchiebung operators of $W_\star\Omega^\bullet_{U'}$, as shown
  above, it follows that $W_\star\Omega^\bullet_{X'}(\log E)$ is a Witt complex over
  $X'_{\log}$.

  By virtue of the above, we get a unique morphism of Witt complexes
  $\phi_{X'} \colon W_\star\Omega^\bullet_{X'_{\log}} \to  W_\star\Omega^\bullet_{X'}(\log E)$
  over $X'_{\log}$.
  On the other hand, $f$ induces a morphism of Witt complexes
  \begin{equation}\label{eqn:RP-4-0}
 W_\star\Omega^\bullet_{X}(\log E) \xrightarrow{\cong} W_\star\Omega^\bullet_{X_{\log}}
 \xrightarrow{f^*} f_* W_\star\Omega^\bullet_{X'_{\log}}
  \end{equation}
  over $X_\log$, where the first isomorphism follows from \cite[Thm.~2.8]{KM-3}.
  By adjointness, this yields a unique morphism
  of $W_m\sO_{X'}$-modules $\psi_{X'} \colon W_\star\Omega^\bullet_{X'}(\log E) \to
  W_\star\Omega^\bullet_{X'_{\log}}$. We have shown above that the left hand side of this  
  map is a Witt complex over $X'_{\log}$, and it follows from ~\eqref{eqn:RP-4-0} that
  $\psi_{X'}$ is a morphism of Witt complexes over $X'_{\log}$ (see the proof of
  \cite[Thm.~C]{Hesselholt-Acta}). Since
  $\phi_{X'} \circ \psi_{X'}$ is identity, the universal property of
  $W_\star\Omega^\bullet_{X'_{\log}}$ implies that $\phi_{X'}$ and $\psi_{X'}$ are
  inverses of each other. This completes the proof.
\end{proof}

\subsection{Some properties of $\Fil_D W_m \Omega^\bullet_{U'}$}\label{sec:RP-fil-basic}
In the next few lemmas, we describe some properties of the filtered de Rham--Witt
complex on $X_\rp$. We begin with the following properties of
various operators on $\Fil_D W_m \Omega^\bullet_{U'}$.

\begin{lem}\label{lem:RP-5}
  Let $f \colon X' \to X$ be as in \lemref{lem:RP-3}. We then have the following.
  \begin{enumerate}
  \item
  $d(\Fil_D W_m \Omega^q_{U'}) \subset \Fil_D W_m \Omega^{q+1}_{U'}$.
\item
  $F(\Fil_D W_{m+1} \Omega^q_{U'}) \subset \Fil_D W_m \Omega^{q}_{U'}$.
\item
$V(\Fil_D W_{m} \Omega^q_{U'}) \subset \Fil_D W_{m+1} \Omega^{q}_{U'}$.
\item
  $R(\Fil_D W_{m+1} \Omega^q_{U'}) = \Fil_{D/p} W_m \Omega^{q}_{U'}$.
\item
There exists a short exact sequence of sheaves of $W_m\sO_{X'}$-modules
  \[
  0 \to V^{m-1}\left(\Fil_D \Omega^q_{U'}\right) +
  dV^{m-1}\left(\Fil_D \Omega^{q-1}_{U'}\right) \to \Fil_D W_{m} \Omega^q_{U'}
  \xrightarrow{R} \Fil_{D/p} W_{m-1} \Omega^{q}_{U'} \to 0.
  \]
  \end{enumerate}
 \end{lem}
\begin{proof}
The first three items easily follow from \lemref{lem:RP-3}(6) and
  \propref{prop:RP-4} (and its proof). To prove item (4), we base change the
  surjective $W_{m+1}\sO_X$-linear map
  $R \colon \Fil_D W_{m+1} \Omega^q_{U} \surj  \Re_*\Fil_{D/p} W_m \Omega^{q}_{U}$
  (cf. \thmref{thm:Global-version}(1)) to
  $X'$ to get a surjective $W_{m+1}\sO_{X'}$-linear map
  $\Fil_D W_{m+1} \Omega^q_{U} \otimes_{W_{m+1}\sO_X} W_{m+1}\sO_{X'} \surj
  \Re_*\Fil_{D/p} W_m \Omega^{q}_{U} \otimes_{W_{m+1}\sO_X} W_{m+1}\sO_{X'}$.
  We now use \lemref{lem:RP-1} and Definition~\ref{defn:RP-2} to conclude (4).

To prove (5), we look at the exact sequence (cf. \thmref{thm:Global-version}(1))
  \begin{equation}\label{eqn:RP-5-0}
    \fF^{m}_* \Fil_D \Omega^q_{U} \bigoplus  \fF^{m}_* \Fil_D \Omega^{q-1}_{U}
    \xrightarrow{(V^{m-1}, \ dV^{m-1})}
    \fF_*  \Fil_D W_{m} \Omega^q_{U} \xrightarrow{R}
    \fF_* \Re_* \Fil_{D/p} W_{m-1} \Omega^{q}_{U} \to 0.
  \end{equation}
  Since this is an exact sequence of $W_{m+1}\sO_{X}$-modules by \cite[Lem.~5.7]{KM-3},
  we can apply \propref{prop:RP-0} to get an exact sequence of
  $W_{m+1}\sO_{X'}$-modules
  \[
  \left(\fF^{m}_* \Fil_D \Omega^q_{U} \bigoplus
  \fF^{m}_* \Fil_D \Omega^{q-1}_{U}\right)
 \otimes W_{m+1}\sO_{X'}   \xrightarrow{(V^{m-1}, \ dV^{m-1})}
 \fF_*  \Fil_D W_{m} \Omega^q_{U} \otimes W_{m+1}\sO_{X'} \hspace{2cm}
 \]
 \[
 \hspace{8cm} \xrightarrow{R}
    \fF_* \Re_* \Fil_{D/p} W_{m-1} \Omega^{q}_{U} \otimes  W_{m+1}\sO_{X'} \to 0,
    \]
where the tensor products are taken over $W_{m+1}\sO_{X}$.
By \lemref{lem:RP-3}, this 
yields an exact sequence of $W_{m+1}\sO_{X'}$-modules
\begin{equation}\label{eqn:RP-5-1}
\fF^m_* \Fil_D \Omega^q_{U'} \bigoplus  \fF^m_* \Fil_D \Omega^{q-1}_{U'}
    \xrightarrow{(V^{m-1}, \ dV^{m-1})}
    \fF_*  \Fil_D W_{m} \Omega^q_{U'} \xrightarrow{R}
    \fF_* \Re_* \Fil_{D/p} W_{m-1} \Omega^{q}_{U'} \to 0,
\end{equation}
    which proves (5). This completes the proof.
  \end{proof}

\begin{lem}\label{lem:RP-6}
 Let $f \colon X' \to X$ be as in \lemref{lem:RP-3}. We then have the following.
  \begin{enumerate}
  \item
    $Z_1\Fil_D W_{m} \Omega^q_{U'} =
    \Ker\left(F^{m-1}d \colon \Fil_D W_{m} \Omega^q_{U'} \to \Omega^{q+1}_{U'}\right) = F(\Fil_D W_{m+1} \Omega^q_{U'})$.
  \item
 $V^{m-1}\left(\Fil_D \Omega^q_{U'}\right) = \Ker(\Fil_D W_m\Omega^q_{U'} 
    \xrightarrow{F} \Fil_D W_{m-1}\Omega^q_{U'})$.
  \item
    There exists a $W_{m+1}\sO_{X'}$-linear monomorphism
    $\ov{p} \colon \Re_*\Fil_{D/p} W_m\Omega^q_{U'} \inj \Fil_D W_{m+1}\Omega^q_{U'}$
    such that $\ov{p} \circ R = p$.
  \item
    The map $\Fil_D W_{m+1}\Omega^q_{U'} \xrightarrow{F} \Fil_D W_{m}\Omega^q_{U'}$
    induces an isomorphism of $W_{m+1}\sO_{X'}$-modules
    $\ov{F} \colon \Re_*\Fil_{D/p} W_{m}\Omega^q_{U'} \xrightarrow{\cong}
    \frac{\fF_* Z_1\Fil_D W_{m} \Omega^q_{U'}}{\fF_*dV^{m-1}\Fil_D\Omega^{q-1}_{U'}}$. 
  \item
    There exists a Cartier homomorphism $C \colon \fF_*Z_1\Fil_D W_{m} \Omega^q_{U'} \to \Re_*\Fil_{D/p} W_{m}\Omega^q_{U'}$
    of $W_{m+1}\sO_{X'}$-modules  such that $\ov{p} \circ C = V$ and $C \circ F = R$.
  \item
    The map $\fF_*Z_1\Fil_D W_{m} \Omega^q_{U'} \xrightarrow{C} \Re_*\Fil_{D/p} W_{m}\Omega^q_{U'}$
    induces an isomorphism of $W_{m+1}\sO_{X'}$-modules 
    $\ov{C} \colon \frac{\fF_*Z_1\Fil_D W_{m} \Omega^q_{U'}}{\fF_*dV^{m-1}\Fil_D\Omega^{q-1}_{U'}} \xrightarrow{\cong} \Re_*\Fil_{D/p} W_{m}\Omega^q_{U'}$ such that $\ov{F} \circ \ov{C}
    = \id$ and $\ov{C} \circ \ov{F} = \id$.
\end{enumerate}
\end{lem}
\begin{proof}

The first equality in item $(1)$ follows from the definition of $Z_1\Fil_D W_{m} \Omega^q_{U'}$ (see \defref{defn:RP-2}). To prove the second equality, we use \thmref{thm:Global-version}(2) to get an exact sequence
  \begin{equation}\label{eqn:RP-6-0}
\Fil_D W_{m+1} \Omega^q_{U} \xrightarrow{F}
\fF_*\Fil_D W_{m} \Omega^q_{U} \xrightarrow{F^{m-1}d} \fF^m_*\Omega^{q+1}_U
\end{equation}
By ~\eqref{eqn:Mor-W_mX*-3}, this is an exact sequence of $W_{m+1}\sO_X$-modules.
In particular, 
\[
\Fil_D W_{m+1} \Omega^q_{U} \otimes W_{m+1}\sO_{X'} \xrightarrow{F\otimes \id}
\fF_*\Fil_D W_{m} \Omega^q_{U}\otimes W_{m+1}\sO_{X'} \xrightarrow{F^{m-1}d \otimes \id}  \fF^m_*\Omega^{q+1}_U\otimes W_{m+1}\sO_{X'}
\]
is an exact sequence of $W_{m+1}\sO_{X'}$-modules, where the tensor products are
taken over $W_{m+1}\sO_X$. By Proposition~\ref{prop:RP-0} and ~\lemref{lem:RP-3}(7),
we can write the latter sequence as
\begin{equation}\label{eqn:RP-6-1}
 \Fil_D W_{m+1} \Omega^q_{U'} \xrightarrow{F}  \fF_*\Fil_D W_{m} \Omega^q_{U'} \xrightarrow{F^{m-1}d} 
   \fF^m_*\Omega^{q+1}_{U'}.
    \end{equation}
This proves the second equality of item (1). Moreover, we have the following isomorphism: 
\begin{equation}\label{eqn:Z1}
    \fF_*Z_1\Fil_DW_m\Omega^q_U \otimes _{W_{m+1}\sO_X}W_{m+1}\sO_{X'}\cong \fF_*Z_1\Fil_DW_m\Omega^q_{U'}.
\end{equation}

To prove item (2), we use the exact sequence of $W_{m+1}\sO_X$-modules
(this is a part of the proof of \cite[Lem.~6.11]{KM-3})
\begin{equation}\label{eqn:RP-6-2}
  \fF^m_*\Fil_D\Omega^q_{U} \xrightarrow{V^m} \Fil_D W_{m+1} \Omega^q_{U}
  \xrightarrow{F}  \fF_*W_{m} \Omega^q_{U}.
\end{equation}
Taking the tensor product of this sequence with $W_{m+1}\sO_{X'}$ over
$W_{m+1}\sO_X$ and using Proposition~\ref{prop:RP-0} and ~\lemref{lem:RP-3}(7), we get
an analogous exact sequence of $W_{m+1}\sO_{X'}$-modules with $U$ replaced by
$U'$. This proves (2). 
Since the map $dV^{m-1} \colon \fF^{m}_* \Fil_D\Omega^r_U \to \fF_* \Fil_D W_m\Omega^{r+1}_U$ is $W_{m+1}\sO_X$-linear, by taking the tensor product with $W_{m+1}\sO_{X'}$ over
$W_{m+1}\sO_X$ and using Proposition~\ref{prop:RP-0} and ~\lemref{lem:RP-3}(7), we get
\begin{equation}\label{eqn:dV}
   \fF_* dV^{m-1}(\Fil_DW_m\Omega^q_U) \otimes _{W_{m+1}\sO_X}W_{m+1}\sO_{X'}\cong \fF_*dV^{m-1}(\Fil_DW_m\Omega^q_{U'}).
\end{equation}

Using \eqref{eqn:Z1} and \eqref{eqn:dV}, we deduce items
(3) and (4) from the corresponding assertions for $U$ established in
\thmref{thm:Global-version}, by tensoring with $W_{m+1}\sO_{X'}$ over $W_{m+1}\sO_X$
and applying Proposition~\ref{prop:RP-0} and ~\lemref{lem:RP-3}(7).
Finally, we deduce items (5) and (6) by using the first four items of the lemma,
together with \lemref{lem:RP-5}(5), and by mimicking
the proof of \cite[Prop.~10.1]{GK-Duality-1} mutatis mutandis.
\end{proof}

Using \lemref{lem:RP-6}(4), we get a filtration
\begin{equation}\label{eqn:RP-8}
  \{0\} \subseteq  B_1\Fil_D\Omega^q_{U'} \subseteq
  B_2\Fil_D\Omega^q_{U'} \subseteq
  \cdots \subseteq Z_2\Fil_D\Omega^q_{U'} \subseteq  Z_1\Fil_D\Omega^q_{U'} \subseteq
 \Fil_D\Omega^q_{U'}
\end{equation}
of $\Fil_D\Omega^q_{U'}$ analogous to the definition of $\Fil_D\Omega^q_U$
given after \lemref{lem:Log-fil-4}
so that $B_1\Fil_D\Omega^q_{U'} = d(\Fil_D\Omega^{q-1}_{U'})$ and
$Z_1\Fil_D\Omega^q_{U'} =
\Ker(\Fil_D\Omega^q_{U'} \xrightarrow{d} \Fil_D\Omega^{q+1}_{U'})$.
We let $B_0\Fil_D\Omega^q_{X'} =0$ and
$Z_0\Fil_D\Omega^q_{U'} := \Fil_D\Omega^q_{U'}$.
By iterating $\ov{F}$ and $\ov{C}$, we get isomorphisms
\begin{equation}\label{eqn:RP-8-0}
\ov{F}^r \colon \Fil_{D/{p^r}}\Omega^q_{U'} \xrightarrow{\cong} 
  \frac{Z_r\Fil_{D}\Omega^q_{U'}}{B_r\Fil_{D}\Omega^q_{U'}}; \ \
    \ov{C}^r \colon  \frac{Z_r\Fil_{D}\Omega^q_{U'}}{B_r\Fil_{D}\Omega^q_{U'}}
      \xrightarrow{\cong}  \Fil_{D/{p^r}}\Omega^q_{U'}
      \end{equation}
whose compositions are the identity maps.

\begin{lem}\label{lem:RP-9}
  Let $f \colon X' \to X$ be as in \lemref{lem:RP-3}. We then have the following.
  \begin{enumerate}
  \item
   $\Ker(\Fil_{D}\Omega^q_{U'} \xrightarrow{V^{m-1}}
    \frac{\Fil_{D}W_m\Omega^q_{U'}}{dV^{m-1}\Fil_{D}\Omega^{q-1}_{U'}}) =
     B_{m}\Fil_{D}\Omega^q_{U'}$.
  \item
    $\Ker(\Fil_{D}\Omega^q_{U'} \xrightarrow{V^{m-1}} \Fil_{D}W_m\Omega^q_{U'})
    = B_{m-1}\Fil_{D}\Omega^q_{U'}$.
  \end{enumerate}
\end{lem}
\begin{proof}
  We have exact sequences of $W_{m+1}\sO_X$-modules
  \begin{equation}\label{eqn:RP-9-0}
    0 \to {\fF}^m_*Z_m\Fil_{D}\Omega^q_{U} \to
    {\fF}^{m}_*Z_{m-1}\Fil_{D}\Omega^q_{U}
    \xrightarrow{{C}^{m-1}} \frac{\Re^{m}_*\Fil_{D/{p^{m-1}}}\Omega^q_{U}}
                { {\fF}^{m}_*Z_1\Fil_{D/{p^{m-1}}}\Omega^q_{U}} \to 0;
  \end{equation}
  \begin{equation}\label{eqn:RP-9-1}
    0 \to  {\fF}^m_*B_m\Fil_{D}\Omega^q_{U} \to {\fF}^m_*Z_m\Fil_{D}\Omega^q_{U}
    \xrightarrow{{C}^m} \Re^m_*\Fil_{D/{p^{m}}}\Omega^q_{U} \to 0.
  \end{equation}
  It follows from these exact sequences and induction on $m$ that
  the canonical maps ${\fF}^{m}_*Z_{m}\Fil_{D}\Omega^q_{U} \otimes_{W_{m+1}\sO_X}
  W_{m+1}\sO_{X'} \to {\fF}^{m}_*Z_{m}\Fil_{D}\Omega^q_{U'}$ (cf. \eqref{eqn:Z1}) and
${\fF}^{m}_*B_{m}\Fil_{D}\Omega^q_{U} \otimes_{W_{m+1}\sO_X}
  W_{m+1}\sO_{X'} \to {\fF}^{m}_*B_{m}\Fil_{D}\Omega^q_{U'}$ are isomorphisms
  of $W_{m+1}\sO_{X'}$-modules.

The lemma now follows by base changing the exact sequences of $W_m\sO_X$-modules
\[
  0 \to {\fF}^{m-1}_*B_{m-1}\Fil_{D}\Omega^q_{U} \to \fF^{m-1}_* \Fil_{D}\Omega^q_{U}
  \xrightarrow{V^{m-1}} \Fil_{D}W_m\Omega^q_{U};
  \]
  and of $W_{m+1}\sO_X$-modules
  \[
  0 \to {\fF}^{m}_*B_{m}\Fil_{D}\Omega^q_{U} \to \fF^{m}_* \Fil_{D}\Omega^q_{U}
  \xrightarrow{V^{m-1}}
  \frac{ {\fF}_*\Fil_{D}W_m\Omega^q_{U}}{dV^{m-1}\fF^{m}_*\Fil_{D}\Omega^{q-1}_{U}}
  \]
  (cf. items (7) and (11) of \thmref{thm:Global-version}) to $W_m\sO_{X'}$
  (resp. $W_{m+1}\sO_{X'}$) and  applying Proposition~\ref{prop:RP-0} and
  ~\lemref{lem:RP-1}.
\end{proof}

Most of the content of Lemmas~\ref{lem:RP-6} and ~\ref{lem:RP-9} can be illustrated
by the following commutative diagram, in which the  top row is exact and
$\ov{F} = {(\ov{C})}^{-1}$.
\begin{equation}\label{eqn:RP-10}
  \xymatrix@C1.5pc{
    0 \ar[r] & B_{m}\Fil_{D}\Omega^q_{U'} \ar[r] &
    \Fil_{D}\Omega^q_{U'} \ar[r]^-{V^{m}} &
    \Fil_{D}W_{m+1}\Omega^q_{U'} \ar[r]^-{F} \ar@{->>}[d]_-{R} &
    Z_1\Fil_{D}W_{m}\Omega^q_{U'} \ar[dr]^-{V} \ar@{->>}[d]^-{\can} \ar[dl]_-{C} \ar[r]
    & 0 \\
  & & & \Fil_{D/p}W_{m}\Omega^q_{U'} \ar@<1ex>[r]^-{\ov{F}} \ar@/_1.1cm/[rr]^-{\ov{p}} &
  \frac{Z_1\Fil_{D}W_m\Omega^q_{U'}}{dV^{m-1}\Fil_{D}\Omega^{q-1}_{U'}} \ar@<1ex>[l]^-{\ov{C}}
  \ar@{^{(}->}[r]^-{\ov{V}} & \Fil_{D}W_{m+1}\Omega^q_{U'} \\
  & & & & & \Fil_{D}\Omega^q_{U'} \ar[u]_-{V^m} \ar[lu]^-{\ \ V^{m-1}}.}
  \end{equation}

\section{Motivic complexes on the RP site}\label{sec:Mot-RP}
In this section, we continue to work under the setting of \S~\ref{sec:FDR-RP}.
We introduce the motivic complexes and establish two distinguished triangles
in the derived category of sheaves of abelian groups on $X_\rp$
associated to these complexes.

\begin{defn}\label{defn:RP-11} 
  Let $f \colon X' \to X$ be as in \lemref{lem:RP-3}.
  If $q \in \N_0$, $m \in \N$, and if $D \in \Div_E(X)$ is such that $D + E \ge 0$,
  we define
\[
W_m\sF^{q,\bullet}_{X'|D} =
\left[Z_1\Fil_{D} W_m\Omega^q_{U'} \xrightarrow{1-C} \Fil_{D} W_m\Omega^q_{U'}\right]
\]
as an object of $\sD^m_\et(X')$.
If $D \in \Div_E(X)$ is such that $D \le 0$, we define
\[
W_m\sG^{q,\bullet}_{X'|D} = \left[\Fil_{D} W_m\Omega^q_{U'} \xrightarrow{1-\ov{F}} 
\frac{\Fil_{D} W_m\Omega^q_{U'}}{dV^{m-1}(\Fil_{D}\Omega^{q-1}_{U'})}\right]
\]
as an object of $\sD^m_\et(X')$.
\end{defn}

It follows from \lemref{lem:RP-6} that the above complexes are well-defined.
They are the modulus versions of the classical {\'e}tale motivic complex
${\Z}/{p^m}(q)_X[q] \cong W_m\Omega^q_{X,\log}$ on $X_\et$.
Our goal is to establish two distinguished triangles associated to these
complexes.

\subsection{A distinguished triangle for $W_m\sG^{q,\bullet}_{X'|D}$}\label{sec:G-DT}
We begin with the following result for $W_m\sG^{q,\bullet}_{X'|D}$.

\begin{prop}\label{prop:RP-12}
    Let $f \colon X' \to X$ be as in \lemref{lem:RP-3} and assume that
    $D \le 0$. Then one has a distinguished triangle in $\sD^m_\et(X')$ of the form
  \begin{equation}\label{eqn:RP-12-0}
  W_1\sG^{q,\bullet}_{X'|{D/{p^{m-1}}}} \xrightarrow{\ov{p}^{m-1}}
     W_m\sG^{q,\bullet}_{X'|D} \xrightarrow{F}  W_{m-1}\sG^{q,\bullet}_{X'|D}
     \xrightarrow{+1}  W_1\sG^{q,\bullet}_{X'|{D/{p^{m-1}}}}[1].
  \end{equation}
  \end{prop}
  \begin{proof}
We let
$W_m\sE^{q,\bullet}_{X'|D} =
\left[Z_1\Fil_{D} W_m\Omega^q_{U'} \xrightarrow{1-C} \Fil_{D/p} W_m\Omega^q_{U'}\right]$.
It is easy to check using ~\eqref{eqn:RP-10} that the diagram
\begin{equation}\label{eqn:RP-11-0}
  \xymatrix@C1.5pc{
    0 \ar[r] & Z_1\Fil_{D} W_m\Omega^q_{U'} \ar[d]_-{1 - C} \ar[r] &
    \Fil_{D} W_m\Omega^q_{U'} \ar[r] \ar[d]^-{1-\ov{F}} &
    \frac{\Fil_{D} W_m\Omega^q_{U'}}{Z_1\Fil_{D} W_m\Omega^q_{U'}} \ar[d]^-{\id}
    \ar[r] & 0 \\
   0 \ar[r] & \Fil_{D/p} W_m\Omega^q_{U'} \ar[r]^-{-\ov{F}} & 
    \frac{\Fil_{D}W_m\Omega^q_{U'}}{dV^{m-1}(\Fil_{D}\Omega^{q-1}_{U'})} \ar[r] &
    \frac{\Fil_{D}W_m\Omega^q_{U'}}{Z_1\Fil_{D} W_m\Omega^q_{U'}} \ar[r] & 0}
\end{equation}
is commutative and its rows are exact. It follows that the horizontal arrows in the
left square of this diagram gives rise to a morphism of complexes of sheaves
\begin{equation}\label{eqn:RP-11-1}
  \Phi_m \colon W_m\sE^{q,\bullet}_{X'|D} \to W_m\sG^{q,\bullet}_{X'|D},
\end{equation}
on $X'_\et$ which is a quasi-isomorphism.

We next consider the diagram
\begin{equation}\label{eqn:RP-11-2}
  \xymatrix@C1.5pc{
    0 \ar[r] & \frac{\Fil_{D}\Omega^q_{U'}}{B_{m-1}\Fil_{D}\Omega^q_{U'}}
    \ar[r]^-{\ov{V}^{m-1}} \ar[d]_-{1-\ov{F}} &
      \Fil_{D} W_m\Omega^q_{U'} \ar[r]^-{F} \ar[d]^-{1-\ov{F}} &
      Z_1\Fil_{D} W_{m-1}\Omega^q_{U'} \ar[d]^-{1-C} \ar[r] & 0 \\
      0 \ar[r] & \frac{\Fil_{D}\Omega^q_{U'}}{B_{m}\Fil_{D}\Omega^q_{U'}}
      \ar[r]^-{\ov{V}^{m-1}} &
    \frac{\Fil_{D}W_m\Omega^q_{U'}}{dV^{m-1}(\Fil_{D}\Omega^{q-1}_{U'})} \ar[r]^-{-R} &
    \Fil_{D/p} W_{m-1}\Omega^q_{U'} \ar[r] & 0.}
    \end{equation}
It follows from Lemmas~\ref{lem:RP-5}, ~\ref{lem:RP-6} and ~\ref{lem:RP-9} that
this diagram is commutative and its rows are exact.

Letting $W_m\sL^{q,\bullet}_{X'|D} :=
\left(\frac{\Fil_{D}\Omega^q_{U'}}{B_{m-1}\Fil_{D}\Omega^q_{U'}} \xrightarrow{1-\ov{F}}
\frac{\Fil_{D}\Omega^q_{U'}}{B_{m}\Fil_{D}\Omega^q_{U'}}\right)$, we get a short
exact sequence of complexes of sheaves
\begin{equation}\label{eqn:RP-11-3}
0 \to W_m\sL^{q,\bullet}_{X'|D} \xrightarrow{\ov{V}^{m-1}} W_m\sG^{q,\bullet}_{X'|D} 
\xrightarrow{\Psi_m}  W_{m-1}\sE^{q,\bullet}_{X'|D} \to 0
\end{equation}
on $X'_\et$, where $\Psi_m$ is a morphism of 2-term complexes such that 
$\Psi^0_{m} = F$ and $\Psi^1_{m} = -R$.

We now consider the diagram
\begin{equation}\label{eqn:RP-11-4}
  \xymatrix@C1.5pc{
    0 \ar[r] & \Fil_{D/{p^{m-1}}} \Omega^q_{U'} \ar[r]^-{\ov{F}^{m-1}} \ar[d]_-{1- \ov{F}} &
    \frac{\Fil_{D}\Omega^q_{U'}}{B_{m-1}\Fil_{D}\Omega^q_{U'}} \ar[r] \ar[d]^-{1-\ov{F}} &
    \frac{\Fil_{D}\Omega^q_{U'}}{Z_{m-1}\Fil_{D}\Omega^q_{U'}} \ar[r]
    \ar[d]^-{1-\ov{F}} & 0 \\
    0 \ar[r] & \frac{\Fil_{D/{p^{m-1}}} \Omega^q_{U'}}{B_1\Fil_{D/{p^{m-1}}} \Omega^q_{U'}}
    \ar[r]^-{\ov{F}^{m-1}} & \frac{\Fil_{D}\Omega^q_{U'}}{B_{m}\Fil_{D}\Omega^q_{U'}}
    \ar[r] & \frac{\Fil_{D}\Omega^q_{U'}}{Z_{m-1}\Fil_{D}\Omega^q_{U'}} \ar[r] & 0.}
  \end{equation}
Using \lemref{lem:RP-6} and iterating $\ov{F}$ (cf. ~\eqref{eqn:RP-8-0}), it is
easy to check that this diagram is commutative and its rows are exact.

Since $\ov{F}^{m-1}(\Fil_{D}\Omega^q_{U'}) =
\frac{Z_{m-1}\Fil_{p^{m-1}D}\Omega^q_{U'} + Z_{m-1}\Fil_{D}\Omega^q_{U'}}{Z_{m-1}\Fil_{D}\Omega^q_{U'}} = 0$ in $\frac{\Fil_{D}\Omega^q_{U'}}{Z_{m-1}\Fil_{D}\Omega^q_{U'}}$, we get that
  $\ov{F} \colon \frac{\Fil_{D}\Omega^q_{U'}}{Z_{m-1}\Fil_{D}\Omega^q_{U'}} \to
  \frac{\Fil_{D}\Omega^q_{U'}}{Z_{m-1}\Fil_{D}\Omega^q_{U'}}$ is nilpotent.
  This implies that the right vertical arrow in ~\eqref{eqn:RP-11-4} is an
  isomorphism. It follows that the map of complexes
  $\ov{F}^{m-1} \colon W_1\sG^{q,\bullet}_{X'|{D/{p^{m-1}}}} \to W_m\sL^{q,\bullet}_{X'|D}$
  is a quasi-isomorphism.

 We finally look at the diagram of complexes
    \begin{equation}\label{eqn:RP-12-1}
\xymatrix@C1.5pc{ 
  0 \ar[r] & W_m\sL^{q,\bullet}_{X'|D} \ar[r]^-{\ov{V}^{m-1}} 
  & W_m\sG^{q,\bullet}_{X'|D} \ar[dr]_-{F} \ar[r]^-{\Psi_m} & W_{m-1}\sE^{q,\bullet}_{X'|D}
  \ar[r] \ar[d]^-{\Phi_{m-1}} & 0 \\
  & W_1\sG^{q,\bullet}_{X'|{D/{p^{m-1}}}} \ar[u]^-{\ov{F}^{m-1}} \ar[ru]_-{\ov{p}^{m-1}} & &
  W_{m-1}\sG^{q,\bullet}_{X'|D}. &}
\end{equation}

    One easily checks from the definitions of $\Phi_m$ and $\Psi_{m-1}$ that the
    right triangle in ~\eqref{eqn:RP-12-1} is commutative.
    Since $\ov{V} \circ \ov{F} = \ov{F} \circ \ov{V}$ as
    one checks using the relation $VF = FV$ and the
    diagram ~\eqref{eqn:RP-10}, it follows that $\ov{V}^{m-1}\circ \ov{F}^{m-1} =
    (\ov{V} \circ \ov{F})^{m-1} =\ov{p}^{m-1}$. Since the left and right vertical arrows
    in ~\eqref{eqn:RP-12-1} are quasi-isomorphisms of complexes, we get the
    desired  distinguished triangle in $\sD^m_\et(X')$. This completes the proof.
\end{proof}

We shall now derive some consequences of \propref{prop:RP-12}.
For $f \colon X' \to X$ as in \lemref{lem:RP-3} with $D \ge 0$, let
$D' = f^*(D)$ and let $j' \colon U' \inj X'$ be the inclusion.  We let
\begin{equation}\label{eqn:RP-12-2}
W_m\Omega^q_{X'|D, \log} :=  \Ker\left(\Fil_{-D} W_m\Omega^q_{U'} \xrightarrow{1- \ov F} 
\frac{\Fil_{-D} W_m\Omega^q_{U'}}{dV^{m-1}(\Fil_{-D}\Omega^{q-1}_{U'})}\right).
\end{equation}

The notation $W_m\Omega^q_{X|D, \log}$ was introduced in \cite[Defn.~1.1.1]{JSZ}. 
Our reason for using this notation in ~\eqref{eqn:RP-12-2}
is the following.

\begin{lem}\label{lem:RP-21}
   \begin{enumerate}
  \item
  If $E = \emptyset$, then
  $W_m\Omega^q_{X'|D, \log}$ coincides with $W_m\Omega^q_{X',\log}$, where the
  latter is the image of the map of {\'e}tale sheaves
$\dlog \colon \sK^M_{q,X'} \to W_m\Omega^q_{X'}$.
\item
  If $X'$ is Noetherian, we have the following.
  \begin{listabc}
  \item
 If $D \ge E$, then $W_m\Omega^q_{X'|D, \log}$ coincides with
  the image of the map of {\'e}tale sheaves
$\dlog \colon \Ker(\sK^M_{1,X'} \surj \sK^M_{1,D'}) \otimes j'_* (\sK^M_{q-1,U'})
  \to W_m\Omega^q_{X'}$.
\item
  If $E \neq \emptyset$ but $D = 0$, then
  $W_m\Omega^q_{X'|D, \log}$ coincides with $j'_*W_m\Omega^q_{U',\log}$.
  \end{listabc}
  \end{enumerate}
 \end{lem}
\begin{proof}
 Item (1) follows from \cite[Cor.~4.2]{Morrow-ENS} and item (2) follows by
  combining \lemref{lem:RP-20}, \cite[Prop.~1.1.3]{JSZ} (which implies that
  the image of $\dlog$ lies in $W_m\Omega^q_{X'}$),
  and items (1) and (5) of \cite[Lem.~10.2]{KM-3}.
\end{proof}

\begin{lem}\label{lem:RP-13}
  The inclusion $\vartheta^q_m \colon W_m\Omega^q_{X'|D, \log}
  \inj W_m\sG^{q,\bullet}_{X'|-D}$ of complexes of {\'e}tale sheaves is
    a quasi-isomorphism. 
  \end{lem}
  \begin{proof}
    The lemma is equivalent to the assertion that
    $\sH^1(W_m\sG^{q,\bullet}_{X'|-D}) = 0$. To prove this latter assertion,
    we use \propref{prop:RP-12} (especially, the diagram ~\eqref{eqn:RP-12-1}),
    which gives us an exact sequence of {\'e}tale sheaves
\[
    \sH^1(W_1\sG^{q,\bullet}_{X'|{-D/{p^{m-1}}}}) \to
    \sH^1(W_{m}\sG^{q,\bullet}_{X'|-D}) \to \sH^1(W_{m-1}\sG^{q,\bullet}_{X'|-D}).
    \]
    Using this exact sequence and induction on $m \ge 1$, the problem is reduced
    to showing the desired
    vanishing when $m = 1$. Being a local question in the {\'e}tale
    topology, the latter case follows from \lemref{lem:RP-13-0}.
\end{proof}

  The following result is a variant of \cite[Claim~1.2.2]{JSZ}.
  
\begin{lem}\label{lem:RP-13-0}
    Let $R$ be a regular local $\F_p$-algebra with maximal ideal
    $(t_1, \ldots , t_d)$.
    We let $\pi = t_1\cdots t_r$ and $t = t^{n_1}_1\cdots t^{n_r}_r$ with
    $n_i \ge 0$ for each $1 \le i \le r \le d$. 
 Let $A$ be a strict henselian local $\F_p$-algebra with maximal ideal $\fm_A$ and
 residue field $k$. Let $R \to A$ be a relatively perfect local ring homomorphism
 and let $A' = A[\pi^{-1}]$. Then the map
 $t\Omega^q_A(\log (\pi)) \xrightarrow{1-\ov{F}}
 \frac{t\Omega^q_A(\log (\pi))}{d(t\Omega^{q-1}_A(\log (\pi)))}$ is surjective.
\end{lem}
  \begin{proof}
 By \lemref{lem:RP-3}(4) and the definition of $\Omega^q_R(\log(\pi))$
    (cf. \cite[\S~2.1]{KM-3}), a typical generating element of
    $t\Omega^q_A(\log (\pi))$ is of the form
    $w = at\dlog(x_1)\cdots \dlog(x_q)$, where $a \in A$ and $x_i \in A'^{\times}$.
    It suffices to show that any such element $w$ has a pre-image under
    $1-\ov{F}$ modulo $d(t\Omega^{q-1}_A(\log (\pi)))$. To that end, suppose
    $\alpha = at$ and there exists an element $b \in tA$ such that
    $b - b^p = \alpha$. We then get
    $(1-\ov{F})(b\dlog(x_1)\cdots \dlog(x_q)) =
    (1-\ov{F})(b)\dlog(x_1)\cdots \dlog(x_q) = w$.
    We are therefore reduced to showing that there exists an element $b \in tA$
    such that $b - b^p = \alpha$. 

Since $A$ is strict henselian, there exists $a' \in A$ such that
  $a'- a'^p = \alpha$. If $t =1$, we are already done.
    We assume now that $t \in \fm_A$. Letting $\ov{b}$ denote the image in $k$ of
    any element $b \in A$, we get that $\ov{a'} \in \F_p \subset k$. We can
    therefore choose $b' \in \F_p \subset A$ such that $a'-b' \in \fm_A$.
    Letting $b = a'-b'$, we get $b - b^p = a'-b' -a'^p + b'^p = a' - a'^p = \alpha$
    (since $b' \in \F_p$). Furthermore, $b-b^p = \alpha = at$ implies that
    $b(1 - b^{p-1}) = at$. Since $1 - b^{p-1} \in A^\times$, we get $b \in tA$.
    We thus get $b \in tA$ such that $b - b'^p = \alpha$. This completes the
    proof.
\end{proof}

  The following application of \propref{prop:RP-12} provides an extension of
  \cite[Thm.~1.1.6]{JSZ} to relatively perfect schemes over $F$-finite regular
  schemes.

\begin{thm}\label{thm:RP-13-1}
   Let $f \colon X' \to X$ be as in \lemref{lem:RP-3} with $D \ge 0$, and let
   $1 \le i \le m-1$. Then one has an exact sequence of {\'e}tale sheaves on $X'$:
   \begin{equation}\label{eqn:RP-13-1-0}
     0 \to W_i\Omega^q_{X'|{\lceil{D/{p^{m-i}}\rceil}}, \log} \xrightarrow{\ov{p}^{m-i}}
     W_m\Omega^q_{X'|D, \log}
   \xrightarrow{R^i} W_{m-i}\Omega^q_{X'|D, \log} \to 0.
   \end{equation}
\end{thm}
\begin{proof}
  Using induction on $m \ge 1$, one easily reduces the proof of the theorem to the
  case $i =1$. 
By \lemref{lem:RP-13}, we have
\begin{equation}
    \sH^i(W_m\sG^{q,\bullet}_{X'|-D})
\cong
\begin{cases}
W_m\Omega^q_{X'|D, \log} & \text{if } i=0,\\
0 & \text{if } i\neq 0.
\end{cases}
\end{equation}
Therefore, using \propref{prop:RP-12} and the long exact sequence of the cohomology sheaves associated with the distinguished triangle \eqref{eqn:RP-12-0}, it remains to show that $F = R$ when restricted to the subsheaf $W_m\Omega^q_{X'|D, \log} \subset \Fil_{-D}W_m\Omega^q_{U'}$. Since $\Fil_{-D}W_m\Omega^q_{U'} \subset \Fil_{0}W_m\Omega^q_{U'}$, it suffices to prove the claim for $W_m\Omega^q_{X'|0, \log}$.

We now look at the commutative diagram

\begin{equation}
    \xymatrix@C1.7pc{
   \Fil_{0}W_{m+1}\Omega^q_{U'} \ar[r]^{R-F}\ar@{->>}[d]_{R}& \Fil_{0}W_m\Omega^q_{U'}\ar@{->>}[d]\\
    \Fil_{0}W_m\Omega^q_{U'} \ar[r]^{1- \ov F} & \frac{\Fil_{0}W_m\Omega^q_{U'}}{dV^{m-1}(\Fil_{0}\Omega^{q-1}_{U'}),}
    }
\end{equation}
where the left vertical arrow is surjective by \lemref{lem:RP-5}.
Since $(R-F) \circ dV^m(\Fil_0\Omega_{U'}^{q-1})=-dV^{m-1}(\Fil_0\Omega_{U'}^{q-1})$, the restriction of $R-F$ to the kernels of the vertical arrows is surjective. Using
a diagram chase, we conclude that the induced map $R:\Ker (R-F) \to W_m\Omega^q_{X'|0, \log}$ is surjective. Combining this with the fact that $RF=FR$, we obtain the claim.  
\end{proof}

\subsection{A distinguished triangle for $W_{m}\sF^{q,\bullet}_{X'|D}$}\label{sec:F-DT}
Let $f \colon X' \to X$ be as in \lemref{lem:RP-3} with $D \in \Div_E(X)$
  such that $D + E \ge 0$. Let $q \in \N_0$ and $m \in \N$.
  Our next goal is to extend \thmref{thm:Global-version}(9) to
  $X_{\rp}$. We make the following definitions before we state the result.

\begin{defn}\label{def:two-term*}
We let
  \[
  Z_1R^{q,m}_{X'|D} :=
\Ker\left(Z_1\Fil_{D}\Omega^q_{U'} \oplus \Fil_{D}\Omega^{q-1}_{U'}
\xrightarrow{(V^{m},\ dV^{m})} Z_1\Fil_{D}W_{m+1}\Omega^q_{U'}\right) \text{ and}
\]
\[
R^{q,m}_{X'|D} := \Ker\left(\Fil_{D}\Omega^q_{U'} \oplus \Fil_{D}\Omega^{q-1}_{U'}
\xrightarrow{(V^{m}, \ dV^{m})} \Fil_{D}W_{m+1}\Omega^q_{U'}\right).
\]
\end{defn}

\begin{lem}\label{lem:two-term*-0}
  The canonical inclusion $Z_1R^{q,m}_{X'|D} \inj R^{q,m}_{X'|D}$ is an isomorphism.
\end{lem}
\begin{proof}
  It suffices to show that
  $\wt{R}^{q,m}_{X'|D} := R^{q,m}_{X'|D} \bigcap \Fil_{D}\Omega^q_{U'}
  \subset Z_1\Fil_{D}\Omega^q_{U'}$. But this is clear since
  $\wt{R}^{q,m}_{X'|D} \subset \Ker\left(\Fil_{D}\Omega^q_{U'} \xrightarrow{V^m}
  \frac{\Fil_{D}W_{m+1}\Omega^q_{U'}}{dV^{m}(\Fil_{D}\Omega^{q-1}_{U'})}\right)$,
  and the latter is contained in $B_{m}\Fil_{D}\Omega^{q}_{U'}$
  (hence in $Z_1\Fil_{D}\Omega^q_{U'}$) by \lemref{lem:RP-9}.
\end{proof}

\begin{prop}\label{prop:RP-14}
  There exists a distinguished triangle in $\sD^m_{\et}(X')$ of the form
  \begin{equation}\label{eqn:RP-14-0}
    W_1\sF^{q,\bullet}_{X'|D} \xrightarrow{V^{m-1}} W_m\sF^{q,\bullet}_{X'|D}
    \xrightarrow{R} W_{m-1}\sF^{q,\bullet}_{X'|D/p} \xrightarrow{+1}
    W_1\sF^{q,\bullet}_{X'|D}[1].
  \end{equation}
  \end{prop}
\begin{proof}
We begin by looking at the diagram
\begin{equation}\label{eqn:RP-14-1}
    \xymatrix@C1.5pc{
{R}^{q,m-1}_{X'|D} \ar[d]_-{\phi} \ar@{^{(}->}[r]^-{\alpha} &  
 Z_1\Fil_{D}\Omega^q_{U'} \oplus \Fil_{D}\Omega^{q-1}_{U'} \ar[r]^-{\beta} 
 \ar[d]^-{\phi} & Z_1\Fil_{D}W_m\Omega^q_{U'} \ar[d]^-{(1-C)} \ar@{->>}[r]^-{R} &
 Z_1\Fil_{D/p}W_{m-1}\Omega^q_{U'} \ar[d]^-{(1-C)} \\
 {R}^{q,m-1}_{X'|D} \ar@{^{(}->}[r]^-{\alpha} &  
 \Fil_{D}\Omega^q_{U'} \oplus \Fil_{D}\Omega^{q-1}_{U'} \ar[r]^-{\beta} 
 & \Fil_{D}W_m\Omega^q_{U'} \ar@{->>}[r]^-{R} & \Fil_{D/p}W_{m-1}\Omega^q_{U'},}
\end{equation}    
where we let $\alpha$ denote the canonical inclusion,
$\wt{C}(a, b) = C(a), \ \phi = 1-\wt{C}$ and $\beta = (V^{m-1},\ dV^{m-1})$.
The term ${R}^{q,m-1}_{X'|D}$ in the top row has been identified with
$Z_1{R}^{q,m-1}_{X'|D}$ using \lemref{lem:two-term*-0}.
It is clear that this diagram is commutative.

The bottom row in ~\eqref{eqn:RP-14-1} is exact by \lemref{lem:RP-5}(5). The map
$R$ in the top row is surjective by its surjectivity in the bottom row and
\lemref{lem:RP-6}(1). To show the exactness of the top row, we only have to
check its exactness at its third term, which we can do {\'e}tale
locally. Since the bottom row is exact, it suffices to show that
if $(a,b) \in \Fil_{D}\Omega^q_{U'} \oplus \Fil_{D}\Omega^{q-1}_{U'}$
such that $V^{m-1}(a) + dV^{m-1}(b) \in Z_1\Fil_{D}W_m\Omega^q_{U'}$,
then $a \in Z_1\Fil_{D}\Omega^q_{U'}$.
Now, by \lemref{lem:RP-6}(1), our assumption implies that
$F^{m-1}d(V^{m-1}(a) + dV^{m-1}(b)) = 0$. In turn, this implies that
$d(a) = F^{m-1}dV^{m-1}(a) =0$. That is, $a \in Z_1\Fil_{D}\Omega^q_{U'}$.

We let $\sF^{q, \bullet}_{X'|D}$ (resp. $\wt{\sF}^{q,\bullet}_{X'|D}$)
    denote the 2-term complex consisting of the second (resp. first)
    column of ~\eqref{eqn:RP-14-1} from the left. To prove the proposition,
    it remains to show that the composite map
   \begin{equation}\label{eqn:RP-14-2}
     W_1\sF^{q,\bullet}_{X'|D} \inj \sF^{q, \bullet}_{X'|D} \surj
     \frac{\sF^{q, \bullet}_{X'|D}}{\wt{\sF}^{q,\bullet}_{X'|D}}
     \end{equation}
is a quasi-isomorphism.

To that end, we note that the injective arrow
    in ~\eqref{eqn:RP-14-2} is a quasi-isomorphism
    by the definition of $\phi$. It therefore suffices to show that the surjective
    arrow in  ~\eqref{eqn:RP-14-2} is also a quasi-isomorphism.
    Equivalently, we need to
    show that the left vertical arrow in ~\eqref{eqn:RP-14-1} is bijective.     
   To this end, observe that the composite inclusion
    ${R}^{q,m-1}_{X'|D} \inj \Fil_{D}\Omega^q_{U'} \oplus \Fil_{D}\Omega^{q-1}_{U'}
    \inj \Omega^q_{U'} \oplus \Omega^{q-1}_{U'}$ factors through the inclusion
    ${R}^{q,m-1}_{X'|D} \inj B_m\Omega^q_{U'} \oplus \Omega^{q-1}_{U'}$ by
    \lemref{lem:RP-9}(2), as we showed earlier. 
 Furthermore, $C^{m} = 0$ on $B_m\Omega^q_{U'}$.
Since $\wt{C} = (C,0)$ vanishes on $\Omega^{q-1}_{U'}$, it follows that
${(\wt{C})}^{m} = 0$ on ${R}^{q,m-1}_{X'|D}$. This implies the required
bijectivity and completes the proof of the proposition.
\end{proof}

\begin{cor}\label{cor:RP-16}
  For $1 \le i \le m-1$, there exists a distinguished triangle in $\sD^m_\et(X')$ of
  the form
  \[
  W_{i}\sF^{q, \bullet}_{X'|D} \xrightarrow{V^{m-i}}  W_m\sF^{q, \bullet}_{X'|D}
  \xrightarrow{R^{i}}  W_{m-i}\sF^{q, \bullet}_{X'|{D/{p^{i}}}} \xrightarrow{+1}
    W_{i}\sF^{q, \bullet}_{X'|D}[1].
    \]
  \end{cor}
\begin{proof}
  This is easily deduced from \propref{prop:RP-14} using induction on $m$ and
  a diagram chase.
  \end{proof}

\subsection{Pairing between $\sF$ and $\sG$ complexes}\label{sec:Pair*}
Let $f \colon X' \to X$ be as in \lemref{lem:RP-3} with $D \ge 0$. 
Let $q \in \N_0, \ m \in \N, \ q' = d-q$ and $D' = D+E$, where recall that
$d$ is the $p$-rank of $X$.
We shall now define a pairing between $W_m\sF^{q, \bullet}_{X'|D}$ and
$W_m\sG^{q, \bullet}_{X'|-D'}$.
We let

\begin{equation}\label{eqn:RP-15}
W_m\sF^{q,\bullet}_{X'} =
\left[Z_1W_m\Omega^q_{X'} \xrightarrow{1-C} W_m\Omega^q_{X'}\right]; \
W_m\sG^{q,\bullet}_{X'} = \left[W_m\Omega^q_{X'} \xrightarrow{1-\ov{F}}
  \frac{W_m\Omega^q_{X'}}{dV^{m-1}(\Omega^{q-1}_{X'})}\right];
\end{equation}
\[
W_m\sH^{d,\bullet}_{X'} =
\left[W_m\Omega^d_{X'} \xrightarrow{1-C} W_m\Omega^d_{X'}\right]
\]
as objects of $\sD^m_\et(X')$.

By \lemref{lem:RP-3}(8) and the equality $Z_1\Fil_{-E}W_m\Omega^d_{U'} =
\Fil_{-E}W_m\Omega^d_{U'}$, we see that $W_m\sH^{d,\bullet}_{X'}$ coincides with the
complex
$\left[Z_1\Fil_{-E}W_m\Omega^d_{U'} \xrightarrow{1-C} \Fil_{-E}W_m\Omega^d_{U'}\right]$.
By Lemmas~\ref{lem:RP-21} and ~\ref{lem:RP-13}, and ~\eqref{eqn:RP-11-0},
we also have that there is a canonical inclusion
$\vartheta^d_m \colon W_m\Omega^d_{X', \log} \inj  W_m\sH^{d,\bullet}_{X'}$
which is a quasi-isomorphism. 
We thus get

\begin{equation}\label{eqn:RP-15-0}
W_m\Omega^d_{X', \log}
  \xrightarrow{\simeq} W_m\sH^{d,\bullet}_{X'}
  \xrightarrow{\cong}
  \left[Z_1\Fil_{-E}W_m\Omega^d_{U'} \xrightarrow{1-C} \Fil_{-E}W_m\Omega^d_{U'}\right].
\end{equation}

We now define the pairing between the complexes $W_m\sF^{q,\bullet}_{X'|D}$ and
$W_m\sG^{q',\bullet}_{X'|D'}$ as follows.
We let
\begin{equation}\label{eqn:RP-15-1}
  \< \ , \ \>^{0,0}_0 \colon Z_1\Fil_{D}W_m\Omega^q_{U'} \times
  \Fil_{-D'}W_m\Omega^{q'}_{U'}
  \to \Fil_{-E}W_m\Omega^d_{U'} = W_m\Omega^d_{X'}
\end{equation}
and
\begin{equation}\label{eqn:RP-15-2}
  \< \ , \ \>^{1,0}_1 \colon \Fil_{D}W_m\Omega^q_{U'} \times
  \Fil_{-D'}W_m\Omega^{q'}_{U'}
  \to \Fil_{-E}W_m\Omega^d_{U'} = W_m\Omega^d_{X'}
\end{equation}
be given by $\<w_1,w_2\>^{0,0}_0 := w_1 \wedge w_2$ and $\<w_1,w_2\>^{1,0}_1 :=
w_1 \wedge w_2$.

We let
\begin{equation}\label{eqn:RP-15-3}
  \< \ , \ \>^{0,1}_1 \colon Z_1\Fil_{D}W_m\Omega^q_{U'} \times
  \frac{\Fil_{-D'}W_m\Omega^{q'}_{U'}}{dV^{m-1}(\Fil_{-D'}\Omega^{q'-1}_{U'})}
  \to \Fil_{-E}W_m\Omega^d_{U'} = W_m\Omega^d_{X'}
\end{equation}
be defined by $\<w_1, \ov{w}_2\>^{0,1}_1 := -C(w_1 \wedge w_2)$.

The first two pairings are well-defined by \lemref{lem:RP-3}(3). To
show that the third pairing is also defined, we only need to prove that
$C(w_1 \wedge w_2) = 0$ if $w_1 \in Z_1\Fil_{D}W_m\Omega^q_{U'}$ and
$w_2 \in dV^{m-1}(\Fil_{-D'}\Omega^{q'-1}_{U'})$.
To prove it, we can write $w_1 = F(w'_1)$ (cf. \lemref{lem:RP-6}(1))
and $w_2 = dV^{m-1}(w'_2)$, where
$w'_1 \in \Fil_{D}W_{m+1}\Omega^q_{U'}$ and $w'_2 \in \Fil_{-D'}\Omega^{q'-1}_{U'}$.
Then
\[
C(w_1 \wedge w_2) = C(F(w'_1) \wedge dV^{m-1}(w'_2)) = CF(w'_1 \wedge dV^m(w'_2)) =
R(w'_1 \wedge dV^m(w'_2))
\]
\[
\hspace*{4cm} = R(w'_1)\wedge dRV^m(w'_2) = 0,
\]
where the last equality holds because $RV^m =0$ by \lemref{lem:RP-5}(5).
For the definition of a pairing of two complexes, we refer to
\cite[\S~1, p.~175]{Milne-ENS}.

\begin{lem}\label{lem:RP-15-4}
The above pairings of {\'e}tale sheaves induce a pairing of
  the complexes
  \begin{equation}\label{eqn:RP-15-4-0}
    W_m\sF^{q, \bullet}_{X'|D} \times W_m\sG^{q', \bullet}_{X'|-D'}
    \xrightarrow{{\< \ , \ \>}}
    W_m\sH^{d,\bullet}_{X'}
  \end{equation}
  such that one has a commutative diagram of pairings
  \begin{equation}\label{eqn:RP-15-4-1}
  \xymatrix@C1.5pc{  
    W_{m}\sF^{q, \bullet}_{X'|D} \times W_{m}\sG^{q', \bullet}_{X'|-D'}  \ar[r]^-{\< \ , \ \>}
    \ar@<5ex>[d] &
  W_m\sH^{d,\bullet}_{X'} \ar@{=}[d] \\
  W_{m} \sF^{q, \bullet}_{X'} \times W_{m}\sG^{q', \bullet}_{X'}  \ar[r]^-{\< \ , \ \>}
  \ar@<7ex>[u] & W_{m}\sH^{d, \bullet}_{X'},}
  \end{equation}
  in which the vertical arrows are the canonical maps.
\end{lem}
\begin{proof}
To show that ~\eqref{eqn:RP-15-4-0} is well-defined, it suffices to verify
  (cf. \cite[\S~1, p.~175]{Milne-ENS}) that
  \[
  (1-C)(w_1 \wedge w_2) = - C(w_1 \wedge (1-\ov{F})(w_2)) + (1-C)(w_1) \wedge w_2
  \]
  for all $w_1 \in Z_1\Fil_{D}W_m\Omega^q_{U'}$ and
  $w_2 \in \Fil_{D'}W_m\Omega^{q'}_{U'}$. 
 By \cite[Lem.~1.1(c)]{Milne-ENS} and the identity $\ov{C} \circ \ov{F} = \id$
  (cf. ~\eqref{eqn:RP-10}), the above identity follows readily once we
  establish that $C(w_1 \wedge w_2) = C(w_1) \wedge C(w_2)$ if
  $w_1, w_2 \in Z_1W_m\Omega^q_{U'}$. To prove this, we use \lemref{lem:RP-6}(1)
  to write $w_i = F(w'_i)$. Then
  $C(w_1 \wedge w_2) = R(w'_1 \wedge w'_2)$ while $C(w_1) \wedge C(w_2) =
  R(w'_1) \wedge R(w'_2) = R(w'_1 \wedge w'_2)$.
  This establishes the desired identity.
The commutativity of ~\eqref{eqn:RP-15-4-1} is clear.
\end{proof}

\section{Kato's duality for modulus pairs}\label{sec:KD*}
In this section, we prove \thmref{thm:Main-1}.  We continue to work under the
setting and assumptions of \S~\ref{sec:FDR-RP}.

\subsection{Duality pairing on RP site}\label{sec:DRP-site}
For $m \in \N$ and a quasi-coherent sheaf $\sF$ of $W_m\sO_X$-modules on $X_\et$, we
denote by $\sF^{\rp}$ the pull-back of $\sF$ along the morphism of sites
$\tau \colon X_{\rp} \to X_\et$. More generally, for any sheaf of $W_m\sO_X$-modules
$\sF$, we denote by $\sF^{\rp}$ the sheafification on $X_{\rp}$ of the presheaf
$(f \colon X' \to X) \mapsto H^0_\et(X', \sF \otimes_{W_m\sO_X} W_m\sO_{X'})$.
It follows from Lemmas~\ref{lem:RP-1} and ~\ref{lem:RP-3}, 
together with the proof of \lemref{lem:RP-5}, that
the presheaf $X' \mapsto \Fil_DW_m\Omega^q_{U'}$ is a sheaf on $X_{\rp}$ and is
canonically isomorphic to $\Fil_DW_m\Omega^{q, \rp}_{U}$. The same holds for
the presheaves $X' \mapsto Z_1\Fil_DW_m\Omega^q_{U'}$ and $X' \mapsto
\frac{\Fil_DW_m\Omega^q_{U'}}{dV^{m-1}(\Fil_D\Omega^{q-1}_{U'})}$.

If $D \in \Div_E(X)$ is such that $D+E \ge 0$, we let
  \begin{equation}\label{eqn::RP-16-1}
    W_m\sF^{q,\rp}_{X|D} =
\left[Z_1\Fil_{D} W_m\Omega^{q, \rp}_{U} \xrightarrow{1-C}
  \Fil_{D} W_m\Omega^{q, \rp}_{U}\right];
\end{equation}
and if $D \le 0$, we let
\[
W_m\sG^{q,\rp}_{X|D} = \left[\Fil_{D} W_m\Omega^{q,\rp}_{U} \xrightarrow{1-\ov{F}} 
\frac{\Fil_{D} W_m\Omega^{q,\rp}_{U}}{dV^{m-1}(\Fil_{D}\Omega^{q-1, \rp}_{U})}\right].
\]

It is clear that the restriction of $W_m\sF^{q,\rp}_{X|D}$ (resp. $W_m\sG^{q,\rp}_{X|D}$)
to $X'_\et$, for any $X' \in X_{\rp}$, coincides with
$W_m\sF^{q}_{X'|D}$ (resp. $W_m\sG^{q}_{X'|D}$) as defined in Definition~\ref{defn:RP-11}.
It also follows that
the presheaf $X' \mapsto W_m\Omega^{q}_{X'|D, \log}$ on $X_{\rp}$
is a sheaf, which we denote by $W_m\Omega^{q, \rp}_{X|D, \log}$.
We write $W_m\Omega^{q, \rp}_{X|D, \log}$ as $W_m\Omega^{q, \rp}_{X, \log}$ when
$E = \emptyset$.

From \thmref{thm:RP-13-1} and \corref{cor:RP-16}, we deduce the following.
\begin{thm}\label{thm:RP-16-3}
For $1 \le i \le m-1$ and $D \in \Div_E(X)$ such that $D+E \ge 0$, there exists a distinguished
  triangle in $D_{m}(X_\rp)$ of the form
  \[
  W_{i}\sF^{q, \rp}_{X|D} \xrightarrow{V^{m-i}}  W_m\sF^{q, \rp}_{X|D}
  \xrightarrow{R^{i}}  W_{m-i}\sF^{q, \rp}_{X|{D/{p^{i}}}} \xrightarrow{+1}
    W_{i}\sF^{q, \rp}_{X|D}[1].
    \]
    If $D \ge 0$, there exists an exact sequence of sheaves on $X_{\rp}$ of the
    form
    \[
0 \to W_i\Omega^{q,\rp}_{X|{\lceil{D/{p^{m-i}}\rceil}}, \log} \xrightarrow{\ov{p}^{m-i}}
     W_m\Omega^{q, \rp}_{X|D, \log}
     \xrightarrow{R^i} W_{m-i}\Omega^{q, \rp}_{X|D, \log} \to 0.
     \]
\end{thm}

Let $D_m(X_{\rp})$ denote the derived category of sheaves of
$\Lambda_m$-modules on $X_{\rp}$ and let $D^c_m(X_{\rp})$ denote the
full triangulated subcategory of $D_m(X_{\rp})$
generated by those objects of $D_m(X_{\rp})$ 
which are represented by coherent $\sO^{\rp}_X$ modules regarded as complexes
of $\Lambda_m$-modules supported in degree zero. 
For any object $E \in D_m(X_{\rp})$, we let
$\D_{X,m}(E) := {\bf R}\sHom_{D_m(X_{\rp})}(E, W_m\Omega^{d, \rp}_{X, \log})
\xrightarrow{\cong} {\bf R}\sHom_{D_m(X_{\rp})}(E, W_m\sH^{d, \rp}_X)$,
where the isomorphism is induced by $\vartheta^d_m$ (cf. \S~\ref{sec:Pair*}).

\begin{lem}\label{lem:RP-16-2}
  Each of the objects $W_m\sF^{q,\rp}_{X|D}, \ W_m\sG^{q,\rp}_{X|D}$ and $W_m\sH^{d, \rp}_X$
   in $D_m(X_{\rp})$ lies in $D^c_m(X_{\rp})$. In particular,
    $W_m\Omega^{q, \rp}_{X|D, \log}, W_m\Omega^{d, \rp}_{X, \log} \in D^c_m(X_{\rp})$.
  \end{lem}
  \begin{proof}
    The first part of the lemma for $m =1$ follows from \lemref{lem:Can-filt-5}.
    For $m \ge 1$, this follows by induction on $m$ together with
    \thmref{thm:RP-16-3}.
    The second part follows from \lemref{lem:RP-13}.
    \end{proof}

One checks directly from its construction that the pairing ~\eqref{eqn:RP-15-4-0} 
is compatible with respect to morphisms of schemes in $X_\rp$, and it
is also compatible with $D \in \Div^{\eff}_E(X)$ and $m \in \N$. In particular,
we get for $q \in \N_0$ and $D \ge 0$, a pairing
\begin{equation}\label{eqn:RP-15-4-2}
 W_m\sF^{q, \rp}_{X|D} \times W_m\sG^{q', \rp}_{X|-D'}
    \xrightarrow{{\< \ , \ \>}}
    W_m\sH^{d,\rp}_{X}
  \end{equation}
on $X_\rp$, where $q' = d-q$ and $D' = D+E$.
This induces the duality maps (cf. \cite[\S~1, p.~175]{Milne-ENS})
\begin{equation}\label{eqn:RP-15-4-3}
  \Phi^{q,m}_{X|D} \colon W_m\sF^{q,  \rp}_{X|D} \to
  \sHom^\bullet_{\Lambda_m}(W_{m}\sG^{q', \rp}_{X|-D'}, W_m\sH^{d,\rp}_{X})
  \to \D_{X,m}(W_{m}\sG^{q',\rp}_{X|-D'});
\end{equation}
\[
\Psi^{q,m}_{X|D} \colon W_m\sG^{q', \rp}_{X|-D'} \to
\sHom^\bullet_{\Lambda_m}(W_{m}\sF^{q, \rp}_{X|D}, W_m\sH^{d,\rp}_{X}) \to
\D_{X,m}(W_{m}\sF^{q \rp}_{X|D}).
\] 
Composing $\Psi^{q,m}_{X|D}$ with the inclusion
$W_m\Omega^{q', \rp}_{X|D',\log} \xrightarrow{\vartheta^{q'}_m} W_m\sG^{q', \rp}_{X|-D'}$,
and $\Phi^{q,m}_{X|D}$ with the map
${\bf R}\sHom^\bullet_{D_m(X_\rp)}(W_{m}\sG^{q', \rp}_{X|-D'}, W_m\sH^{d,\rp}_{X})
\xrightarrow{(\vartheta^{q'}_m)^\star}
{\bf R}\sHom^\bullet_{D_m(X_\rp)}(W_m\Omega^{q', \rp}_{X|D',\log}, W_m\sH^{d,\rp}_{X})$,
we obtain the modulus extension of Kato's duality maps
(cf. \cite[\S~4]{Kato-duality-1}):

\begin{equation}\label{eqn:RP-15-4-5}
  \Phi^{q,m}_{X|D} \colon W_m\sF^{q, \rp}_{X|D} \to
  \D_{X,m}(W_m\Omega^{q', \rp}_{X|D',\log}); \ \
  \Psi^{q,m}_{X|D} \colon W_m\Omega^{q', \rp}_{X|D', \log} \to \D_{X,m}(W_m\sF^{q, \rp}_{X|D}).
\end{equation}

Kato proved in Theorem~4.3 of op. cit. that these maps are isomorphisms
when $E = \emptyset$. Our goal in this section is to extend Kato's theorem
to arbitrary $E$ and $D \ge 0$. We begin with some preliminary results.

\subsection{Duality functor under change of coefficient $\Lambda_m$}\label{sec:DCC}
  For any $1 \le i \le m$, the canonical surjection $\Lambda_m \surj \Lambda_i$
induces an exact  functor $\theta_i \colon D_i(X_{\rp}) \to D_m(X_{\rp})$.
We shall need the following result which is shown in
\cite[(4.2.4)]{Kato-duality-1} when $i =1$.
  
\begin{lem}\label{lem:RP-17}
  For $E \in D_i(X_{\rp})$, there is a natural isomorphism
  \[
  \kappa_i \colon
  \theta_i\left(\D_{X,i}(E)\right) \xrightarrow{\cong} \D_{X,m}(\theta_i(E)).
  \]
\end{lem}
\begin{proof}
We consider the free $\Lambda_m$-resolution 
$(\cdots \to \Lambda_m \xrightarrow{p^{m-i}} \Lambda_m \xrightarrow{p^{i}}
\Lambda_m \xrightarrow{\epsilon_{m,i}} \Lambda_i \to 0)$ of $\Lambda_i$,
where $\epsilon_{m,i}$ is the canonical quotient map. Letting $A^\bullet_{m,i} :=
(\cdots \to \Lambda_m \xrightarrow{p^{m-i}} \Lambda_m \xrightarrow{p^{i}}
\Lambda_m)$, this resolution yields an isomorphism
\begin{equation}\label{eqn:RP-17-0}
  \epsilon^*_{m,i} \colon
  {\bf R}\sHom_{D_m(X_{\rp})}(\Lambda_i,  W_m\Omega^{d, \rp}_{X, \log})
  \xrightarrow{\cong}
  \sHom^\bullet_{K_m(X_{\rp})}(A^\bullet_{m,i}, W_m\Omega^{d, \rp}_{X, \log}) = 
  \sE^\bullet_{m,i},
  \end{equation}
   where
$\sE^\bullet_{m,i} := (W_m\Omega^{d, \rp}_{X, \log} \xrightarrow{p^{i}} 
   W_m\Omega^{d, \rp}_{X, \log} \xrightarrow{p^{m-i}} \cdots)$, and
   $K_m(X_{\rp})$
   is the homotopy category of chain complexes of sheaves of $\Lambda_m$-modules
   on $X_{\rp}$. Since the map
   $W_i\Omega^{d, \rp}_{X, \log} \xrightarrow{\ov{p}^{m-i}} \sE^\bullet_{m,i}$
   is a quasi-isomorphism by \thmref{thm:RP-13-1}, we get that
   $\ov{p}^{m-i} \colon W_i\Omega^{d, \rp}_{X, \log}  \to W_m\Omega^{d, \rp}_{X, \log}$
   induces in $D_m(X_{\rp})$ a canonical isomorphism
\begin{equation}\label{eqn:RP-17-1}
  \lambda_{m,i} \colon W_i\Omega^{d, \rp}_{X, \log} \xrightarrow{\cong}
         {\bf R}\sHom_{D_m(X_{\rp})}(\Lambda_i,  W_m\Omega^{d, \rp}_{X, \log}).
         \end{equation}

Suppose now that $E \in D_i(X_{\rp})$ and denote $\theta_i(E)$ also by $E$.
We then get
\[
\begin{array}{lll}
  \theta_i\left(\D_{X,i}(E)\right)
          & \stackrel{\lambda_{m,i} \circ (-)}{\underset{\cong}\longrightarrow} &
        \theta_i\left({\bf R}\sHom_{D_i(X_{\rp})}(E,
   {\bf R}\sHom_{D_m(X_{\rp})}(\Lambda_i,  W_m\Omega^{d, \rp}_{X, \log}))\right) \\
            & {\cong}^{\dagger} & 
{\bf R}\sHom_{D_m(X_{\rp})}(E \stackrel{\L}{\otimes}_{\Lambda_i} {\Lambda_i},
  W_m\Omega^{d, \rp}_{X, \log}) \\
  & \cong & {\bf R}\sHom_{D_m(X_{\rp})}(E, W_m\Omega^{d, \rp}_{X, \log}) \\
    & \cong & \D_{X,m}(E),
\end{array}
\]
where the isomorphism ${\cong}^\dagger$ follows from
\cite[Prop.~12.10.12]{Yekutieli}. This finishes the proof.
\end{proof}

In view of \lemref{lem:RP-17}, we shall henceforth use the notation $E$ to denote
$\theta_i(E)$ for any $E \in D_i(X_\rp)$.
 
\begin{lem}\label{lem:RP-17-2}
  For $1 \le i \le m$ and $E \in D_i(X_{\rp})$,  the inclusion $\vartheta^d_{i} \colon
  W_i\Omega^{d, \rp}_{X, \log} \inj W_i\sH^{d,\rp}_X$ induces a commutative diagram
  \begin{equation}\label{eqn:RP-17-2-0}
    \xymatrix@C2.5pc{
      \D_{X,i}(E) \ar[r]^-{\lambda_{m,i} \circ (-)} \ar[d]_-{\vartheta^d_{i}\circ (-)} &
      \D_{X,m}(E) \ar[d]^-{\vartheta^d_{m}\circ (-)} \\
        {\bf R}\sHom_{D_i(X_\rp)}(E, W_i\sH^{d,\rp}_X)
        \ar[r]^-{V^{m-i} \circ (-)} & {\bf R}\sHom_{D_m(X_\rp)}(E, W_m\sH^{d,\rp}_X)}
  \end{equation}
  in $D_m(X_\rp)$. The same is also true if we replace
  $V^{m-i} \circ (-)$ by $\ov{p}^{m-i} \circ (-)$.
\end{lem}
\begin{proof}
By the construction of various maps in ~\eqref{eqn:RP-17-2-0}, it suffices to show
  that the diagram of complexes of sheaves
  \begin{equation}\label{eqn:RP-17-2-1}
    \xymatrix@C1.5pc{
      W_i\Omega^{d, \rp}_{X, \log} \ar[r]^-{\ov{p}^{m-i}} \ar[d]_-{\vartheta^d_{i}} &
      W_m\Omega^{d, \rp}_{X, \log} \ar[d]^-{\vartheta^d_{m}} \\
      W_i\sH^{d,\rp}_X \ar[r]^-{V^{m-i}} & W_m\sH^{d,\rp}_X}
  \end{equation}
  is commutative. However, this can be checked locally, in which case, one deduces it
  from diagram~\eqref{eqn:RP-10}, which implies that, for any 
$w \in \Ker(1-C)$, one has $\ov{p}(w) = \ov{p}C(w) = V(w)$. The proof for
  $\ov{p}^{m-i} \circ (-)$ is already clear.
\end{proof}

We now fix $q \in \N_0, \ m \in \N$. We let $D \ge 0$ and  $D' = D+E$ be in
$\Div_E(X)$. 
We let $D'_1 = \lceil{{D'}/{p^{m-1}}}\rceil$ and $D_1 = {D}/{p^{m-1}} :=
\lfloor{{D}/{p^{m-1}}}\rfloor$ so that $D'_1 = D_1+E$.
For $i \in J_m$, we let $\ov{\Psi}^{q',i}_{X|D_m}
\colon W_i\Omega^{q', \rp}_{X|D'_1, \log} \to
\D_{X,m}(W_{i}\sF^{q, \rp}_{X|D_1})$ denote the composite map
$W_i\Omega^{q', \rp}_{X|D'_1, \log} \to
\D_{X,i}(W_{i}\sF^{q, \rp}_{X|D_1}) \xrightarrow{\kappa_i}
\D_{X,m}(W_{i}\sF^{q, \rp}_{X|D_1})$, in which the first arrow is
${\Psi}^{q',i}_{X|D_1}$.
For $E \in  D^c_m(X_{\rp})$, we let
$\D^0_{X,m}(E) := {\bf R}\sHom_{D^c_m(X_{\rp})}(E, W_m\Omega^{d, \rp}_{X, \log})
\xrightarrow{\cong} {\bf R}{\sHom}_{D^c_m(X_\rp)}(E, W_m\sH^{d, \rp}_X)$.
Note that the canonical map $\D^0_{X,m}(E) \to \D_{X,m}(E)$ is an isomorphism.
For a complex of the sheaves of $\Lambda_m$-modules $F$ on $X_\rp$, we let
${\sD}_{X.m}(F) = {\sHom}^\bullet_{\Lambda_m}(F, W_m\sH^{d, \rp}_X)$.

We now look at the diagram
\begin{equation}\label{eqn:RP-18-1}
  \xymatrix@C2.4pc{
    W_1\Omega^{q', \rp}_{X|{D'_1}, \log} \ar[d]_-{\vartheta^{q'}_1}  \ar[r]^-{\ov{p}^{m-1}} &
    W_m\Omega^{q', \rp}_{X|D', \log}
    \ar[r]^-{R} \ar[d]^-{\vartheta^{q'}_m}  & W_{m-1}\Omega^{q', \rp}_{X|D', \log} 
    \ar[d]^-{\vartheta^{q'}_{m-1}}  \\
  W_{1}\sG^{q', \rp}_{X|-D'_1} \ar[r]^-{\ov{p}^{m-1}} \ar[d]_-{\alpha_1}
    & W_{m}\sG^{q', \rp}_{X|-D'}
    \ar[r]^-{R} \ar[d]^-{\alpha_m} & W_{m-1}\sG^{q', \rp}_{X|-D'} 
    \ar[d]^-{\alpha_{m-1}}  \\
    {\sD}_{X,1}(W_{1}\sF^{q, \rp}_{X|D_1}) \ar[d]_-{\ov{p}^{m-1} \circ (-)}
    \ar[r]^-{(R^{m-1})^\star} & 
\sD_{X,m}(W_{m}\sF^{q, \rp}_{X|D}) \ar[d]^-{\id}
\ar[r]^-{V^\star} & 
{\sD}_{X,m-1}(W_{m-1}\sF^{q, \rp}_{X|D}) \ar[d]^-{V \circ (-)} \\
\sD_{X,m}(W_{1}\sF^{q, \rp}_{X|D_1}) \ar[r]^-{(R^{m-1})^\star} \ar[d] &
  \sD_{X,m}(W_{m}\sF^{q, \rp}_{X|D}) \ar[r]^-{V^\star} \ar[d] &
   \sD_{X,m}(W_{m-1}\sF^{q, \rp}_{X|D}) \ar[d] \\
{\D}^0_{X,m}(W_{1}\sF^{q, \rp}_{X|D_1}) \ar[r]^-{(R^{m-1})^\star} &
    {\D}^0_{X,m}(W_{m}\sF^{q, \rp}_{X|D}) \ar[r]^-{V^\star} &
  {\D}^0_{X,m}(W_{m-1}\sF^{q, \rp}_{X|D}),} 
\end{equation}
where $\alpha_i$ is the duality map induced by ~\eqref{eqn:RP-15-4-2}.
The image of $\alpha_i$ lies in $D^c_i(X_{\rp})$ by \lemref{lem:RP-16-2}. 
The vertical arrows on the lower floor are the canonical maps.
The upper and the lower floors are clearly commutative.

\begin{lem}\label{lem:RP-18-2}
  The diagram~\eqref{eqn:RP-18-1} is commutative.
\end{lem}
\begin{proof}
  We only have to show that
  $(R^{m-1})^\star \circ (\ov{p}^{m-1} \circ (-)) \circ \alpha_1
  = \alpha_m \circ \ov{p}^{m-1}$ and
    $V^* \circ \id \circ \alpha_m = (V\circ (-)) \circ \alpha_{m-1} \circ R$. 
To show the first identity, we note that in terms of the pairing of
complexes, the map $(R^{m-1})^\star \circ ({\ov{p}^{m-1} \circ (-)}) \circ \alpha_1$
corresponds to the pairing 
$W_{m}\sF^{q, \rp}_{X|D} \times  W_{1}\sG^{q', \rp}_{X|-D'_1} \to W_m\sH^{d, \rp}_X$,
which in different bi-degrees is given as follows.
\begin{equation}\label{eqn:RP-18-2-0}
  \begin{array}{lll}
\deg (0,0): & \alpha \times \beta \mapsto \ov{p}^{m-1}(R^{m-1}(\alpha) \wedge \beta) \\
\deg (0,1): & x \times y \mapsto \ov{p}^{m-1}(-C[R^{m-1}(x) \wedge y]) \\    
\deg (1,0): & a \times b \mapsto \ov{p}^{m-1}(R^{m-1}(a) \wedge b). \\
  \end{array}
\end{equation}

On the other hand, $\alpha_m \circ \ov{p}^{m-1}$ is given by the sheaf level pairings
\begin{equation}\label{eqn:RP-18-2-1}
  \begin{array}{lll}
    \deg (0,0): & \alpha \times \beta \mapsto \alpha \wedge
    \ov{p}^{m-1}(\beta) \\
\deg (0,1): & x \times y \mapsto -C(x \wedge \ov{p}^{m-1}(y)) \\    
\deg (1,0): & a \times b \mapsto a \wedge \ov{p}^{m-1}(b). \\
  \end{array}
\end{equation}
We thus need to show $\ov{p}^{m-1}(R^{m-1}(\alpha) \wedge \beta) = \alpha \wedge
    \ov{p}^{m-1}(\beta)$ and $\ov{p}^{m-1}(-C[R^{m-1}(x) \wedge y]) = 
    -C(x \wedge \ov{p}^{m-1}(y))$. But both are easy to verify by choosing some lifts
    of $\beta$ and $y$ via $R^{m-1}$. We omit the details.

To show the second identity, we note that in terms of the pairing of
complexes, the map $(V\circ (-)) \circ \alpha_{m-1} \circ R$
corresponds to the pairing 
$W_{m-1}\sF^{q, \rp}_{X|D} \times  W_{m}\sG^{q', \rp}_{X|-D'} \to W_m\sH^{d, \rp}_X$,
which in different bi-degrees is given as follows.
\begin{equation}\label{eqn:RP-18-2-2}
  \begin{array}{lll}
 \deg (0,0): & \alpha \times \beta \mapsto V(\alpha \wedge F(\beta)) \\
\deg (0,1): & x \times \ov{y} \mapsto V(-C[x \wedge F(y)]) \\    
\deg (1,0): & a \times b \mapsto V(a \wedge F(b)). \\
  \end{array}
\end{equation}

On the other hand, $V^\star \circ \id \circ \alpha_m$ is given by the sheaf level
pairings
\begin{equation}\label{eqn:RP-18-2-3}
  \begin{array}{lll}
    \deg (0,0): & \alpha \times \beta \mapsto V(\alpha) \wedge \beta \\
\deg (0,1): & x \times y \mapsto -C(V(x) \wedge y) \\    
\deg (1,0): & a \times b \mapsto V(a) \wedge b. \\
  \end{array}
\end{equation}
We thus need to show $V(\alpha \wedge F(\beta)) = V(\alpha) \wedge \beta$ and
$V(-C[x \wedge F(y)]) = -C(V(x) \wedge y)$. But the first identity is obvious and
the second follows easily using the identity $VC = CV$ on $\Ker(F^{m-1}d)$, which holds
because every element of $\Ker(F^{m-1}d)$ lies in the image of $F$ and $VCF = CVF$.
This completes the proof.
\end{proof}

Our main duality theorem is the following.

\begin{thm}\label{thm:RP-19}
  There are canonical isomorphisms
  \[
  \Psi^{q,m}_{X|D} \colon W_m\Omega^{q', \rp}_{X|D', \log}
  \xrightarrow{\cong}
  \D^0_{X,m}(W_{m}\sF^{q, \rp}_{X|D}); \ \
  \Phi^{q,m}_{X|D} \colon W_{m}\sF^{q, \rp}_{X|D} \xrightarrow{\cong} 
  \D^0_{X,m}(W_m\Omega^{q', \rp}_{X|D', \log})
  \]
  in $D^c_m(X_\rp)$.
\end{thm}
\begin{proof}
  We first recall that $\Psi^{q,m}_{X|D}$ is the composition of the middle vertical
  arrows in ~\eqref{eqn:RP-18-1} and
  $\Phi^{q,m}_{X|D}$ is the composition of the canonical maps
  $W_{m}\sF^{q, \rp}_{X|D} \to  \D^0_{X,m}({\D}^0_{X,m}(W_{m}\sF^{q, \rp}_{X|D})) \to
 {\D}^0_{X,m}(W_m\Omega^{q', \rp}_{X|D', \log})$,
  where the first arrow is the evaluation map and the second arrow is the one
  induced by ~\eqref{eqn:RP-15-4-3}. By \cite[Thm.~4.3]{Kato-duality-1},
  showing that $\Phi^{q,m}_{X|D}$ is an isomorphism is equivalent to  showing
  that the same is true for $\Psi^{q,m}_{X|D}$.
  Since $W_m\Omega^{q', \rp}_{X|D', \log}$ is a sheaf, our task is 
  equivalent to showing that $\Psi^{q,m}_{X|D}$ induces an isomorphism of sheaves
\begin{equation}\label{eqn:RP-19-0}
\Psi^{q,m}_{X|D} \colon W_m\Omega^{q', \rp}_{X|D', \log} \xrightarrow{\cong}
    {\bf R}^0{\sHom}_{D_m(X_\rp)}(W_{m}\sF^{q, \rp}_{X|D}, W_m\sH^{d,\rp}_X), \ \ \mbox{and}
\end{equation}
\begin{equation}\label{eqn:RP-19-1}
  {\bf R}^i{\sHom}_{D_m(X_\rp)}(W_{m}\sF^{q, \rp}_{X|D}, W_m\sH^{d,\rp}_X) = 0 \ \
  \forall \ \ i \neq 0.
  \end{equation}
  In the remainder of the proof, we shall show this by induction on $m$.
  Let us first assume that $m =1$.

By the definition of the pairing ~\eqref{eqn:RP-15-4-2},
  we have a commutative diagram of distinguished triangles
\begin{equation}\label{eqn:RP-19-2}
    \xymatrix@C1.5pc{
      \frac{\Fil_{-D'}\Omega^{q',\rp}_{X}}{B_1\Fil_{-D'}\Omega^{q'-1,\rp}_{X}}[-1] \ar[r]
      \ar[d]_{\beta[-1]} & W_1\sG^{q', \rp}_{X|-D'} \ar[r] \ar[d]^-{\Psi^{q,1}_{X|D}}
      \ar[r] &
      \Fil_{-D'}\Omega^{q', \rp}_{X} \ar[r]^-{1-\ov{F}} \ar[d]^-{\alpha} &
      \frac{\Fil_{-D'}\Omega^{q',\rp}_{X}}{B_1\Fil_{-D'}\Omega^{q'-1,\rp}_{X}}
      \ar[d]^{\beta}
      \\
      \D_{X,1}(Z_1\sA) \ar[r] &
      \D_{X,1}(W_{1}\sF^{q, \rp}_{X|D}) \ar[r] &  \D_{X,1}(\sA)[1]
      \ar[r]^-{(1-C)^\star} & \D_{X,1}(Z_1\sA)[1]}
  \end{equation}
in $D_1(X_\rp)$ (cf. \cite[Lem.~11.13]{Ren-Rulling}), where
  $\alpha$ (resp. $\beta$) is induced by the pairing
  $\<,\>^{0,0}_0$ (resp. $\<,\>^{0,1}_1$) in ~\eqref{eqn:RP-15-1} (resp.
  ~\eqref{eqn:RP-15-3}). In the above, we have written $\sA$ (resp. for $Z_1\sA$)
for $\Fil_D\Omega^{q,\rp}_{X}$ (resp. $Z_1\Fil_D\Omega^{q,\rp}_{X}$)
as a shorthand. Using ~\eqref{eqn:RP-19-0}, it suffices to show that $\alpha$ and
  $\beta$ are isomorphisms.

To prove that $\alpha$ is an isomorphism, we note that it factors through the
  maps
  \[
  \Fil_{-D'}\Omega^{q', \rp}_{X} \xrightarrow{\cong}
(\Fil_{-D'}\Omega^{q'}_{X})^{\rp} \xrightarrow{\cong}
\left(\sHom_{\sO_{X_\et}}(\Fil_D\Omega^{q}_{X}, \Omega^d_X)\right)^{\rp}
\xrightarrow{\cong} \D_{X,1}(\Fil_D\Omega^{q,\rp}_{X}[-1]).
\]
The first arrow is clearly an isomorphism, and the second arrow is
an isomorphism by \lemref{lem:Can-filt-5} and the proof of
\cite[Thm.~4.1.4]{JSZ} (see its Step~2). The third arrow is an isomorphism
by \cite[Thm.~3.2(ii)]{Kato-duality-1}. It follows that $\alpha$ is an
isomorphism. An identical argument shows that $\beta$ is also an isomorphism. 
By compositing with the inclusion (which is a quasi-isomorphism)
$\Omega^{q'\rp}_{X|D', \log} \inj W_1\sG^{q', \rp}_{X|-D'}$, we get that
$\Psi^{q,1}_{X|D} \colon \Omega^{q'\rp}_{X|D', \log} \to \D_{X,1}(W_{1}\sF^{q, \rp}_{X|D})$
is an isomorphism. But this is equivalent to proving ~\eqref{eqn:RP-19-0}
and ~\eqref{eqn:RP-19-1} for $m =1$.

We now prove the general case. The vanishing claim ~\eqref{eqn:RP-19-1}
follows by induction, using the long exact cohomology
sequence associated to (the Kato dual of) the first distinguished
triangle of \thmref{thm:RP-16-3}.
To prove ~\eqref{eqn:RP-19-0} for $m \ge 2$,
we apply \lemref{lem:RP-18-2} to get a commutative diagram of sheaves
\begin{equation}\label{eqn:RP-19-3}
  \xymatrix@C1.5pc{
    0 \ar[r] & \Omega^{q', \rp}_{X|D'_1, \log}
      \ar[d] \ar[r]^-{\ov{p}^{m-1}} &
    W_m\Omega^{q', \rp}_{X|D', \log}
    \ar@{->>}[r]^-{R} \ar[d] & W_{m-1}\Omega^{q', \rp}_{X|D', \log} 
    \ar[d] \\
    0 \ar[r] & H^0({\D}_{X,1}(W_{1}\sF^{q, \rp}_{X|D_1}))
  \ar[r]^-{(R^{m-1})^\star} &
    H^0({\D}_{X,m}(W_{m}\sF^{q, \rp}_{X|D})) \ar@{->>}[r]^-{V^\star} &
    H^0({\D}_{X,m-1}(W_{m-1}\sF^{q, \rp}_{X|D})).}
\end{equation}

The two rows are exact by \thmref{thm:RP-16-3}, \lemref{lem:RP-17},
and ~\eqref{eqn:RP-19-1}, shown above.
The left and the right vertical arrows are isomorphisms by induction. It
follows that the middle vertical arrow is also an isomorphism. This completes the
proof of the theorem.
\end{proof}

\section{Duality for logarithmic Hodge--Witt cohomology}\label{sec:DTHW}
In this section, we prove Theorems~\ref{thm:Main-2} and ~\ref{thm:Main-8}.
We begin with the extension of Milne's duality to modulus pairs.

\subsection{The duality theorem over finite fields}\label{sec:DT-0}
In the setting of \S~\ref{sec:KD*}, we assume additionally in this subsection that
$X$ is a smooth, proper and connected scheme of
dimension $d$ over a finite field $k$. We assume that $D \ge 0, \ D_\red = E$ and
$D' = D+E$. We let $q \in \N_0, m \in \N$ and $q' = d-q$.

By combining \thmref{thm:RP-19} with \cite[Thm.~5.1]{Kato-duality-1} (applied to
the structure map $f \colon X \to \Spec(k)$),  the natural isomorphism
${\bf R}\Gamma(X_\rp, E) \xrightarrow{\cong} {\bf R}\Gamma(X_\et, E|_{X_\et})$ (for
$E \in D_m(X_\rp)$), and the trace map
$\Tr_m \colon H^{d+1}_\et(X, W_{m}\Omega^d_{X,\log}) \xrightarrow{\cong} \Lambda_m$
(cf. \cite[Cor.~1.12]{Milne-Zeta}), we obtain a cohomology pairing
\begin{equation}\label{eqn:Pair-4-2}
 \H^i_\et(X, W_{m}\sF^{q, \bullet}_{X|D}) \times
    H^{d+1-i}_\et(X, W_{m}\Omega^{q'}_{X|D', \log})  \xrightarrow{\cup}
    \Lambda_m,
    \end{equation}
which is compatible with Milne's pairing for the cohomology of the sheaves
$W_m\Omega^{\bullet}_{X,\log}$ (cf. \cite[Thm.~1.14]{Milne-Zeta}).

The following application of \thmref{thm:RP-19} extends Milne's
duality to modulus pairs.

\begin{thm}\label{thm:Perfect-finite}
For $i \in \N_0$,  \eqref{eqn:Pair-4-2} is a perfect pairing of finite abelian groups.
  \end{thm}
\begin{proof}
  Using \thmref{thm:RP-13-1} and \corref{cor:RP-16}, the proof is readily reduced to
  the $m=1$ case which is already shown in the proof (see its Step~ 2) of
  \cite[Thm.~4.1.4]{JSZ}. Recall here that
  $\Fil_D\Omega^q_X \cong \Omega^q_X(\log E)(D)$ (cf. \lemref{lem:Log-fil-4}(8)).
  \end{proof}
 
\vskip.2cm

The following is \corref{cor:Main-4}.

    \begin{cor}\label{cor:Duality-log}
  For $X$ as in \thmref{thm:Perfect-finite}, there is a perfect pairing of
  finite abelian groups
  \[
  H^i_\et(X, W_m\Omega^q_{X_\log, \log}) \times H^{d+1-i}_\et(X, W_m\Omega^{d-q}_{X|E, \log})
  \xrightarrow{\cup} \Lambda_m.\]
  \end{cor}
\begin{proof}
  Combine \propref{prop:RP-4}, \thmref{thm:Perfect-finite} and
  \cite[Lem.~10.2]{KM-3}(4).
\end{proof}

\subsection{Duality theorem over DVR}\label{sec:DT-dvr}
We prove \thmref{thm:Main-8} in this subsection. 
We let $R$ be an equicharacteristic complete discrete valuation ring with
quotient field $K$ and finite residue field $\ff$.
We let $S = \Spec(R)$ and let $f \colon X \to S$ be a projective
strictly semistable connected
scheme of relative dimension $d$ in the sense of \cite[Defn.~1.4.2]{Zhao}.
In particular, $X$ has $p$-rank $d+1$.
We let $X_\eta$ (resp. $X_s := X \times_S \Spec(\ff)$) denote the generic
(resp. special) fiber of $f$. Our assumptions imply that $X_\eta$ is smooth over $K$
and $X_s$ is reduced and an snc divisor on $X$.
We let $E \subset X$ be an snc divisor
with irreducible components $E_1, \ldots , E_r$.  We let $j \colon
X_\eta \inj X$ and $\iota \colon X_s \inj X$ denote the inclusions.
We let $D \in \Div^{\eff}_E(X)$ such that $D_\red = E$ and let $D' = D+E$.

We have the canonical maps 
\begin{equation}\label{eqn:Pair-5-dvr}
  \begin{array}{lll}
 {\bf R}\iota_*(\iota^* (W_{m}\sF^{q, \bullet}_{X|D}) \otimes^{\mathbb{L}}
{\bf R}\iota^{!}(W_{m}\Omega^{d+1-q}_{X|D',\log})) & \to &
W_{m}\sF^{q, \bullet}_{X|D} \otimes^{\mathbb{L}} 
{\bf R}\iota_*{\bf R}\iota^{!}(W_{m}\Omega^{d+1-q}_{X|D',\log}) \\
& \to &W_{m}\sF^{q, \bullet}_{X|D}  \otimes^{\mathbb{L}} 
W_{m}\Omega^{d+1-q}_{X|D',\log} \\
&\to & W_{m}\sH^{d+1, \bullet}_X, \\
\end{array}
\end{equation}
where the first arrow is an isomorphism by the projection formula, the second arrow is induced by the counit of adjunction ${\bf R}\iota_*{\bf R}\iota^{!} \to \id$, and the third arrow is obtained by
using the pairing in ~\eqref{eqn:RP-15-4-0}.

By comparing the diagram ~\eqref{eqn:RP-15-4-1} with the
pairing of ~\eqref{eqn:Pair-5-dvr}, we get a commutative diagram of pairings
\begin{equation}\label{eqn:Pair-3-1-dvr}
  \xymatrix@C1.5pc{  
    W_{m}\sF^{q, \bullet}_{X|D} \times {\bf R}\iota^{!}(W_{m}\Omega^{d+1-q}_{X|D',\log})
    \ar[r]^-{\< \ , \ \>} \ar@<5ex>[d] &
  {\bf R}\iota^{!}(W_m\sH^{d+1,\bullet}_X) \ar@{=}[d] \\
  W_{m} \sF^{q, \bullet}_X \times {\bf R}\iota^{!}(W_m\Omega^{d+1-q}_{X,\log})
  \ar[r]^-{\< \ , \ \>}
  \ar@<9ex>[u] & {\bf R}\iota^{!}(W_{m}\sH^{d+1,\bullet}_X)}
  \end{equation}
on $X_\et$ in which the vertical arrows are the canonical maps.

Considering the induced hypercohomology pairings and composing them
with the trace map ${\rm Tr} \colon
 H^{d+2}_{\et, X_s}(X, W_{m}\Omega^{d+1}_{X,\log}) \xrightarrow{\cong} \Lambda_m$
 of \cite[Cor.~1.4.10]{Zhao}, we get a cohomology pairing
 \begin{equation}\label{eqn:Pair-4-2-dvr}
 \H^i_\et(X, W_{m}\sF^{q, \bullet}_{X|D}) \times
    H^{d+2-i}_{\et, X_s}(X, W_{m}\Omega^{d+1-q}_{X|D', \log})  \xrightarrow{\cup}
    \Lambda_m
    \end{equation}
which is compatible with the unramified duality pairing of
\cite[(3.1)]{Zhao}. Furthermore, as $D \ge 0$ varies, this gives a strict
pairing between pro and ind abelian groups.

\vskip.2cm

The main result of this subsection is the following duality theorem for
semistable schemes with modulus. This provides a refinement of the
main result of \cite{Zhao}.

\begin{thm}\label{thm:Perfect-dvr}
For $i \in \N_0$,  ~\eqref{eqn:Pair-4-2-dvr} induces an isomorphism of
  $\Lambda_m$-modules
  \begin{equation}\label{eqn:Perfect-dvr-1}
    \H^i_\et(X, W_{m}\sF^{q, \bullet}_{X|D}) \xrightarrow{\cong}
    \Hom_{\Lambda_m}(H^{d+2-i}_{\et, X_s}(X, W_{m}\Omega^{d+1-q}_{X|D', \log}), \Lambda_m).
  \end{equation}
  \end{thm}
\begin{proof}
  By iterating the reduction step in the proof of \thmref{thm:Perfect-finite}, we
  reduce the proof to the case $m=1$, which was already established in
  \cite[Prop.~3.4.5]{Zhao}. 
\end{proof}

\vskip.2cm

If we endow $H^{d+2-i}_{\et, X_s}(X, W_{m}\Omega^{d+1-q}_{X|D', \log})$ with the discrete
topology and $\H^i_\et(X, W_{m}\sF^{q, \bullet}_{X|D})$ with the compact-open topology
using \thmref{thm:Perfect-dvr}, then it follows that
~\eqref{eqn:Pair-4-2-dvr} is a perfect pairing of topological abelian groups
so that
\[
H^{i}_{\et, X_s}(X, W_{m}\Omega^{d+1-q}_{X|D', \log})
\xrightarrow{\cong} \Hom_{\cont}(\H^{d+2-i}_\et(X, W_{m}\sF^{q, \bullet}_{X|D}), \Lambda_m).
\]
Since the above isomorphisms are compatible as $D \ge 0$ varies, passing to
the limit over $D$, one recovers the main result of
\cite{Zhao}.

\end{document}